# Geometric Zeta Functions, $L^2$-Theory, and Compact Shimura Manifolds

Anton Deitmar

# Contents

# Introduction

Zeta functions encoding geometric information such as zeta functions of algebraic varieties over finite fields or zeta functions of finite graphs will loosely be called geometric zeta functions in the sequel. Sometimes the geometric situation gives one tools at hand to prove analytical continuation, functional equation and an adapted form of the Riemann hypothesis. This works as follows: The geometrical data is related to the fixed point set of an operator, this fixed point set is considered as the local data. There is a Lefschetz formula relating these to the action of this operator on a global object such as a suitable cohomology theory, thus expressing the zeta function by determinants of the operator. This is what we will call a determinant formula.

This could have been the guiding idea when Selberg in 1956 [Sel] first defined a geometric zeta function for compact Riemannian surfaces by now known as the Selberg zeta function. To play the rôle of a Lefschetz formula Selberg designed his by now famous trace formula which soon was generalized to higher dimensional Shimura manifolds, i.e. quotients of symmetric spaces. The generalization of the zeta function took somewhat more time.

In 1977 R. Gangolli [Gang] defined a zeta functions for all compact rank one Shimura manifolds. In 1985 M. Wakayama defined twisted versions of these [Wak].

Up to this point all the work was formulated in terms of the harmonic analysis not obeying the above philosophy and there was no such a thing as a determinant formula which would have given a more comprehensive way to formulate results and deeper insights into the connection between zeta functions and the topology of Shimura manifolds. Since the generalized Selberg zeta functions have infinitely many zeroes it was clear that a determinant formula would require a notion of a determinant of infinite rank operators. The latter was provided by the Ray-Singer regularized determinant [RS-AT]. In 1986 D. Fried [Fr] proved a determinant formula for the real hyperbolic spaces which are special symmetric spaces of rank one. Fried further showed that certain products of the generalized Selberg functions in these cases give the Ruelle zeta function of the geodesic flow. D. Ruelle had defined the zeta function carrying his name ten years earlier in a more general context [Rue] where he was able to show the existence of a meromorphic continuation but not very much more than that. Fried's work however gave a determinant formula and in that way a functional equation and the generalized Riemann hypothesis. Fried showed that the special value at the center of the functional equation is given by a topological invariant: the analytical torsion.

The extension to higher rank seemed impenetrable until in 1991 H. Moscovici and R. Stanton [MS-2] found a way to circumvent the analysis of orbital integrals in that they used higher supersymmetry arguments to provoke a cancellation in the terms of higher rank. Their approach produced a kind of a Ruelle zeta function which still had analytic torsion as a special value. Unfortunately this only worked for the symmetric spaces attached to $SL_3(\mathbb{R})$ or $SO(p,q)$ with $pq$ odd.

In [D-Hitors] the author was able to extend their results to all Shimura manifolds of odd dimension, defining higher twisted versions of the analytic torsion which showed up as special values there. The present work will extend this to all Shimura manifolds and also give various other types of zeta functions in the higher rank situation. Apparently in even dimensions the invariants to expect are not the analytic (resp. holomorphic) torsion numbers but quotients of those and their analogues in the $L^2$-theory of Atiyah [At]. The term $L^2$-theory here means that in constructing the invariants the usual trace is replaced by a certain "canonical trace". The invariants thus arising are called $L^2$-invariants.

Further sources for Selberg and Ruelle zeta functions are [BuOl], [Fr2], [Fr3], [Fr4], [Ju], [MS-1], [Wak2], [Wak3], [Wak4]. For related other zeta functions read [Hash1], [Hash2], [Hash3].

The first Chapter contains an introduction to the theory of Ray-Singer determinants, their $L^2$-analogues and convergence theorems relating the one to the other. In general it is an open question whether Ray-Singer determinants in a tower of coverings do converge to their $L^2$ analogues, but this can be shown for determinant functions which also might be called characteristic functions.

The second chapter introduces the reader to the classical, twisted and $L^2$-torsion. In the third chapter one finds some folklore concerning Shimura manifolds.

Chapter 4 contains a treatment of the two kinds of zeta functions in the rank one case. They are expressed in a "complete" determinant formula which also contains the "factor at infinity".

In chapter 5 one finds a definition of a zeta function for holomorphic Shimura manifolds that has the holomorphic torsion as special value. Of particular neatness is the case where the symmetric space is a power $\mathcal{H}^m$, $m \in \mathbb{N}$ of the hyperbolic plane $\mathcal{H}$, in which case the Euler product is strikingly simple (sec. 5.2).

Holomorphic Shimura manifolds come from symmetric spaces of fundamental rank zero. Chapter 6 deals with spaces of higher fundamental rank. In section 6.4 we give the appropriate Ruelle zeta function and in 6.5 we find a generalization of the Selberg function to this setting.

In chapter 7 we present a different approach to the geometric zeta functions in the rank one case which gives rise to some interesting speculations about possible generalizations.

In the 8[th] chapter we give the analytic continuation of theta functions which are somewhat dual objects to the zeta functions in the following way: The geometrical zeta functions considered here are defined by geometrical data (closed geodesics) and appear to have their singularities (zeroes and poles) in spectral points (Laplace eigenvalues). The theta functions in turn are defined by the very same spectral data and happen to have their singularities in the geometrical data. The duality carries on in the proof that we give in that the functional equation of the theta function is not gotten from the functional equation of the zeta function but from its Euler product.

All Shimura manifolds considered here are compact. There is, however also a Selberg zeta function in the noncompact case. The extension of the contents of this paper to the noncompact case will be the subject of future work.

*I thank U. Bunke, A. Juhl, M.Olbrich and K. Takase for valuable comments. This paper was written during my stay at the MSRI at Berkeley in the special year on automorphic forms 1994/95. I thank the organizers for inviting me and all people there for generating a concentrated*

*working atmosphere. Further I thank the staff for their continuous support and for giving me an office with a terrific view to the San Francisco Bay.*

*I finally thank the Deutsche Forschungsgemeinschaft for financing this research.*

# Notation

We will write $\mathbb{N}, \mathbb{Z}, \mathbb{Q}, \mathbb{R}$ and $\mathbb{C}$ for the natural, integer, rational, real and complex numbers.

By a **virtual vector space** we understand a formal difference $V = V_+ - V_-$ of vector spaces. If the components are finite dimensional, we define the **dimension** of $V$ to be $\dim(V) = \dim(V_+) - \dim(V_-)$. An **endomorphism** $A$ **of virtual a virtual vector space** $V$ is a formal difference $A = A^+ - A^-$ where $A^{\pm}$ is an endomorphism of $V_{\pm}$. We then define the **trace** and the **determinant** of $A$ as $\operatorname{tr}(A) = \operatorname{tr}(A^+) - \operatorname{tr}(A^-)$ and $\det(A) = \det(A^+)/\det(A^-)$.

When $V = \oplus_{k \in \mathbb{Z}} V_k$ is a $\mathbb{Z}$-graded vector space we will automatically consider it as a virtual vector space

$$V_* = V_{\text{even}} - V_{\text{odd}} = \oplus_{k \text{even}} V_k - \oplus_{k \text{odd}} V_k$$

and if both are finite dimensional we define

$$\dim(V_*) = \dim(V_{\text{even}}) - \dim(V_{\text{odd}}).$$

For any Hilbert $H$ space we denote by $B(H)$ the algebra of bounded linear operators on $H$. When we speak of an operator on a Hilbert H space this will usually mean a densely defined possibly unbounded linear operator on H.

For a smooth vector bundle $E$ over some smooth manifold we denote by $C^{\infty}(E)$ the space of smooth sections of $E$ and if E is defined over a Riemannian manifold and equipped with a metric we denote by $L^2(E)$ the Hilbert space of measurable sqare integrable sections of E modulo the space of measurable sections vanishing outside a set of measure zero.

Let $\mathcal{A}$ be an abelian category. A **complex** over $\mathcal{A}$ is a sequence $\ldots \to E_{-1} \to E_0 \to E_1 \to \ldots$, where the $E_j$ are objects of $\mathcal{A}$ and the **differentials** $d_j : E_j \to E_{j+1}$ satisfy $d_j d_{j+1} = 0$. A complex is said to be **positive**, if $E_{-j} = 0$ for $j \in \mathbb{N}$ and **bounded from below**, if there is a $j_0 \in \mathbb{N}$ such that $E_{-j} = 0$ for $j \geq j_0$. A complex $E$ is said to be **finite**, if there are $a, b \in \mathbb{Z}$, $a \leq b$ such that $E_j = 0$ for $j < a$ or $j > b$. The number $b - a + 1$ is then called the **length** of the complex $E$.

Let $\mathcal{A}$ denote an abelian category. The **Grothendieck group** $\mathcal{G}$ of $\mathcal{A}$ is the abelian group generated by the set of isomorphism classes of objects in $\mathcal{A}$, subject to the relations $[A] - [B] -$ $[C] = 0$ for every exact sequence $0 \to B \to A \to C \to 0$ in $\mathcal{A}$. Every element $x$ in $\mathcal{G}$ can be written as $x = x_+ - x_-$, where $x_\pm \in \mathcal{A}$.

Lie algebras of real Lie groups $G, H, K, \ldots$ are denoted by the corresponding small case german letters with index zero: $\mathfrak{g}_0, \mathfrak{h}_0, \mathfrak{k}_0$ and their complexifications are denoted without any index, so $\mathfrak{g} = \mathfrak{g}_0 \otimes \mathbb{C}$.

For any locally compact group $G$ we write $\hat{G}$ for the **unitary dual**, i.e.: $\hat{G}$ is the set of equivalence classes of irreducible unitary representations of $G$.

# Chapter 1

# Determinants, $L^2$-Determinants

In this chapter we give an introduction to the classical theory of Ray-Singer determinants and the not so classical $L^2$-determinants. We give a comparison of the Ray-Singer determinant to the seemingly more direct Fredholm determinant which a fortiori shows that they amount to the same. We then consider the $L^2$-theoretical counterpart in the case of positive Gromov-Novikov-Shubin invariants. We finally show that in a tower of coverings the determinant functions converge to the corresponding $L^2$-determinant functions.

## 1.1  Determinants

### 1.1.1  Ray-Singer Determinants

To extend the notion of a determinant to infinite dimensions we follow the regularizing procedure by Ray and Singer [RS-RM]. So let $H$ be a Hilbert space and consider a positive operator $A$, densely defined on $H$. The operator $A$ is called **admissible**, if there is a natural number n such that $A^{-n}$ is of trace class, and the **zeta function of** $A$:

$$\zeta_A(s) \stackrel{\text{def}}{=} \mathrm{tr}(A^{-s}), \quad \mathrm{Re}(s) > n,$$

admits an analytic continuation to some function holomorphic at $s = 0$. For an admissible operator $A$ we define its **determinant** by

$$\det(A) \stackrel{\text{def}}{=} \exp(-\frac{d}{ds}\mid_{s=0} \zeta_A(s)).$$

**Example 1.1.1** *On $l^2(\mathbb{N})$ with the canonical basis $(e_j)_{j \in \mathbb{N}}$ define an operator $A$ by $A(e_j) = je_j$. Then $A^{-2}$ is if trace class and $\zeta_A(s) = \zeta(s)$, the Riemann zeta function. We have $\zeta'(0) = -\frac{1}{2}\log(2\pi)$ and thus*

$$\det(A) = \sqrt{2\pi}.$$

Since $\overline{\zeta_A(\bar{s})} = \zeta_A(s)$ the determinant is real and positive and there are some easy rules which we list as follows:

- $\det(A \oplus B) = \det(A)\det(B)$ if both are admissible.

- $\det(cA) = c^{\zeta_A(0)}\det(A)$ for $c > 0$.

- $\det(A^s) = \det(A)^s$ for $s > 0$.

If the operator $A$ has a nontrivial kernel we will write $\det(A) = 0$, but if $A' = A \mid_{(ker(A))^\perp}$ is admissible we will call $A$ admissible as well and **define**:

$$\det'(A) \stackrel{\text{def}}{=} \det(A').$$

**Proposition 1.1.1** *A positive elliptic differential operator $D$ defined on a hermitian vector bundle over a compact smooth Riemannian manifold is admissible and the determinant is independent of the choice of the metrics.*

**Proof:** See ([Gilk], Lemma 1.10.1) for the admissibility, for the independence of metrics recall that $\det(D)$ only depends on the eigenvalues of $D$. ∎

Since the determinant of a finite dimensional operator is the product of the eigenvalues, the regularized determinant may give the definition of a regularized product . So let $(a_j)_{j \in \mathbb{N}}$ be a sequence of real numbers tending to $+\infty$. The sequence $(a_j)_{j \in \mathbb{N}}$ is called **admissible**, if the operator $D$, which has the $a_j$ as eigenvalues, is admissible. Under these circumstances we define the **regularized product** as

$$\prod_{j \in \mathbb{N}} a_j \stackrel{\text{def}}{=} det(D).$$

The above example thus gives

$$\prod_{n \in \mathbb{N}} n = \sqrt{2\pi}.$$

### 1.1.2 Fredholm Determinants

There is also a more direct way to define determinants of operators on infinite dimensional spaces, namely the Fredholm determinant, which is defined as follows: Let T be a trace class operator, then we define

$$\det_{Fr}(1 + T) \stackrel{\text{def}}{=} \sum_{k=0}^{\infty} \text{tr} \wedge^k T.$$

The sum is absolutely convergent. Moreover, assume T normal with eigenvalues $(\lambda_n)_{n \in \mathbb{N}}$ then it follows that the infinite product

$$\prod_{n \in \mathbb{N}} (1 + \lambda_n)$$

is convergent and equals the Fredholm determinant of $1 + T$. If the operator norm $\| T \|$ of $T$ is less than 1 we also have

$$\det_{Fr}(1 + T) = \exp(-\sum_{n=1}^{\infty} \frac{(-1)^n}{n} \mathrm{tr} T^n).$$

The connection between Fredholm determinant and zeta regularized determinant is:

**Proposition 1.1.2** *Let $A$ and $A + 1$ be admissible with $\ker(A) = 0$. Assume $A^{\epsilon-1}$ is of trace class for some $\epsilon > 0$, then*

$$\det_{Fr}(1 + A^{-1}) = \frac{\det(A + 1)}{\det(A)}.$$

**Proof:** Since the claim is certainly true for finite dimensions we can divide both sides by the factors stemming from the finitely many eigenvalues of A which are less than 1 in value. This means, we may assume $\| A^{-1} \| < 1$ so that we have at hand the exponential series description of the Fredholm determinant. Let $(\lambda_n)$ denote the eigenvalues of $T$ then we have:

$$\begin{aligned}
\frac{\partial}{\partial s} \zeta_{A+1}(s) &= -\sum_{n=1}^{\infty} \log(\lambda_n + 1)(\lambda_n + 1)^{-s} \\
&= -\sum_{n=1}^{\infty} \log(\lambda_n)(\lambda_n + 1)^{-s} - \sum_{n=1}^{\infty} \log(1 + \frac{1}{\lambda_n})(\lambda_n + 1)^{-s}.
\end{aligned}$$

The second term converges to $-\log \det_{Fr}(1 + A^{-1})$ as $s \to 0$. The first equals

$$-\sum_{n=1}^{\infty} \log(\lambda_n) \lambda_n^{-s} - \sum_{n=1}^{\infty} log(\lambda_n) \lambda_n^{-s}((1 + \frac{1}{\lambda_n})^{-s} - 1),$$

in which the first summand is $\frac{\partial}{\partial s} \zeta_A(s)$. Let $f_s(x) = (1+x)^{-s} - 1$ then the Taylor series expansion of $f_s$ around $x = 0$ is

$$f_s(x) = \sum_{n=1}^{\infty} \frac{(-1)^n}{n!} s(s+1) \ldots (s+n-1) x^n,$$

so $f_s(x) = sxh(s, x)$ where $h$ is differentiable around $(x, s) = (0, 0)$, so we have $|f_s(x)| \leq c_1 sx$ for some $c > 0$ thus

$$\begin{aligned}
\left| \sum_{n=1}^{\infty} \log(\lambda_n) \lambda_n^{-s}((1 + \frac{1}{\lambda_n})^{-s} - 1) \right| &\leq c_2 s \sum_{n=1}^{\infty} \log(\lambda_n) \lambda_n^{-s-1} \\
&\leq c_3 s \sum_{n=1}^{\infty} \lambda_n^{\epsilon-s-1}
\end{aligned}$$

so this contribution vanishes at $s = 0$. ∎

### 1.1.3 Characteristic Functions

We will from now on restrict our attention to generalized Laplace operators. Much of what we are going to say will hold in a more general context but our restriction gives us some serious technical simplifications.

Let $E$ denote a smooth vector bundle over a smooth Riemannian manifold $M$. A **generalized Laplacian** on $E$ is a second order differential operator $D$ such that for the principal symbol of $D$ we have

$$\sigma_2(x, \xi) = \| \xi \|^2 Id.$$

Furthermore we insist that $D$ is semipositive, i.e. $D \geq 0$.

It then follows that there is a connection $\nabla^E$ on $E$ such that $D$ equals the connection Laplacian up to an operator of order zero ([BGV], Prop. 2.5). It further follows that the heat operator $e^{-tD}$, $t > 0$ is of trace class on $L^2(E)$ and that

$$\mathrm{tr}(e^{-tD}) \sim \sum_{k=0}^{\infty} c_{\alpha_k} t^{\alpha_k},$$

as $t \to \infty$, where $\alpha_k \in \mathbb{R}$, $\alpha_k \to \infty$.

**Proposition 1.1.3** *Let $D$ denote a generalized Laplacian over a compact smooth Riemannian manifold $M$. Then the **characteristic function**: $\lambda \mapsto \det(D + \lambda)$, $\lambda > 0$ of $D$ extends to an entire function with zeroes at the eigenvalues of $-D$, the order of a zero being the multiplicity of the eigenvalue. For $\lambda \to \infty$ we have the following asymptotics:*

$$-logdet(A + \lambda) \sim \sum_{\alpha_k \neq 0, -1, -2, \ldots} c_{\alpha_k} \Gamma(\alpha_k) \lambda^{-\alpha_k}$$
$$- \sum_{k=0, -1, -2, \ldots} c_k \left( log\lambda - \sum_{j=1}^{k} \frac{1}{j} \right) \frac{(-\lambda)^k}{k!},$$

**Proof:** Let $M(s, \lambda) = \Gamma(s)\zeta_{A+\lambda}(s)$. We have the differential equation:

$$\frac{\partial}{\partial \lambda} M(s, \lambda) = -M(s+1, \lambda).$$

Now let $m \in \mathbb{N}$, we get

$$(\frac{\partial}{\partial \lambda})^{m+1} \zeta_{A+\lambda}(s) = (-1)^{m+1}(s+m)(s+m-1) \ldots s\zeta_{A+\lambda}(s+m+1),$$

so that for m large enough it follows

$$(\frac{\partial}{\partial \lambda})^{m+1} \zeta_{A+\lambda}(0) = 0.$$

It is clear that we have log det$(D + \lambda) = \lim_{s \to 0}(M(s, \lambda) - \zeta_{A+\lambda}(0))/s$ and the limit in $s$ may be interchanged with the $\lambda$-derivation. It follows that for $m$ large enough we have

$$(\frac{\partial}{\partial \lambda})^{m+1}\log \det(D + \lambda) = (-1)^{m+1}M(m+1, \lambda)$$
$$= (-1)^{m+1}\Gamma(m+1)\sum_{n=1}^{\infty} \frac{1}{(\lambda + \lambda_m)^{m+1}}.$$

From this the claim follows. For the asymptotics see [Voros]. ∎

## 1.2   $L^2$-determinants

### 1.2.1   Von Neumann Dimension and Trace

From now on we will fix the following situation: $X_\Gamma$ will denote a smooth compact oriented Riemannian manifold with universal covering $X$ and fundamental group $\Gamma$. We will assume $\Gamma$ infinite. Over $X_\Gamma$ we will have a smooth hermitian vector bundle $E_\Gamma$ with pullback $E$ over $X$. On $E_\Gamma$ we fix a generalized Laplacian $D_\Gamma$ and we denote its pullback to $E$ by $D$.

Let $L^2(\Gamma)$ denote the $L^2$-space over $\Gamma$ with respect to the counting measure. The group $\Gamma$ acts on $L^2(\Gamma)$ by left translations. Inside the algebra $B(L^2(\Gamma))$ of all bounded linear operators on the Hilbert space $L^2(\Gamma)$ we consider the **von Neumann algebra of** $\Gamma$, denoted $VN(\Gamma)$ generated by the left translations $(L_\gamma)_{\gamma \in \Gamma}$. There is a canonical trace on $VN(\Gamma)$ defined by tr$(\sum_\gamma c_\gamma L_\gamma) = c_e$ which makes $VN(\Gamma)$ a type II$_1$ von Neumann algebra [GHJ], which is a factor if and only if every nontrivial conjugacy class in $\Gamma$ is infinite.

It is easy to see that the commutant of $VN(\Gamma)$ is the von Neumann algebra generated by the right translations.

Fix a fundamental domain $\mathcal{F}$ of the $\Gamma$-action on $X$. We get isomorphisms of unitary $\Gamma$-modules:

$$L^2(E) \cong L^2(\Gamma)\hat{\otimes}L^2(E|_\mathcal{F}) \cong L^2(\Gamma)\hat{\otimes}L^2(E_\Gamma).$$

So for the von Neumann algebra $B(L^2(E))^\Gamma$ of operators commuting with the $\Gamma$-action we get

$$B(L^2(E))^\Gamma \cong VN(\Gamma)\hat{\otimes}_w B(L^2(E_\Gamma)),$$

where $\hat{\otimes}_w$ means the weak closure of the algebraic tensor product, i.e. the tensor product in the category of von Neumann algebras.

Taking the canonical trace on $VN(\Gamma)$ and the usual trace on $B(L^2(E_\Gamma))$ we obtain a type II$_\infty$ trace on the algebra $B(L^2(E))^\Gamma$ which we will denote by tr$_\Gamma$. The corresponding dimension function is denoted by dim$_\Gamma$. Assume for example, a $\Gamma$-invariant operator $T$ on $L^2(E)$ is given as integral operator with a smooth kernel $k_T$, then a computation shows

$$\text{tr}_\Gamma(T) = \int_\mathcal{F} \text{tr}(k_T(x, x))dx.$$

The reader familiar with the Selberg trace formula will immediately recognize this as the "term of the identity".

### 1.2.2   The Heat Kernel

Since the operator $D$ is the pullback of an operator on $\Gamma \backslash X$, it commutes with the $\Gamma$-action and so does its heat operator $e^{-tD}$. This heat operator $e^{-tD} \in B(L^2(E))^\Gamma$ has a smooth kernel $< x \mid e^{-tD} \mid y >$. We get:

$$\mathrm{tr}_\Gamma e^{-tD} = \int_{\mathcal{F}} \mathrm{tr} < x \mid e^{-tD} \mid x > dx.$$

From this we read off that $\mathrm{tr}_\Gamma e^{-tD}$ satisfies the same small time asymptotics as $\mathrm{tr} e^{-tD_\Gamma}$ [BGV].

Let $D' = D \mid_{(ker(D))^\perp}$. Unfortunately very little is known about large time asymptotics of $\mathrm{tr}_\Gamma e^{-tD'}$ (see [LL]). Let

$$\mathrm{GNS}(D) \overset{\text{def}}{=} \sup\{\alpha \in \mathbb{R} \mid \mathrm{tr}_\Gamma e^{-tD'} = O(t^{-\alpha/2}) \text{ as } t \to \infty\}$$

denote the **Gromov-Novikov-Shubin invariant** of $D$ ([GrSh], [LL]). Then $\mathrm{GNS}(D)$ is always $\geq 0$. J. Lott showed in [Lo] that the GNS-invariants of Laplacians are homotopy invariants of a manifold. J.Lott and W. Lück conjecture in [LL] that the Gromov-Novikov-Shubin invariants of Laplace operators are always positive rational or $\infty$.

### 1.2.3   The $L^2$-Zeta Function

Throughout we will **assume** that the Gromov-Novikov-Shubin invariant of $D$ is positive. We consider the integral

$$\zeta^1_{D_\Gamma}(s) \overset{\text{def}}{=} \frac{1}{\Gamma(s)} \int_0^1 t^{s-1} \mathrm{tr}_\Gamma e^{-tD'} dt,$$

which converges for Re $(s) >> 0$ and extends to a meromorphic function on the entire plane which is holomorphic at $s = 0$, as is easily shown by using the small time asymptotics ([BGV],Thm 2.30).

Further the integral

$$\zeta^2_{D_\Gamma(s)}(s) \overset{\text{def}}{=} \frac{1}{\Gamma(s)} \int_1^\infty t^{s-1} \mathrm{tr}_\Gamma e^{-tD'} dt$$

converges for Re $(s) < \mathrm{GNS}(D)$, so in this region we define the $L^2$-**zeta function** of $D_\Gamma$ as

$$\zeta^{(2)}_{D_\Gamma}(s) \overset{\text{def}}{=} \zeta^1_{D_\Gamma}(s) + \zeta^2_{D_\Gamma}(s).$$

Assuming the Gromov-Novikov-Shubin invariant of $D$ to be positive we define the $L^2$-**determinant** of $D_\Gamma$ as

$$\det{}^{(2)}(D_\Gamma) \stackrel{\mathrm{def}}{=} \exp(-\frac{d}{ds}\mid_{s=0} \zeta^{(2)}_{D_\Gamma}(s)).$$

**Proposition 1.2.1** *Let $D_\Gamma$ denote a generalized Laplacian over the manifold $X_\Gamma$. Then the $L^2$-**characteristic function***: $\lambda \mapsto \det{}^{(2)}(D_\Gamma + \lambda)$, $\lambda > 0$, *extends to a holomorphic function on $\mathbb{C}\backslash(-\infty, 0]$.*

**Proof:** Let $\lambda \in \mathbb{C}\backslash(-\infty, 0]$ and $s \in \mathbb{C}$ with Re $(s) >> 0$. Then the operator $(D_\Gamma + \lambda)^{-s}$, defined by the spectral theorem, has a continuous kernel whose diagonal, as function in $S$, continues to a meromorphic function in $s$, holomorphic at $s = 0$ [Sh]. So the same counts for the zeta function

$$\begin{aligned}
\zeta^{(2)}_{D_\Gamma + \lambda}(s) &= \mathrm{tr}_\Gamma(D_\Gamma + \lambda)^{-s} \\
&= \int_{\mathcal{F}} \mathrm{tr} < x \mid (D_\Gamma + \lambda)^{-s} \mid x > dx.
\end{aligned}$$

By the proof of Thm 12.1 in [Sh] it is clear that $\zeta^{(2)}_{D_\Gamma + \lambda}(s)$ depends holomorphically on $\lambda$. Whence the proposition. ■

## 1.3 Convergence of Determinants

In this section we will show that when the fundamental group is residually finite, the characteristic functions of a tower converge to the $L^2$-characteristic function.

### 1.3.1 Residually Finite Groups

Let $\Gamma$ be an infinite group. A **tower** of subgroups is a sequence $\Gamma_1 \supset \Gamma_2 \supset \ldots$ of normal subgroups of $\Gamma$, each of which has finite index in $\Gamma$ and the sequence satisfies $\cap_j \Gamma_j = \{1\}$. A group $\Gamma$ that has a tower is called **residually finite**.

Finitely generated subgroups of $GL_n(\mathbb{C})$ are residually finite [Alp]. See [Kir] for more criteria of groups to be residually finite.

As a special example from arithmetic consider the following: Let $H$ denote a nonsplit quaternion algebra over $\mathbb{Q}$ which splits at $\infty$. Consider the group schemes $H^*$ and $H^1$ of units and norm one elements in $H$. Since via the splitting homomorphism $\pi : H^*(\mathbb{R}) \to GL_2(\mathbb{R})$ the norm becomes the determinant, $H^1(\mathbb{R})$ is mapped to $SL_2(\mathbb{R})$. Let $R$ denote an order in $H$ then $\Gamma(1) = \pi(H^1(R))$ is a cocompact discrete subgroup of $SL_2(\mathbb{R})$. For an ideal $I$ in $R$ let $\Gamma(I) = \{\gamma \in \Gamma(1) \mid \gamma \equiv 1 \bmod I\}$ the principal congruence subgroup. Now let $(I_n)_{n \in \mathbb{N}}$ be a decreasing sequence of nontrivial ideals in $R$ with $\cap_j I_j = 0$ then the groups $\Gamma_j = \Gamma(I_j)$ form a tower.

### 1.3.2   The Convergence Theorem

Now assume $\Gamma = \pi_1(X_\Gamma)$ is residually finite and fix a tower $(\Gamma_j)_{j \in \mathbb{N}}$. For any continuous function $f$ on $\mathbb{R}$ we define the operator $f(D)$ by means of the spectral theorem as follows: Let $\mu$ denote a spectral resolution of $D$, i.e. $\mu$ is a projection valued measure on the spectrum $\operatorname{Spec}(D)$ of the operator $D$ such that on $L^2(E)$ we have

$$D = \int_{\operatorname{Spec}(D)} \lambda \, d\mu(\lambda).$$

We then define

$$f(D) \stackrel{\mathrm{def}}{=} \int_{\operatorname{Spec}(D)} f(\lambda) \, d\mu(\lambda).$$

It follows

$$\operatorname{tr}_\Gamma(f(D)) = \int_{\operatorname{Spec}(D)} (f(\lambda)) \, d\operatorname{tr}_\Gamma \mu(\lambda),$$

where the integral will converge if $f(D)$ has a sufficiently smooth kernel.

Choose $n \in \mathbb{N}$ so large that $(1 + D_{\Gamma_j})^{-n}$ is of trace class and has a Schwartz kernel which is $(\dim X + 1)$-times continuously differentiable. Let $F_{ev}^n$ denote the space of even $C^\infty$-functions on $\mathbb{R}$ which decrease more rapidly than $(1 + x^2)^{-n}$, i.e. for $f \in C^\infty(\mathbb{R})_{even}$ we have

$$f \in F_{ev}^n \Leftrightarrow \lim_{x \to \infty} f(x)(1 + x^2)^n = 0.$$

Note that the even Schwartz functions lie in $F_{ev}^n$.

**Theorem 1.3.1**   *(Compare [Don].) For any $f \in F_{ev}^n$ we have*

$$\frac{\operatorname{tr} f(\sqrt{D_{\Gamma_j}})}{[\Gamma : \Gamma_j]} \; \stackrel{j \to \infty}{\longrightarrow} \; \operatorname{tr}_\Gamma f(\sqrt{D}).$$

*For any Schwartz function $g \in \mathcal{S}(\mathbb{R})$ we have*

$$\frac{\operatorname{tr} g(D_{\Gamma_j})}{[\Gamma : \Gamma_j]} \; \stackrel{j \to \infty}{\longrightarrow} \; \operatorname{tr}_\Gamma g(D).$$

Before proving the theorem we give some immediate applications:

**Corollary 1.3.1**   *(Kazhdan inequality, compare [Lü]) Let $h(D_{\Gamma_j})$ denote the dimension of the kernel of $D_{\Gamma_j}$ and $h^{(2)}(D_\Gamma)$ the $\Gamma$-dimension of $\ker(D)$ then we have*

$$\limsup_j \frac{h(D_{\Gamma_j})}{[\Gamma : \Gamma_j]} \leq h^{(2)}(D_\Gamma).$$

We say that $D$ has a **spectral gap at zero** if there is $\epsilon > 0$ such that $\mathrm{Spec}D \cap ]0, \epsilon[= \emptyset$. Let

$$N_j(x) \stackrel{\mathrm{def}}{=} \sum_{\lambda \leq x} \dim \mathrm{Eig}(D, \lambda).$$

**Corollary 1.3.2** *Assume there is an $x > 0$ such that*

$$\frac{N_J(x) - N_j(0)}{[\Gamma : \Gamma_j]} \stackrel{j \to \infty}{\longrightarrow} 0.$$

*then $D$ has a spectral gap at zero.*

**Proof of Corollary 1:** Let $F$ denote the set of all $f \in \mathcal{S}(\mathbb{R})$ with $f(0) = 1$, $f(x) \geq 0$ for all $x \in \mathbb{R}$, then we have

$$h(D_{\Gamma_j}) \leq \mathrm{tr} f(D_{\Gamma_j})$$

for all $f \in F$. It follows

$$\limsup_j \frac{h(D_{\Gamma_j})}{[\Gamma : \Gamma_j]} \leq \lim_j \frac{\mathrm{tr} f(D_{\Gamma_j})}{[\Gamma : \Gamma_j]} = \mathrm{tr}_\Gamma f(D)$$

for all $f \in F$. So that

$$\limsup_j \frac{h(D_{\Gamma_j})}{[\Gamma : \Gamma_j]} \leq \inf_{f \in F} \mathrm{tr}_\Gamma f(D) = h^{(2)}(D_\Gamma). \blacksquare$$

**Proof of Corollary 2:** Assume $x > 0$ as in the corollary. Let $f \in C_c^\infty(0, x)$ be positive, $f(y) \leq 1$ for all $y$, then

$$N_j(x) - N_j(0) \geq \mathrm{tr} f(D_{\Gamma_j}) \geq 0.$$

Hence for all such $f$ we have $\mathrm{tr}_\Gamma f(D) = 0$ which gives the claim. $\blacksquare$

**Proof of the theorem:** Fix $x_0 \in X$ and $v \in E_{x_0}^*$ the dual space to $E_{x_0}$. Recall that the distribution $u_s = \cos(s\sqrt{D})v(\delta_{x_0}) \in C^\infty(E)'$ satisfies the wave equation

$$(\frac{\partial^2}{\partial s^2} + D)u_s = 0.$$

By general results on hyperbolic equations ([Tayl], chap IV) it follows that $\cos(s\sqrt{D})$ has propagation speed $\leq | s |$. Now for $f \in F_{ev}^n$ the formula

$$f(\sqrt{D}) = \frac{1}{2\pi} \int_{\mathbb{R}} \hat{f} \cos(s\sqrt{D}) ds$$

implies that $f(\sqrt{D})$ has finite propagation speed if its Fourier Transform $\hat{f}$ has compact support. Let $PW$ denote the space of all $f \in F_{ev}^n$ such that $\hat{f}$ has compact support. Then $PW$ is dense in $f \in F_{ev}^n$ with respect to the norm

$$\| f \| = \sup_{x \in \mathbb{R}} | f(x) | (1 + x^2)^n.$$

Let $f \in PW$, we will prove the assertion of the theorem for $f$.

Schwartz kernels of operators on sections of the bundle $E_\Gamma$ are distributional sections of the bundle $E_\Gamma \otimes_{ext} E_\Gamma^*$ over $X_\Gamma \times X_\Gamma$. Where $\otimes_{ext}$ means the **exterior tensor product**, i.e.: $(E \otimes_{ext} F)_{(x,y)} \cong E_x \otimes F_y$. Since $X_\Gamma = \Gamma \backslash X$, these can be identified with $\Gamma \times \Gamma$-invariant sections of $E \otimes_{ext} E^*$.

**Lemma. 1.3.1** *Modulo the above identification we have for $f \in PW$ the identity of Schwartz kernels*

$$< \Gamma x \mid f(\sqrt{D_\Gamma}) \mid \Gamma y > = \sum_{\gamma \in \Gamma} < x \mid f(\sqrt{D}) \mid \gamma y > .$$

**Proof of the lemma:** Take $\varphi \in C^\infty(E_\Gamma)$ and identify $\varphi$ with its pullback to $E$, which is a $\gamma$-invariant section of $E$. The operator $\cos(s\sqrt{D})$ may be applied to $\varphi$ in the distributional sense. Since $\cos(s\sqrt{D})$ is $\Gamma$-invariant we obtain a $\Gamma$-invariant distribution, so $\cos(s\sqrt{D})\varphi \in C^\infty(E_\Gamma)'$. we thus get an operator $\cos(s\sqrt{D})$ from $C^\infty(E_\Gamma)$ to $C^\infty(E_\Gamma)'$ satisfying the same differential equation and initial value conditions as $\cos(s\sqrt{D_\Gamma})$, hence it equals the latter. For $f \in PW$ the formula $f(\sqrt{D}) = \frac{1}{2\pi} \int_\mathbb{R} \hat{f} \cos(s\sqrt{D}) ds$ gives $f(\sqrt{D}) \mid_{C^\infty(E_\Gamma)} = f(\sqrt{D_\Gamma})$. Since $f(\sqrt{D_\Gamma})$ is a smoothing operator we may write

$$\begin{aligned}
f(\sqrt{D})\varphi(x) \quad &= \int_X < x \mid f(\sqrt{D}) \mid y > \varphi(y) dy \\
&= \sum_{\gamma \in \Gamma} \int_{\mathcal{F}_\Gamma} < x \mid f(\sqrt{D}) \mid \gamma y > \varphi(y) dy,
\end{aligned}$$

where $\mathcal{F}_\Gamma$ denotes a fundamental domain of the $\Gamma$-action on $X$. Now the sum

$$\sum_{\gamma \in \Gamma} < x \mid f(\sqrt{D}) \mid \gamma y >$$

is locally finite and this gives the lemma. ∎

To continue the proof of the theorem fix fundamental domains $\mathcal{F}_\Gamma \subset \mathcal{F}_{\Gamma_1} \subset \mathcal{F}_{\Gamma_2} \subset \ldots$ in a way that there are representatives for the classes in $\Gamma/\Gamma_j$ such that $\mathcal{F}_{\Gamma_j} = \cup_{\sigma:\Gamma/\Gamma_j} \sigma \mathcal{F}_\Gamma$. We then get

$$\frac{\operatorname{tr} f(\sqrt{D_{\Gamma_j}})}{[\Gamma : \Gamma_j]} \quad = \quad \frac{1}{[\Gamma : \Gamma_j]} \sum_{\gamma \in \Gamma_j} \int_{\mathcal{F}_{\Gamma_j}} \operatorname{tr} < x \mid f(\sqrt{D}) \mid \gamma x > dx$$

$$= \sum_{\gamma \in \Gamma_j} \frac{1}{[\Gamma : \Gamma_j]} \sum_{\sigma : \Gamma / \Gamma_j} \int_{\mathcal{F}_\Gamma} \mathrm{tr} < \sigma x \mid f(\sqrt{D}) \mid \gamma \sigma x > dx$$

$$= \sum_{\gamma \in \Gamma_j} \frac{1}{[\Gamma : \Gamma_j]} \sum_{\sigma : \Gamma / \Gamma_j} \int_{\mathcal{F}_\Gamma} \mathrm{tr} < x \mid f(\sqrt{D}) \mid \sigma^{-1} \gamma \sigma x > dx.$$

The last equality stems from the $\Gamma$-invariance of the operator $f(\sqrt{D})$. Since $\Gamma_j$ is normal in $\Gamma$ we end up with

$$\frac{\mathrm{tr} f(\sqrt{D_{\Gamma_j}})}{[\Gamma : \Gamma_j]} = \sum_{\gamma \in \Gamma_j} \int_{\mathcal{F}_\Gamma} \mathrm{tr} < x \mid f(\sqrt{D}) \mid \gamma x > dx.$$

As $j \to \infty$ we get less and less summands, only the summand of $\gamma = e$ remains, but this is just $\mathrm{tr}_\Gamma f(\sqrt{D})$, so the claim follows for $f \in PW$.

Now let $f \in F_{ev}^n$ be arbitrary. From Theorem 1.4 in [CGT] we get that there is a constant $C > 0$ such that

$$\frac{\mathrm{tr}(1 + D_{\Gamma_j})^{-n}}{[\Gamma : \Gamma_j]} \leq C \quad \text{for all } j \in \mathbb{N}.$$

Let $\epsilon > 0$. To our given $f \in F_{ev}^n$ there is a $g \in PW$ such that

$$\sup_{x \in \mathbb{R}} \mid f(x) - g(x) \mid < \frac{\epsilon}{(1 + x^2)^n}.$$

We get

$$\begin{aligned}
\mid \frac{\mathrm{tr} f(\sqrt{D_{\Gamma_j}})}{[\Gamma : \Gamma_j]} - \mathrm{tr}_\Gamma f(\sqrt{D}) \mid \quad &\leq \mid \frac{\mathrm{tr} f(\sqrt{D_{\Gamma_j}})}{[\Gamma : \Gamma_j]} - \frac{\mathrm{tr} g(\sqrt{D_{\Gamma_j}})}{[\Gamma : \Gamma_j]} \mid \\
&+ \mid \frac{\mathrm{tr} g(\sqrt{D_{\Gamma_j}})}{[\Gamma : \Gamma_j]} - \mathrm{tr}_\Gamma g(\sqrt{D}) \mid \\
&+ \mid \mathrm{tr}_\Gamma g(\sqrt{D}) - \mathrm{tr}_\Gamma f(\sqrt{D}) \mid.
\end{aligned}$$

The first summand on the right hand side is less than $\frac{\mathrm{tr}((1 + D_{\Gamma_j})^{-n})}{[\Gamma : \Gamma_j]} \leq \epsilon C$. The second summand tends to zero as $j \to \infty$ and the third is less that $\epsilon \mathrm{tr}_\Gamma (1 + D)^{-n}$. The theorem follows. ∎

### 1.3.3 Convergence of Characteristic Functions

For later use we prove:

**Lemma. 1.3.2** *The function*

$$t \mapsto \frac{\mathrm{tr} e^{-t D_{\Gamma_j}}}{[\Gamma : \Gamma_j]} - \mathrm{tr}_\Gamma e^{-t D_\Gamma}, \quad t > 0,$$

*is rapidly decreasing at zero. Moreover this holds uniformly in j, that is, there is a function $f :]0,1[\to [0,\infty[$ which is rapidly decreasing at zero such that*

$$\mid \frac{\mathrm{tr}e^{-tD_{\Gamma_j}}}{[\Gamma:\Gamma_j]} - \mathrm{tr}_\Gamma e^{-tD_\Gamma} \mid \le f(t)$$

*for all $j \in \mathbb{N}$ and all $t \in ]0,1[$.*

**Proof:** Arguing as in the proof of the preceding lemma we see that

$$< \Gamma x \mid e^{-tD_\Gamma} \mid \Gamma y > = \sum_{\gamma \in \Gamma} < x \mid e^{-tD} \mid \gamma y >,$$

where convergence remains to be checked, however. To this end we write $d(x,y)$ for the Riemann distance of the points $x, y \in X$ and we use Theorem 3.1 of [CGT] to get for $d(x,y) \ge a > 0$:

$$\mid < x \mid e^{-tD} \mid \gamma y > \mid \le \frac{C}{\sqrt{t}} e^{-(d(x,y)-a)^2/t},$$

where $\mid . \mid$ here denotes the matrix $L^2$-norm, so $\mid A \mid = \sqrt{\mathrm{tr}AA^*}$.

By Proposition 4.1 in [CGT] we conclude that the function $x \mapsto e^{-(d(x,y)-a)^2/t}$ is in $L^1(X)$ and the $L^1$-norm is clearly bounded by a constant independent of $y$. From this the convergence of the sum above follows. As in the previous section we get:

$$\frac{\mathrm{tr}e^{-tD_{\Gamma_j}}}{[\Gamma:\Gamma_j]} - \mathrm{tr}_\Gamma e^{-tD_\Gamma} = \sum_{1 \ne \gamma \in \Gamma_j} \int_{\mathcal{F}_\Gamma} \mathrm{tr} < x \mid e^{-tD} \mid \gamma x > dx.$$

Now fix $x_0 \in \mathcal{F}_\Gamma$ and $r > 0$ such that the ball of radius $r$ around $x_0$ contains $\mathcal{F}_\Gamma$, i.e. $\mathcal{F}_\Gamma \subset B_r(x_0)$. Choose $R > \max(r,1)$ such that $y \in \mathcal{F}_\Gamma, x \in X$ with $d(x_0,x) > R$ implies $d(x,y) > r$. Choose $j_0 \in \mathbb{N}$ such that for all $j \ge j_0$ and all $\gamma \in \Gamma_j \backslash \{1\}$ we have $\gamma \bar{\mathcal{F}}_\Gamma \cap B_R(x_0) = \emptyset$. Under these circumstances we get for $j \ge j_0$:

$$\begin{aligned}
\mid \frac{\mathrm{tr}e^{-tD_{\Gamma_j}}}{[\Gamma:\Gamma_j]} - \mathrm{tr}_\Gamma e^{-tD_\Gamma} \mid &\le \sum_{\gamma \in \Gamma_j \backslash \{1\}} \frac{C}{\sqrt{t}} \int_{\mathcal{F}_\Gamma} e^{-(d(x,\gamma x)-r)^2/t} dx \\
&\le \frac{C}{\sqrt{t}} \int_{X \backslash B_R(x_0)} e^{-(d(x,B_r(x_0))-r)^2/t} dx \\
&\le \frac{C}{\sqrt{t}} \int_R^\infty e^{-(\xi-r)^2/t} \mathrm{vol}(d(x,x_0)=\xi) d\xi \\
&\le \frac{C_1}{\sqrt{t}} \int_R^\infty e^{-(\xi-r)^2/t} e^{C_2\xi} d\xi,
\end{aligned}$$

for some $C_1, C_2 > 0$ by Proposition 4.1 in [CGT]. But this can be calculated as

$$C_1 \int_{R/\sqrt{t}} e^{-(x-(\frac{r}{\sqrt{t}}+\frac{C_2\sqrt{T}}{2}))^2+r+\frac{C_2^2 t}{4}} dx \leq C_3 \int_{R/\sqrt{t}} e^{-x^2} dx$$

$$\leq C_3 \int_{R/\sqrt{t}} x e^{-x^2} dx \quad \text{for } t < 1$$

$$= \frac{C_3}{2} e^{-R^2/t}$$

which gives the assertion of the lemma. ∎

For the following theorem we need not require the Gromov-Novikov-Shubin invariant to be positive.

**Theorem 1.3.2** *As* $j \to \infty$ *we have*

$$\det(D_{\Gamma_j} + \lambda)^{\frac{1}{[\Gamma:\Gamma_j]}} \longrightarrow \det^{(2)}(D_\Gamma + \lambda)$$

*locally uniformly in* $\lambda \in \mathbb{C} \backslash (-\infty, 0]$.

**Proof:** For Re $\lambda > 0$ and $s \in \mathbb{C}$ we consider

$$F_j(\lambda, s) = \frac{\zeta_{D_{\Gamma_j}+\lambda}(s)}{[\Gamma:\Gamma_j]} - \zeta^{(2)}_{D_\Gamma+\lambda}(s)$$

$$= \frac{1}{\Gamma(s)} \int_0^\infty t^{s-1} (\frac{\text{tr} e^{-tD_{\Gamma_j}}}{[\Gamma:\Gamma_j]} - \text{tr}_\Gamma e^{-tD\Gamma}) e^{-t\lambda} dt$$

which converges by the previous lemma. The formula

$$F_j(\lambda, s) = \frac{\text{tr}(D_{\Gamma_j}+\lambda)^{-s}}{[\Gamma:\Gamma_j]} - \text{tr}_\Gamma (D+\lambda)^{-s} \quad \text{Re } (s) >> 0$$

together with the analytic continuation of the kernel of $(D_{\Gamma_j}+\lambda)^{-s}$ or $(D+\lambda)^{-s}$ as in ([Sh]) shows that $F_j$ extends to a holomorphic function on $(\mathbb{C} \backslash (-\infty, 0]) \times \mathbb{C}$. As Re $(\lambda) \to \infty$ we have $F_j(\lambda, s) \to 0$ uniformly in $\{$Re $(s) < K\}$ for all $K \in \mathbb{R}$. We want to show that this convergence is uniform in $j \in \mathbb{N}$. To this end let Re $(\lambda) > 0$ and split the integral above as $\int_0^1 + \int_1^\infty$. The integrand in the $\int_0^1$-part obeys a uniform estimate by the previous lemma and thus the convergence is uniform in this part.

For the $\int_1^\infty$-part recall that there is a constant $C > 0$ such that $\frac{\text{tr}(D_{\Gamma_j}+1)^{-n}}{[\Gamma:\Gamma_j]} < C$, from which it follows that $\frac{\text{tr} e^{-tD_{\Gamma_j}}}{[\Gamma:\Gamma_j]}$ is bounded for $t \geq 1$ with a bound not depending on j. From this it follows

$$F_j(\lambda, s) \overset{\text{Re } (\lambda) \to \infty}{\longrightarrow} 0$$

locally uniformly in s and uniformly in $j \in \mathbb{N}$.

On the other hand we have a differential equation

$$\frac{\partial}{\partial \lambda} F_j(\lambda, s) = F_j(\lambda, s+1).$$

For $\epsilon > 0$ let $G_\epsilon$ denote the region

$$G_\epsilon \overset{\text{def}}{=} \{z \in \mathbb{C} \mid \text{Re}\ (z) > -\frac{1}{\epsilon} \text{and dist}(z, (-\infty, 0]) > \epsilon\}.$$

From the above it follows that there is a $R > 0$ such that for Re $(s) > R$ we have

$$F_j(\lambda, s) \overset{j \to \infty}{\longrightarrow} 0 \quad \text{uniformly in } G_\epsilon \times \{\text{Re}\ (s) > R\}.$$

We want to show that, having this convergence for Re $(s) > R$ it already follows for Re $(s) > R - 1$. To achieve this let Re $(s) > R - 1$ and consider

$$F_j(\lambda, s) - F_j(\mu, s) = \int_\mu^\lambda F_j(z, s+1) dz,$$

which tends to zero as $j \to \infty$. But on the other hand

$$F_j(\lambda, s) - F_j(\mu, s) \overset{\mu \to \infty}{\longrightarrow} F_j(\lambda, s),$$

as the latter convergence is uniform in $j$ we get the claim. From this we get $F_j(\lambda, s) \overset{j \to \infty}{\longrightarrow} 0$

uniformly in a neighborhood of $s = 0$ and this gives the theorem. ∎

### 1.3.4   Convergence of the Determinant

Write $h(D_\Gamma)$ for the dimension of the kernel of $D_\Gamma$ and $h^{(2)}(D_\Gamma)$ for the $\Gamma$-dimension of the kernel of $D$. Lück showed in [Lü]:

$$\frac{h(\triangle_{p,\Gamma_j})}{[\Gamma : \Gamma_j]} \overset{j \to \infty}{\longrightarrow} h^{(2)}(\triangle_{p,\Gamma})$$

where $\triangle_p$ is the $p$-th Laplacian. Unfortunately the combinatorial arguments he used do not carry over to the analytic situation. So we only can show such a convergence under additional assumptions.

**Theorem 1.3.3** *Suppose that there is an $\epsilon > 0$ such that* $\mathrm{Spec}(D_{\Gamma_j}) \cap ]0, \epsilon[= \emptyset$ *for all $j$, then*

$$\frac{h(D_{\Gamma_j})}{[\Gamma : \Gamma_j]} \overset{j \to \infty}{\longrightarrow} h^{(2)}(\Gamma)$$

*and*

$$\det(D_{\Gamma_j})^{\frac{1}{[\Gamma : \Gamma_j]}} \overset{j \to \infty}{\longrightarrow} \det{}^{(2)}(D_\Gamma).$$

Examples of this are Shimura manifolds $X_\Gamma$, where $\Gamma$ has Kashdan property (T). This for example holds if the universal covering $X$ of $X_\Gamma$ is a symmetric space with every simple factor of rank $> 1$.

**Proof:** Apply the proof of the previous theorem to $D - \epsilon$. ∎

**Theorem 1.3.4** *Assume there exist constants $C, \alpha > 0$ such that*

$$\frac{\mathrm{tr}e^{-tD^r_{\Gamma_j}}}{[\Gamma : \Gamma_j]} \leq Ct^{-\alpha}$$

*for $t > 1$ and all $j \in \mathbb{N}$, then the Gromov-Novikov-Shubin invariant of $D$ is $\geq 2\alpha$ and*

$$\frac{h(D_{\Gamma_j})}{[\Gamma : \Gamma_j]} \overset{j \to \infty}{\longrightarrow} h^{(2)}(\Gamma)$$

*as well as*

$$\det(D_{\Gamma_j})^{\frac{1}{[\Gamma : \Gamma_j]}} \overset{j \to \infty}{\longrightarrow} \det{}^{(2)}(D_\Gamma).$$

Examples for this are The Laplacians of flat tori or Heisenberg manifolds.

**Proof:** We have that $\frac{\mathrm{tr}e^{-tD_{\Gamma_j}}}{[\Gamma:\Gamma_j]}$ tends to $\frac{h(D_{\Gamma_j})}{[\Gamma:\Gamma_j]}$ as $t \to \infty$ and by the assumption this is uniform in j. Thus from $\frac{\mathrm{tr}e^{-tD_{\Gamma_j}}}{[\Gamma:\Gamma_j]} \overset{j \to \infty}{\longrightarrow} \mathrm{tr}_\Gamma e^{-tD}$ and $\mathrm{tr}_\Gamma e^{-tD} \overset{t \to \infty}{\longrightarrow} h^{(2)}(D_\Gamma)$ the second assertion follows. The first is now easy and for the third recall that for $0 < \mathrm{Re}\ (s) < \alpha$ we have

$$\frac{\zeta_{D_{\Gamma_j}}(s)}{[\Gamma : \Gamma_j]} - \zeta_{D_\Gamma}^{(2)}(s) = \frac{1}{\Gamma(s)} \int_0^\infty t^{s-1} \left( \frac{\mathrm{tr}e^{-tD_{\Gamma_j}}}{[\Gamma : \Gamma_j]} - \frac{h(D_{\Gamma_j})}{[\Gamma : \Gamma_j]} - \mathrm{tr}_\Gamma e^{-tD} + h^{(2)}(D_\Gamma) \right) dt$$

We split this integral into a $\int_0^1$- and a $\int_1^\infty$-part. The $\int_0^1$-part defines an entire function converging locally uniformly to zero in a neighborhood of $s = 0$ as $j \to \infty$. In the $\int_1^\infty$-part the integrand is dominated by some constant times $t^{-\alpha}$ which allows us to interchange the integration with the limit as $j \to \infty$, so this part vanishes as $j \to \infty$. Together we get

$$\frac{\zeta_{D_{\Gamma_j}}(s)}{[\Gamma : \Gamma_j]} \overset{j \to \infty}{\longrightarrow} \zeta^{(2)}_{D_\Gamma}(s)$$

uniformly around s=0 which gives the last claim. ■

**Question 1.3.1** *Does for any generalized Laplacian the sequence $\frac{\dim \ker D_{\Gamma_j}}{[\Gamma : \Gamma_j]}$ converge to $\dim_\Gamma \ker D$?*

**Question 1.3.2** *Does for any generalized Laplacian the sequence $\det(D_{\Gamma_j})^{\frac{1}{[\Gamma : \Gamma_j]}}$ converge to $\det^{(2)}(D_\Gamma)$ if $\mathrm{GNS}(D_\Gamma) > 0$?*

Let for $j \in \mathbb{N}$ and $x \in \mathbb{R}$:

$$N_j(x) \overset{\mathrm{def}}{=} \sum_{\lambda \leq x} \dim \mathrm{Eig}(D_{\Gamma_j}; \lambda)$$

the **eigenvalue counting function** of $D_{\Gamma_j}$. The condition of the last theorem would follow from the existence of constants $a, b, C > 0$ such that

$$N_j(x) - N_j(0) \leq C(x^a + x^b) \quad x \geq 0,$$

for all $j \in \mathbb{N}$. If this holds, we say the **the counting functions satisfy a global growth estimate**.

In the following chapters we will specialize to the case where $X$ is a symmetric space without compact factors.

**Question 1.3.3** *Suppose $X$ is a symmetric space of the noncompact type and $D$ a generalized Laplacian which is invariant under the group of orientation preserving isometries of $X$. Do the counting functions of $D$ then satisfy a global growth estimate?*

# Chapter 2

# Torsion and Euler Characteristic

An introduction to the classical notions of analytical and combinatorial torsion in the beginning of this chapter is followed by a section dealing with the higher torsion and Euler characteristics which are the appropriate invariants to look at in the higher rank case.

## 2.1 Classical Torsion

### 2.1.1 General Torsion

A finite complex of Hilbert spaces $E_0 \overset{d}{\longrightarrow} E_1 \overset{d}{\longrightarrow} \ldots \overset{d}{\longrightarrow} E_n$ is called **admissible** if all

**Laplacians** $\triangle_p = dd^* + d^*d \mid_{E_p}$ are admissible in the sense of 1.1.1. Examples are finite dimensional complexes or the complexes of $L^2$-sections of elliptic complexes over Riemannian manifolds. For an admissible complex $E$ we define its **torsion** as

$$\tau_1(E) \overset{\text{def}}{=} \prod_{p=0}^{n} \det{}'(\triangle_p)^{p(-1)^p}.$$

### 2.1.2 Combinatorial and Analytic Torsion

Let $\mathcal{F}$ denote a locally constant sheaf over a connected oriented finite CW complex $X$ of dimension $n$ such that the stalks of $\mathcal{F}$ are finite dimensional unitary vector spaces, the scalar products being locally constant as well. This can be viewed as attaching to each cell $e$ a unitary space $\mathcal{F}_e$ and an unitary isomorphism $d_{e,e'} : \mathcal{F}_e \to \mathcal{F}_{e'}$ if $e'$ lies in the closure of $e$. Let

$$E_j \overset{\text{def}}{=} \bigoplus_{\substack{e \text{ cell} \\ \dim e = n - j}} \mathcal{F}_e$$

then the boundary operator $\partial$ makes $E = E_0 \to \ldots \to E_n$ a finite dimensional complex whose torsion was considered by Reidemeister [Reid], who showed that torsion and Betti numbers are the only numbers attached to $(X, \mathcal{F})$ which are stable under subdivision. It turned out that torsion is a homeomorphy invariant, but not a homotopy invariant and thus can be used to separate homotopy equivalent spaces which are not homeomorphic such as the lens spaces.

To give a locally constant hermitian sheaf amounts to the same as giving a unitary representation $\varphi$ of the fundamental group $\pi_1(X)$. This in turn defines a flat hermitian vector bundle $E_\varphi$ over X provided X is a smooth manifold. We then could also consider the complex of $E_\varphi$-valued differential forms $\Omega_\varphi^*$ and its torsion $\tau_1(\Omega_\varphi^*)$. It was conjectured by Ray and Singer [RS-AT] that the numbers $\tau_1(\Omega_\varphi^*)$ and $\tau_1(\mathcal{F})$ should coincide. This was later proven independently by Cheeger [Che] and Müller [Mü]. By this identity it is justified to write $\tau_1(\varphi)$ from now on.

We now give some properties of the torsion of a flat bundle.

**Proposition 2.1.1** *If the dimension of X is even we have $\tau_1(\varphi) = 1$.*

    **Proof:** [RS-AT]. ∎

**Proposition 2.1.2** *The number $\tau_1(\varphi)$ does not depend on the metrics.*

    **Proof:** [RS-AT]. ∎

**Proposition 2.1.3** *If $X, Y$ are oriented compact manifolds we have*

$$\tau_1(X \times Y, \varphi \otimes \psi) = \tau_1(X, \varphi)^{\chi(Y)} \tau_1(Y, \psi)^{\chi(X)}$$

*where $\chi$ denotes the Euler characteristic.*

    **Proof:** [RS-RM]. ∎

### 2.1.3   Holomorphic Torsion

Now assume the compact Riemannian manifold X is Kähler and let $\Omega_\varphi^{p,\cdot}$ denote the Dolbeault complex of $E_\varphi$-valued $(p, .)$-forms, where $E_\varphi$ is the flat holomorphic hermitian bundle attached to the unitary representation $\varphi$ of $\pi_1(X)$.

The **holomorphic torsion** of $\varphi$ is defined as

$$T_{\text{hol}}^p(\varphi) \overset{\text{def}}{=} \tau_1(\Omega_\varphi^{p,\cdot}).$$

We give some properties of $T_{\text{hol}}^p$:

**Proposition 2.1.4** *We have*

$$\prod_{p=0}^{\dim_{\mathbb{C}} X} T_{\text{hol}}^p(\varphi)^{(-1)^p} = 1.$$

**Proof:** [RS-HT]. ∎

**Proposition 2.1.5** *Let $X = \mathbb{C}/\Gamma$, where $\Gamma$ is the lattice generated by 1 and $z$, $\text{Im}(z) > 0$. Suppose $\varphi$ is a nontrivial character given by*

$$\varphi(mz + n) = e^{2\pi i(mu + nv)} \quad 0 \le u, v < 1.$$

*Then, for $p = 0, 1$ we have*

$$T_{\text{hol}}^p(\varphi) = \left| \frac{e^{-\pi i v^2 z}}{\Theta(u - zv, z)} \right|$$

*where $\Theta$ is the theta function*

$$\Theta(w, z) = -e^{\pi i(w + z/6)} \prod_{k=-\infty}^{\infty} (1 - e^{2\pi i(|k|z - \epsilon_k w)})$$

$$\epsilon_k = \text{sign } (k + \frac{1}{2}) = \frac{k + \frac{1}{2}}{|k + \frac{1}{2}|}.$$

**Proof:** [RS-CM]. ∎

The holomorphic torsion of compact hermitian symmetric spaces is computed in [Koe].

Holomorphic torsion enters into Arakelov theory in the following way: Let $X$ denote a compact Kähler manifold and $E$ a holomorphic hermitian vector bundle over $X$. Consider the cohomology groups $H^q(X, E)$ of the $E$-valued Dolbeault complex and the one dimensional vector space

$$\lambda(E) = \bigotimes_{q \ge 0} \wedge^{\max} H^q(X, E)^{(-1)^q},$$

where $V^{-1}$ denotes the dual of a vector space $V$. The $L^2$-scalar product induces a metric $h_{L^2}$ on $\lambda(E)$. Quillen [Qu] defined a new metric $h_Q$ on $\lambda(E)$ by

$$h_Q = h_{L^2} T_{\text{hol}}^0(E).$$

Now let $f : X \to Y$ be a smooth proper map of Kähler manifolds. Let $E$ be a holomorphic hermitian vector bundle on $X$ and let $f_* \lambda(E) = \det R f_*(E)$ be the determinant of the direct image of $E$ as defined in [KnMf]. This is a holomorphic line bundle $f_* \lambda(E)$ on $Y$ such that for every $y \in Y$

$$f_* \lambda(E)_y = \bigotimes_{y \ge 0} \wedge^{\max} H^q(X_y, E)^{(-1)^q}.$$

It then follows [BGS I-III], that the Quillen metric defined fibrewise is smooth. Further its curvature was computed loc.cit.

Thus analytic torsion became part of the direct image in the arithmetic Chow ring in which setting Gillet and Soulé were able to prove an arithmetic Riemann-Roch Theorem relating the determinant of the direct image in the K-theory to the direct image in the arithmetic Chow ring [GS].

## 2.2  Higher Torsion and Euler Characteristics

### 2.2.1  Determinants of Laplacians

We extend the notion of an admissible complex to infinite complexes $E = E_0 \to E_1 \to \ldots$ of Hilbert spaces in the obvious manner. So $E$ is called **admissible** if all Laplacians $\triangle_p = dd^* + d^*d \mid_{E_p}$ $p = 0, 1, 2, \ldots$ are admissible.

Let $\mathcal{C}_a$ denote the category whose objects are the admissible positive complexes and whose morphisms are the homomorphism of complexes $\varphi : E \to F$ such that the domain of $\varphi_j$ contains all eigenvectors of the Laplacian and $\varphi$ interchanges the Laplace operators, i.e.: $\varphi\triangle_E = \triangle_F\varphi$. Call a virtual complex $E = E_+ - E_-$ **pseudofinite**, if there is a $k_0$ such that for all $k \geq k_0$ we have

$$\det'(\triangle \mid E_k) = 1.$$

(Recall our notational convention that says $\det(\triangle' \mid E_k) = \frac{\det'(\triangle|(E_+)_k)}{\det'(\triangle|(E_-)_k)}$.) For a pseudofinite virtual complex $E$ we define the **determinant** of the Laplacian on $E$ to be

$$\det'(\triangle \mid E) = \prod_{k=0}^{\infty} \det'(\triangle \mid E_k)^{(-1)^k}.$$

**Lemma. 2.2.1** *Assume $E$ is a finite complex in $\mathcal{C}_a$ then*

$$\det'(\triangle \mid E) = 1.$$

**Proof:** Let $v$ be an eigenvector of $d^*d$ for some nonvanishing eigenvalue. Then $dv$ is an eigenvector for $dd^*$ for the same eigenvalue. So we have

$$\det'(d^*d \mid E_p) = \det'(dd^* \mid E_{p+1}).$$

Denoting $E$ as $E_0 \to \ldots \to E_N$ we get

$$
\begin{aligned}
\det'(\triangle \mid E) &= \prod_{p=0}^{N} \det'(dd^* \mid E_p)^{(-1)^p} \det'(d^*d \mid E_p)^{(-1)^p} \\
&= \prod_{p=0}^{N} \det'(dd^* \mid E_p)^{(-1)^p} \det'(dd^* \mid E_{p+1})^{(-1)^p} \\
&= 1. \blacksquare
\end{aligned}
$$

### 2.2.2  Twists and Torsion

For a virtual complex $E$ we consider the **twist** $E'$ of $E$ defined by

$$E' = \sum_{k=0}^{\infty} (-1)^k E[-k],$$

where $E[j]$ denotes the $j$-shifted complex, so $E[j]_k = E_{k+j}$. We also define the **higher twists** by iteration, so $E^{(0)} = E$, $E^{(n+1)} = (E^{(n)})'$.

If $E$ is written as $E_+ - E_-$ we write $E' = E'_+ - E'_-$, where

$$
\begin{aligned}
(E'_+)_j &= \sum_{k=0}^{\infty} (E_+)_{j-2k} + \sum_{k=0}^{\infty} (E_-)_{j-2k+1} \\
(E'_-)_j &= \sum_{k=0}^{\infty} (E_-)_{j-2k} + \sum_{k=0}^{\infty} (E_+)_{j-2k+1}.
\end{aligned}
$$

Since the $E_{\pm}$ are positive complexes the sums are finite.

**Proposition 2.2.1** *Let $E$ be a pseudofinite virtual complex with $\det'(\triangle \mid E) = 1$ then the twist $E'$ is a pseudofinite virtual complex and the determinant of the twist $E'$ equals the inverse torsion of $E$, i.e.*

$$\det'(\triangle \mid E') = \tau_1(E)^{-1}.$$

**Proof:** By our formulas for $(E_{\pm})_j$ we have

$$\det'(\triangle \mid (E')_j) = \prod_{k=0}^{\infty} \frac{\det'(\triangle \mid E_{j-2k})}{\det'(\triangle \mid E_{j-2k+1})}.$$

By $\det(\triangle \mid E) = 1$ it follows that fir $j$ large we have

$$\det'(\triangle \mid (E')_j) = \det'(\triangle \mid E)^{(-1)^j} = 1,$$

so that $E'$ is a pseudofinite complex. It follows that for some $N \in \mathbb{N}$

$$
\begin{aligned}
\det'(\triangle \mid E') &= \prod_{j=0}^{N} \det'(\triangle \mid E'_j)^{(-1)^j} \\
&= \prod_{j=0}^{N} \prod_{k=0}^{\infty} \frac{\det'(\triangle \mid E_{j-2k})^{(-1)^j}}{\det'(\triangle \mid E_{j-2k+1})^{(-1)^j}} \\
&= \prod_{n=0}^{\infty} \det'(\triangle \mid E_n)^{(-1)^n(N-n+1)} \\
&= (\prod_{n=0}^{\infty} \det'(\triangle \mid E_n)^{n(-1)^n})^{-1} \\
&= \tau_1(E)^{-1}. \blacksquare
\end{aligned}
$$

This might lead us to the definition of **higher torsion** as follows: Define the **zeroth torsion** of a pseudofinite complex $E$ as

$$\tau_0(E) \overset{\text{def}}{=} \det{}'(\triangle \mid E)^{-1}.$$

Suppose, the torsion numbers $\tau_0(E), \ldots, \tau_{r-1}(E)$ are defined and are all trivial, then the $r$-th twist $E^{(r)}$ is pseudofinite and we define the **$r$-th torsion** as

$$\tau_r(E) \overset{\text{def}}{=} \tau_0(E^{(r)}).$$

**Proposition 2.2.2** *Let $E$ be a pseudofinite virtual complex whose $r$-th torsion is defined, i.e. the twists $E^{(0)}, \ldots, E^{(r-1)}$ have trivial determinants, then the $r$-th torsion is given by the formula*

$$\tau_r(E) = \prod_{p=0}^{\infty} \det{}'(\triangle \mid E_p)^{((-1)^{p+r-1}\binom{p}{r})},$$

*where the product is in fact finite.*

**Proof:** By induction on $r$. The cases $r = 0, 1$ are clear. For the step $r - 1 \to r$ assume $\tau_0(E^{(0)}) = \ldots = \tau_0(E^{r-1}) = 1$, then

$$
\begin{aligned}
\tau_r(E) &= \tau_0(E^{(}r)) = \tau_0((E')^{(r-1)}) \\
&= \prod_{p=0}^{N} \det{}'(\triangle|(E')_p)^{(-1)^{p+r}\binom{p}{r-1}} \\
&= \prod_{p=0}^{N} \prod_{k=0}^{p} \det{}'(\triangle|E_{p-k})^{-1)^{k+p+r}\binom{p}{r-1}} \\
&= \prod_{p=0}^{N} \det{}'(\triangle|E_n)^{(-1)^{n+r}(\sum_{p=n}^{N}\binom{p}{r-1})}
\end{aligned}
$$

Now because of

$$\sum_{p=n}^{N} \binom{p}{r-1} = \binom{N+1}{r} - \binom{n}{r}$$

the claim follows. ∎

For the moment assume, we are in the situation of 2.1.4. again, so $X_\Gamma$ is a compact Kähler manifold with universal covering $X$ and fundamental group $\Gamma$. Fix a finite dimensional unitary representation $\varphi$ of $\Gamma$. Then we have

**Proposition 2.2.3** *Denoting by $T^p_{\mathrm{hol}}$ the p-th holomorphic torsion we have*

$$\prod_{p=0}^{\dim_{\mathbb{C}} X} T^p_{\mathrm{hol}}(\varphi)^{p(-1)^p} = \tau_2(\varphi)^{-\frac{1}{2}}.$$

**Proof:** We use the previous proposition to calculate

$$
\begin{aligned}
\log(\tau_2(\varphi)^{-\frac{1}{2}}) &= \frac{1}{2}\sum_{n=0}^{\dim X}\binom{n}{2}\log\det(\triangle_{n,\varphi})\\
&= \frac{1}{2}\sum_{n=0}^{\dim X}(-1)^n\sum_{p+q=n}\binom{p+q}{2}\log\det(\triangle_{p,q,\varphi})\\
&= \frac{1}{2}\sum_{n=0}^{\dim X}(-1)^n\sum_{p+q=n}(p^2+2pq+q^2-p-q)\log\det(\triangle_{p,q,\varphi})\\
&= \sum_{n=0}^{\dim X}(-1)^n\sum_{p+q=n}pq\log\det(\triangle_{p,q,\varphi})\\
&= \sum_{p=0}^{\dim_{\mathbb{C}} X}p(-1)^p\sum_{q+0}^{\dim_{\mathbb{C}} X}q(-1)^q\log\det(\triangle_{p,q,\varphi}). \blacksquare
\end{aligned}
$$

### 2.2.3 Higher Euler Characteristics

Let $\mathcal{C}^+$ denote the category of complexes of $\mathbb{C}$-vector spaces which are bounded from below and have degreewise finite dimensional cohomology, i.e. $\dim H^j(E) < \infty$ for all $j$. Let $\mathcal{K}^+$ denote the Grothendieck group of $\mathcal{C}^+$, define the twist as before, so for $E \in \mathcal{K}^+$ the higher twists $E^{(1)}, \ldots$ are defined as in the previous section. We need to extend the notion of an **Euler characteristic** to infinite virtual complexes by

$$\chi(E) \overset{\mathrm{def}}{=} \sum_{k=-\infty}^{\infty}(-1)^k\dim H^k(E),$$

provided $\dim H^k(E) = \dim H^k(E_+) - \dim H^k(E_-)$ vanishes for almost all $k$.

As an example consider the de Rham complex over a compact oriented manifold $X$. Let there be a fibre bundle

$$\mathbb{T} \to X \to B,$$

where $\mathbb{T} = \{z \in \mathbb{C} : |z| = 1\}$ is the circle group. Let $X_T$ denote the vector field on $X$ spanned by the $\mathbb{T}$-action. **Assume** that there is a real closed 1-form $\omega$ such that $\omega(X_T)$ nowhere vanishes. (Note that this assumption forces the first Betti number of $X$ to be nonzero.) By integration we may assume $\omega$ to be $\mathbb{T}$-invariant. Let $F$ denote the subcomplex of $E$ consisting of the $\mathbb{T}$-invariant forms. We get

$$F = \Omega(X)^{\mathbb{T}} = \Omega(B) \oplus \Omega(B) \wedge \omega,$$

where $\Omega(B)$ is identified with its pullback to $X$. Since $\omega$ is closed, this decomposition is stable under the exterior differential. Since furthermore $\mathbb{T}$ is connected, it acts trivially on cohomology, so the cohomology of $X$ is computed by $F$. This gives

$$\chi(E') = \chi(F') = \chi(B).$$

So we get the Euler characteristic of $B$ from the de Rham complex of $X$.

One is thus led to define the higher Euler characteristics as $\chi(E^{(r)})$ for a virtual complex E but unfortunately this does not always make sense, simply because the cohomology of the virtual complex $E^{(r)}$ needn't vanish in high degrees.

So call a complex **cohomologically finite** if $H^j(E) = 0$ for large j, in other words, the total cohomology H(E) is finite dimensional.

**Observation:** *Let the virtual complex E be cohomologically finite and assume that the Euler characteristic $\chi(E)$ vanishes. Then the twist $E'$ is cohomologically finite.*

So start with a cohomologically finite virtual complex E. If $E^{(1)}, \ldots, E^{(r)}$ are cohomologically finite we have

$$\chi(E^{(0)}) = \ldots = \chi(E^{(r-1)}) = 0$$

and

$$\chi(E^{(r)}) = (-1)^r \sum_{j=0}^{\infty} \binom{j}{r} (-1)^j \, dim H^j(E).$$

This is easily proven by induction on r. This motivates the following Definition: The **r-th Euler characteristic** of a cohomologically finite virtual complex E is defined by

$$\chi_r(E) := (-1)^r \sum_{j=0}^{\infty} \binom{j}{r} (-1)^j \, dim H^j(E).$$

The r-th Euler characteristic of a manifold is given by any complex that computes its cohomology. So in case of a fibration

$$\mathbb{T} \to X \to B,$$

as above we have

$$\chi_1(X) = \chi(B).$$

To every compact manifold M we now can attach a sequence of Euler numbers

$$\chi_0(M), \ldots, \chi_n(M),$$

where n is the dimension of M. The most significant of these is, as we shall see, the first nonvanishing one, so define the **generic Euler number** of M as

$$\chi_{gen}(M) = \chi_r(M), \quad \text{where } r \text{ is the least index with } \chi_r(M) \neq 0.$$

**Proposition 2.2.4** *Let M,N be compact Manifolds. We have*

$$\chi_{gen}(M \times N) = \chi_{gen}(M)\chi_{gen}(N).$$

**Proof:** Let $r, s \geq 0$. By the Künneth-rule we get

$$\chi_{r+s}(M \times N) = \sum_{j=0}^{\infty} \left( \begin{array}{c} j \\ r+s \end{array} \right) (-1)^j \sum_{\mu+\nu=j} a_\mu b_\nu,$$

where $a_\mu, b_\nu$ are the Betti numbers of M and N. Suppose now that we have

$$\chi_0(M) = \ldots = \chi_{r-1}(M) = \chi_0(N) = \ldots = \chi_{s-1}(N) = 0,$$

then we have

$$\begin{array}{ll} \chi_r(M)\chi_s(N) & = \sum_{r'+s'=r+s} \chi_{r'}(M)\chi_{s'}(N) \\ & = \sum_{r'+s'=r+s} \sum_{j=0}^{\infty} (-1)^j \sum_{\mu+\nu=j} \left( \begin{array}{c} \mu \\ r' \end{array} \right) \left( \begin{array}{c} \nu \\ s' \end{array} \right) a_\mu b_\nu \\ & = \sum_{j=0}^{\infty} (-1)^j \sum_{\mu+\nu=j} a_\mu b_\nu \sum_{r'+s'=r+s} \left( \begin{array}{c} \mu \\ r' \end{array} \right) \left( \begin{array}{c} \nu \\ s' \end{array} \right) \end{array}$$

because of

$$\sum_{r'+s'=r+s} \left( \begin{array}{c} \mu \\ r' \end{array} \right) \left( \begin{array}{c} \nu \\ s' \end{array} \right) = \left( \begin{array}{c} \mu+\nu \\ r+s \end{array} \right)$$

the claim follows. ∎

### 2.2.4 Lie Algebra Cohomology

Just to give another example of a situation in which higher Euler characteristics occur we will describe a typical situation in Lie algebra cohomology. Further examples are given in chapter 3.

We consider a short exact sequence

$$o \to \mathfrak{n} \to \mathfrak{l} \to \mathfrak{a} \to 0$$

of finite dimensional complex (just for convenience) Lie algebras where $\mathfrak{a}$ is abelian. In such a situation a $\mathfrak{l}$-module $V$ is called **acceptable** , if the $\mathfrak{a}$-module $H^q(\mathfrak{n}, V)$ is a finite dimensional $\mathfrak{a}$-module. Note that $V$ itself needn't be finite dimensional.

**Example 1.** Any finite dimensional $\mathfrak{l}$-module will be acceptable.

**Example 2.** Let $\mathfrak{g}_0$ denote the Lie algebra of a semisimple Lie group $G$ of the Harish-Chandra class, i.e. $G$ is connected and has a finite center. Let $K$ be a maximal compact subgroup of $G$ and let $G = KAN$ be an Iwasawa decomposition of $G$. Write the corresponding decomposition of the complexified Lie algebra as $\mathfrak{g} = \mathfrak{k} \oplus \mathfrak{a} \oplus \mathfrak{n}$. Now let $\mathfrak{l} = \mathfrak{a} \oplus \mathfrak{n}$ with the structure of a subalgebra of $\mathfrak{g}$. Consider an admissible $(\mathfrak{g}, K)$-module $V$. A theorem of [HeSchm] assures us that $V$ then is an acceptable $\mathfrak{l}$-module.

**Theorem 2.2.1** *Let*

$$o \to \mathfrak{n} \to \mathfrak{l} \to \mathfrak{a} \to 0$$

*be an exact sequence of finite dimensional complex Lie algebras with $\mathfrak{a}$ abelian. Let $V$ be an acceptable $\mathfrak{l}$-module then with $r = \dim(\mathfrak{a})$ we have*

$$\chi_0(H^*(\mathfrak{l}, V)) = \ldots = \chi_{r-1}(H^*(\mathfrak{l}, V)) = 0,$$

*and*

$$\chi_r(H^*(\mathfrak{l}, V)) = \chi_0(H^*(\mathfrak{n}, V)^{\mathfrak{a}}),$$

*where $H^*(\mathfrak{n}, V)^{\mathfrak{a}}$ denotes the $\mathfrak{a}$-invariants in $H^*(\mathfrak{n}, V)$.*

**Proof:** There is a Hochschild-Serre spectral sequence to the module $V$, which we denote by $E$. So

$$E_2^{p,q} = H^p(\mathfrak{a}, H^q(\mathfrak{n}, V)),$$

and

$$E_\infty^{p,q} = Gr^q(H^{p+q}(\mathfrak{l}, V))$$

for some filtration on $H^{p+q}(\mathfrak{l}, V)$. We will consider the **twist** of $E$ defined as

$$(E')_s^{p,q} \overset{\text{def}}{=} \sum_{k=0}^{\infty} (-1)^k E_s^{p-k,q},$$

which is a **virtual spectral sequence** , i.e. a formal difference of two spectral sequences. One can iterate this process to obtain a sequence

$$E^{(0)} = E,\ E^{(1)} = E',\ldots, E^{(n)} = (E^{(n-1)})',\ldots$$

of virtual spectral sequences. Of particular interest will be

$$F \overset{\text{def}}{=} E^{(r)}, \quad r = \dim(\mathfrak{a}).$$

To understand the 2-term of our spectral sequences we prove:

**Lemma. 2.2.2** *With* $r = \dim(\mathfrak{a})$ *and* $M$ *any finite dimensional* $\mathfrak{a}$*-module we have*

$$\dim H^p(\mathfrak{a}, M) = \binom{r}{p} \dim(M^{\mathfrak{a}}).$$

**Proof:** At first we consider the case $\dim(M) = 1$, so $\mathfrak{a}$ operates through a linear functional $\lambda$. We will compute the Lie algebra cohomology via the standard complex , i.e. for a Lie algebra $\mathfrak{l}$, an $\mathfrak{l}$-module $V$ and $q \geq 0$ we consider the cochain space

$$C^q \overset{\text{def}}{=} \text{Hom}_{\mathbb{C}}(\wedge^q \mathfrak{l}, V)$$

and the differential $d : C^q \to C^{q+1}$ defined by

$$
\begin{aligned}
df(x_0 \wedge \ldots \wedge x_q) &= \sum_{i=0}^{q-1}(-1)^i x_i.f(x_0 \wedge \ldots \hat{x_i} \ldots \wedge x_q) \\
&+ \sum_{i<j}(-1)^{i+j}f([x_i, x_j] \wedge x_0 \wedge \ldots \hat{x_i} \ldots \hat{x_j} \ldots x_q).
\end{aligned}
$$

Since $\mathfrak{a}$ is abelian, the second sum vanishes in our situation and if $\lambda = 0$ the first does so too, so the differential is zero then and the claim for $\lambda = 0$ follows.

For the case $\lambda \neq 0$ assume $f \in C^q$ with $df = 0$. Let $H \in \mathfrak{a}$ with $\lambda(H) = (-1)^{q+1}$ then define $g \in C^{q-1}$ by

$$g(x_0 \wedge \ldots \wedge x_{q-2}) \overset{\text{def}}{=} f(x_0 \wedge \ldots \wedge x_{q-2} \wedge H).$$

We get

$$
\begin{aligned}
dg(x_0 \wedge \ldots \wedge x_{q-1}) &= \sum_i (-1)^i x_i.g(x_0 \wedge \ldots \hat{x_i} \ldots \wedge x_{q-1}) \\
&= \sum_i (-1)^i x_i.f(x_0 \wedge \ldots \hat{x_i} \ldots \wedge x_{q-1} \wedge H) \\
&= (-1)^{q+1} H.f(x_0 \wedge \ldots \wedge x_{q-1}) \\
&= f(x_0 \wedge \ldots \wedge x_{q-1}),
\end{aligned}
$$

so $dg = f$ and the claim follows.

Now that we have finished the case $\dim(M) = 1$ we will proceed via induction on $\dim(\mathfrak{a})$. Therefore assume $\dim(\mathfrak{a}) = 1$. In this case we use a subinduction on $\dim(M)$. The case $\dim(M) = 1$ has been dealt with. Using the Jordan normal form theorem we may assume that $\dim(M^{\mathfrak{a}}) = 1$. Consider the exact sequence

$$0 \to M^{\mathfrak{a}} \to M \to M/M^{\mathfrak{a}} \to 0,$$

the corresponding long exact sequence of cohomology gives a exact sequence

$$0 \to (M/M^{\mathfrak{a}})^{\mathfrak{a}} \to H^1(\mathfrak{a}, M^{\mathfrak{a}}) \to H^1(\mathfrak{a}, M) \to H^1(\mathfrak{a}, M/M^{\mathfrak{a}}) \to 0.$$

Using the induction assumption we get $\dim(M/M^{\mathfrak{a}})^{\mathfrak{a}} = \dim H^1(\mathfrak{a}, M^{\mathfrak{a}}) = 1$, so the middle arrow is zero, hence $\dim H^1(\mathfrak{a}, M) = \dim H^1(\mathfrak{a}, M/M^{\mathfrak{a}}) = 1$ by the induction assumption.

Now assume the assertion of the lemma proven for $\dim(\mathfrak{a}) \leq r$ and assume we are given an $\mathfrak{a}$ of dimension $r + 1$. Write $\mathfrak{a} = A + B$ with $\dim(A) = 1$. Consider the Hochschild-Serre spectral sequence $E$ with

$$E_2^{p,q} = H^p(A, H^q(B, M)),$$
$$E_\infty^{p,q} = Gr^q(H^{p+q}(\mathfrak{a}, M)).$$

Since $\dim(A) = 1$ we get $E_2^{p,q} \neq 0 \Rightarrow p = 0, 1$ so that this spectral sequence degenerates. Therefore

$$
\begin{aligned}
\dim H^k(\mathfrak{a}, M) &= \dim Gr^k(H^k(\mathfrak{a}, M)) + \dim Gr^{k-1}(H^k(\mathfrak{a}, M)) \\
&= \dim H^0(A, H^q(B, M)) + \dim H^1(A, H^{q-1}(B, M)) \\
&= \dim H^q(B, M^A) + \dim H^{q-1}(B, M^A) \\
&= \dim M^{\mathfrak{a}} \left( \begin{pmatrix} r \\ q \end{pmatrix} + \begin{pmatrix} r \\ q-1 \end{pmatrix} \right). \blacksquare
\end{aligned}
$$

$\blacksquare$

From this lemma we conclude that in the Grothendieck group of vector spaces we have

$$F_2^{p,q} \cong \begin{cases} H^q(\mathfrak{n}, V)^{\mathfrak{a}} & \text{if } p = 0, \\ 0 & \text{else.} \end{cases}$$

This means that we have the identity

$$\chi_0(H^*(\mathfrak{n}, V)^{\mathfrak{a}}) = \chi_0(F_2^{0,*}) = \chi(F_2),$$

where $\chi(F_2)$ is the **Euler characteristic of the spectral sequence** $F$ , which is defined by

$$\chi(F_2) \overset{\text{def}}{=} \sum_{p,q} (-1)^{p+q} \dim(F_2^{p,q}).$$

Note that, since the differential of a spectral sequence always maps a component to a component of the other parity $(= (-1)^{p+q})$ we get

$$\chi(F_2) = \chi(F_3) = \ldots = \chi(F_\infty)$$

and this number is going to be denoted by $\chi(F)$. We conclude

$$
\begin{aligned}
\chi_0(H^*(\mathfrak{n}, V)^{\mathfrak{a}}) &= \sum_{p,q}(-1)^{p+q}\dim(F_\infty^{p+q}) \\
&= \sum_{q=0}^{\infty}(-1)^q\sum_{p=0}^{\infty}(-1)^p\begin{pmatrix} r \\ p \end{pmatrix}\dim(Gr^p(H^{p+q}(\mathfrak{l}, V))) \\
&= \sum_{p=0}^{\infty}(-1)^p\begin{pmatrix} r \\ p \end{pmatrix}\dim(H^p(\mathfrak{l}, V)). \\
&= \chi_r(H^*(\mathfrak{l}, V)). \blacksquare
\end{aligned}
$$

## 2.3   $L^2$-Torsion

### 2.3.1   $L^2$-Torsion

We now return to the notation of 1.2, so $X_\Gamma$ is an oriented compact smooth Riemannian manifold with universal covering $X$ and fundamental group $\Gamma$. Let $\varphi$ denote a finite dimensional unitary representation of $\Gamma$. Let $F_\varphi$ denote the flat hermitian vector bundle over $X_\Gamma$ defined by $\varphi$. We denote by $\triangle_{\varphi,p}$ the Laplace operator on $F_\varphi$-valued $p$-forms. We assume that the Gromov-Novikov-Shubin invariants of all $\triangle_{\varphi,p}$ are positive and define the $L^2$**-torsion** of $\varphi$ by

$$
\tau^{(2)}(\varphi) \overset{\text{def}}{=} \prod_{p=0}^{\dim X}\det{}^{(2)}(\triangle_{\varphi,p}')^{p(-1)^p}.
$$

Assuming $X_\Gamma$ Kählerian we also define the $p$**-th holomorphic torsion** as

$$
T^{(2)}_{\text{hol},p} \overset{\text{def}}{=} \prod_{q=0}^{\dim_\mathbb{C} X}\det{}^{(2)}(\triangle_{\varphi,p,q})^{q(-1)^q},
$$

where $\triangle_{\varphi,p,q}$ is the Hodge-Laplacian on $F_\varphi$-valued $(p,q)$-forms. We do not consider higher twists of $L^2$-torsion since they will play no rôle in the sequel.

**Proposition 2.3.1** *We have* $\tau^{(2)}(\varphi) = 1$ *if* $\dim X$ *is even.*

  **Proof:**  [Lo]. $\blacksquare$

**Proposition 2.3.2** *If* $X, Y$ *are oriented compact manifolds we have*

$$
\tau^{(2)}(X \times Y) = \tau^{(2)}(X)^{\chi(Y)}\tau^{(2)}(Y)^{\chi(X)}
$$

*where* $\chi$ *denotes the Euler characteristic.*

**Proof:** We easily get

$$\mathrm{tr}_{\pi_1(X) \times \pi_1(Y)}(e^- t\triangle_{X \times Y, p}) = \sum_{a+b=p} \mathrm{tr}_{\pi_1(X)}(e^{-t\triangle_X, a})\mathrm{tr}_{\pi_1(Y)}(e^{-t\triangle_Y, b}).$$

From this the claim is derived by a purely formal computation. ∎

# Chapter 3

# The Geometry of Shimura Manifolds

This chapter contains some preparatory remarks on the geometry of Shimura manifolds, their sphere bundles and the geodesic homotopy classes. The higher Euler characteristics of the latter are computed in terms of root systems. This makes an important step in our computation of orbital integrals but it is to some extent a detour since a different choice of Haar measures would directly give the final result identifying the factor that goes with the orbital integral as the higher Euler characteristic of the homotopy class of the corresponding closed geodesic. But unfortunately we are forced to stick to certain standard normalizations since we pick up the results of Harish-Chandra on orbital integrals.

## 3.1 The Geometry of Geodesics

### 3.1.1 Geodesics on the Universal Covering

Let $X_\Gamma$ denote a compact orientable locally symmetric manifold whose universal cover $X$ is globally symmetric without compact factors. We will call such a space $X_\Gamma$ a **Shimura manifold**. If the group $\Gamma$ is arithmetic we say that $X_\Gamma$ is an **arithmetic Shimura manifold**. Do not mix this term with the notion of a **Shimura variety** , which is a special arithmetic Shimura manifold together with an algebraic structure [Del].

The group $G$ of orientation preserving isometries of $X$ is a centerless semisimple Lie group which acts transitively on $X$. More precisely let $K$ denote a maximal compact subgroup of $G$ then $X \cong G/K$. The reflection at the point $eK$ in $X$ defines a **Cartan involution** $\theta$ on $G$ that fixes $K$ pointwise. We get a Cartan decomposition $\mathfrak{g}_0 = \mathfrak{k}_0 + \mathfrak{p}_0$ where $\mathfrak{g}_0$ and $\mathfrak{k}_0$ denote the Lie algebras of $G$ and $K$ and $\mathfrak{p}_0$ is the $-1$-eigenspace of the Cartan involution $\theta$. Let $\mathfrak{a}_0$ denote a maximal abelian subspace of $\mathfrak{p}_0$ and $A_0 = \exp \mathfrak{a}_0$ the corresponding abelian subgroup of $G$. Let $\mathfrak{a}_0^+$ denote a positive Weyl chamber in $\mathfrak{a}_0$ and $A_0^+ = \exp(\mathfrak{a}_0^+)$. Every oriented geodesic $c$ in $X$ can be "turned into $A_0^+$ " the following way: At first we may assume $c(0) = eK$ since $G$ acts transitively on $X$. Then, modulo K we may assume that $c(t) = exp(tH_c)K$ for some $H_c \in \overline{\mathfrak{a}_0^+} =$ closure of the Weyl chamber = interior points plus walls and the origin. Let $W_0, \ldots, W_n$ denote a complete set of nonconjugate walls of $\mathfrak{a}_0^+$ with $W_0 = \mathfrak{a}_0^+$ and let

$$M_j = cent(W_j \mid K)$$

the centralizer of $W_j$ in $K$. If $H_c \in W_j$, then the isotropy group of $c$ is $M_j$ and we thus get a decomposition of the tangent bundle of $X$ as a G-set:

$$TX - 0 \cong \bigcup_{j=0}^{n} (G/M_j) \times W_j \quad \text{as a } G - \text{set.}$$

The sphere bundle thus becomes

$$SX \cong \bigcup_{j=0}^{n} (G/M_j) \times SW_j \quad \text{as a } G - \text{set,}$$

where $SW_j$ denotes the set of norm 1 elements in $W_j$.

### 3.1.2 Geodesics "Downstairs"

Now let $\Gamma \subset G$ be the fundamental group of $X_\Gamma$ then $\Gamma$ is a torsion free lattice in G and

$$X_\Gamma = \Gamma \backslash G/K.$$

Further, let $\gamma \in \Gamma$ then, since $\Gamma$ is the fundamental group of $X_\Gamma$, the element $\gamma$ defines a free homotopy class of closed paths in $X_\Gamma$ and this gives a bijection

$$\{\text{nontrivial conjugacy classes in } \Gamma\} \leftrightarrow \{\text{free homotopy classes of closed paths in } X_\Gamma\}.$$

The sphere bundle over $X_\Gamma$ is given as:

$$SX_\Gamma \cong \bigcup_{j=0}^{n} (\Gamma \backslash G/M_j) \times SW_j$$

Let $\gamma \in \Gamma$ and $G_\gamma, \Gamma_\gamma, K_\gamma$, denote the centralizers of $\gamma$ in $G, \Gamma, K$. Let $\gamma$ be in G conjugate to an element of $W_j K$ then modulo conjugation we have $K_\gamma \subset M_\gamma$. So we have a map

$$X_\gamma = \Gamma_\gamma \backslash G_\gamma / K_\gamma \to \Gamma \backslash G/M \hookrightarrow SX.$$

By Proposition 5.15 of [DKV] we get that $X_\gamma$ actually injects in $SX_\Gamma$ and $X_\Gamma$ and that the image of $X_\gamma$ is the union of all closed geodesics which map to closed paths in $X_\Gamma$ which are in the free homotopy class of $\gamma$.

## 3.2 Euler Characteristics

### 3.2.1 On Higher Euler Characteristics

**Proposition 3.2.1** *Let $\chi_n$ denote the n-th Euler characteristic. If $\gamma$ is conjugate to an element of $WK$, where $W$ is a wall of dimension $d$ then $\chi_n(X_\gamma) = 0$ for $n = 0, \ldots, d - 1$.*

**Proof:** Let $A \cong \mathbb{R}^d$ the group generated by W. There is a Cartan subgroup which contains $A$ and $\gamma$, hence we have $A \subset G_\gamma$. Furthermore the split part of $\gamma$ is regular in $A$ so we have $G_\gamma \subset Z(A)$, the centralizer of $A$. From this we deduce that $A$ acts on $X_\gamma \cong \Gamma_\gamma \backslash G_\gamma / K_\gamma$ by $a : \Gamma_\gamma \, g \, K_\gamma \mapsto \Gamma_\gamma \, ga \, K_\gamma$. Moreover A is the central split component of $G_\gamma$. Now by the considerations in [MS] we see that there is a torsionfree normal subgroup $\Gamma_\gamma^0$ of finite index in $\Gamma_\gamma$ such that the A-action turns $X_\gamma^0 = \Gamma_\gamma^0 \backslash G_\gamma / K_\gamma$ into a $\mathbb{T}^d$-fibre bundle. Let p denote the projection $X_\gamma^0 \to X_\gamma$ and let S denote the finite group $\Gamma_\gamma / \Gamma_\gamma^0$. The cohomology of $X_\gamma$ is computed by the complex of S-invariant forms on $X_\gamma^0$, which we denote by $\Omega(X_\gamma^0)^S$. Since the torus $\mathbb{T}^d$ is compact and connected and the actions of $\mathbb{T}^d$ and S on $X_\gamma^0$ commute, this cohomology is also computed by $\Omega(X_\gamma^0)^{S\mathbb{T}^d}$. We will show that the latter can be written as

$$\Omega(X_\gamma^0)^{S\mathbb{T}^d} \cong \bigoplus_{0 \leq i_1 < \ldots < i_n \leq d} \omega_{i_1} \wedge \ldots \wedge \omega_{i_n} \wedge p^*\Omega(X_\gamma^0/\mathbb{T}^d)^S,$$

where the $\omega_i$ are closed 1-forms dual to the vector fields that span the $\mathbb{T}^d$-action. Since the $\omega_j$ are closed this is a decomposition of the de Rham complex and gives an equivalent decomposition of the cohomology which gives the claim.

To show the existence of $\omega_1, \ldots, \omega_d$ recall that differentiable 1-forms on the homogeneous space $G_\gamma / K_\gamma$ are given by

$$\Omega_1^\infty(G_\gamma/K_\gamma) \cong (C^\infty(G_\gamma) \otimes \mathfrak{p}_\gamma)^{K_\gamma}.$$

So the $G_\gamma$-invariant forms are in one to one correspondence with $\mathfrak{p}_\gamma^{K_\gamma}$ and this space contains the wall $\mathfrak{a}$ with $\exp(\mathfrak{a}) = A$. Let now $\omega_1, \ldots \omega_d$ correspond to a basis of $\mathfrak{a}$ then these forms will be harmonic hence closed. ∎

The number $\chi_d(X_\gamma)$ is computed in terms of root systems in [D-Hitors] and to state the result we have to introduce some notation.

Assume $\gamma$ lies in a Cartan $H$ which is supposed to have minimal splitrank among those Cartans that contain $\gamma$. We say $H$ is **fundamental relative** $\gamma$.

We write $\mathfrak{g}$ and $\mathfrak{k}$ for the complexifications of the Lie algebras $\mathfrak{g}_0$ and $\mathfrak{k}_0$. Let $\mathfrak{h}_0$ be the Lie algebra of $H$ and $\mathfrak{h}$ its complexification then $\mathfrak{h}$ is a Cartan subalgebra of $\mathfrak{g}$ and we write $\Phi(\mathfrak{h}, \mathfrak{g})$ for the **root system** of the pair $(\mathfrak{h}, \mathfrak{g})$. Let $x \mapsto x^c$ denote the complex conjugation on $\mathfrak{g}$ with respect to the real form $\mathfrak{g}_0$. The complex conjugation acts on roots via $\alpha^c(x) = \overline{\alpha(x^c)}$. A root is called **imaginary** if $\alpha^c = -\alpha$, **real** if $\alpha^c = \alpha$ and **complex** if neither holds.

Let $\mathfrak{g} = \mathfrak{h} \oplus \oplus_\alpha \mathfrak{g}_\alpha$ denote the **root space decomposition**. To every root space $\mathfrak{g}_\alpha$ there exists a generator $X_\alpha$ such that

$$[X_\alpha, X_{-\alpha}] = H_\alpha \qquad \text{with } B(H_\alpha, H) = \alpha(H),$$

$$B(X_\alpha, X_{-\alpha}) = 1$$

and $X_\alpha^c = X_{\alpha^c}$ if $\alpha$ complex and $X_\alpha^c = \pm X_{-\alpha}$ if $\alpha$ is imaginary. An imaginary root is called **compact** if $X_\alpha^c = -X_{-\alpha}$ and noncompact otherwise. Let $\Phi_n$ denote the set of noncompact imaginary roots and $\Phi^+$ a set of positive roots as well as $\Phi_n^+ = \Phi_n \cap \Phi^+$.

Let $\Phi_\gamma$ denote the set of roots of $(\mathfrak{g}_\gamma, \mathfrak{h})$ where $\mathfrak{g}_\gamma$ is the complexified Lie algebra of the centralizer $G_\gamma$. Write $\Phi_{\gamma,n}$, $\Phi_{\gamma,n}^+$ accordingly.

Let $\nu_\gamma = \dim(G_\gamma/K_\gamma) - \operatorname{rank}(G_\gamma/K_\gamma)$, the latter denoting the rank as symmetric space.

Denote for any compact subgroup $L$ of $G$ the **standard volume** [HC-HA1, sec.7] of $L$ by $v_B(L)$ and let vol denote the **normalized volume** (loc. cit.).

For any pair of groups $G \supset H$ let $W(G, H)$ denote the group theoretical Weyl group, i.e.:

$$W(G, H) \overset{\text{def}}{=} N_G(H)/Z_G(H),$$

where $N_G(H)$ is the normalizer of $H$ in $G$ and $Z_G(H)$ is the centralizer of $H$ in $G$.

For any topological group $G$ let $G^0$ denote the connected component of the unit element. Finally write $\rho_\gamma$ for the half of the sum of the positive roots of $(\mathfrak{g}_\gamma, \mathfrak{h})$.

Let $A_\gamma$ denote the flat spanned by $\gamma$, i.e. $A_\gamma = \cap_{\gamma \subset A} A$, where the intersection runs over all maximal flats containing $\gamma$. Now let $\lambda_\gamma = \operatorname{vol}(A_\gamma)$.

A discrete subgroup $\Gamma$ of $G$ is called **nice** if for every $\gamma \in \Gamma$ the set of eigenvalues of the adjoint $Ad(\gamma)$ does not contain any root of unity. Every arithmetic group contains a finite index nice subgroup [Bor,17.1].

**Proposition 3.2.2** *Assume that $\Gamma$ is nice and let $\gamma \in \Gamma$ and $H$ as above of splitrank $r$ Assume $\gamma$ is of splitrank $r$ too. The action of the split part $H_\mathbb{R}$ of $H$ on $X_\gamma$ makes $X_\gamma$ the total space of a $(S^1)^r$-fibre bundle and*

$$\chi_r(X_\gamma) = \frac{|\, W(\mathfrak{g}_\gamma, \mathfrak{h}) \,| \prod_{\alpha \in \Phi_\gamma^+} (\rho_\gamma, \alpha)}{\lambda_\gamma c_\gamma [G_\gamma : G_\gamma^0]} \operatorname{vol}(\Gamma_\gamma \backslash G_\gamma).$$

*with*

$$c_\gamma = (-1)^{|\Phi_{\gamma,n}^+|} (2\pi)^{|\Phi_\gamma^+|} 2^{\nu_\gamma/2} \frac{v(H_I)}{v(K_\gamma)} \,|\, W(G_\gamma, H) \,| \,.$$

**Proof:** [D-Hitors]. ∎

### 3.2.2 The Fundamental Rank

We proved Proposition 3.2.1 with direct geometric arguments, i.e. without use of spectral theory and harmonic analysis. Introducing these tools we are now able to extend this proposition considerably.

We will first introduce the notion of the fundamental rank: The **rank** of $G$ is the common dimension of all Cartan subgroups of $G$. The **splitrank** is the maximal dimension of the split part of a Cartan subgroup and the **fundamental rank** is minimal dimension of the split part of a Cartan subgroup. In other words the fundamental rank $FR(G)$ of $G$ is

$$FR(G) = \operatorname{rank}(G) - \operatorname{rank}(K).$$

It is by definition that $G$ has a compact Cartan subgroup if and only if $FR(G) = 0$. By root space decomposition one finds that

$$FR(G) \equiv \dim X \mod (2).$$

As an example consider $G = SL_n(\mathbb{R})$ for $n \geq 2$. Then we have $K = SO(n)$ and so we get $\operatorname{rank}(G) = n - 1$, $\operatorname{rank}(K) = [\frac{n}{2}]$ and so

$$FR(SL_n(\mathbb{R})) = n - [\frac{n}{2}] - 1.$$

### 3.2.3 The Pseudocuspform Condition

Let now $\triangle_k$ denote the $k$-th Laplacian on $X$. The space of smooth $k$-forms on $X$ can be described as

$$\Omega_k^\infty = (C^\infty(G) \otimes \wedge^k \mathfrak{p})^K,$$

where K acts via right translations on $C^\infty(G)$ and via the adjoint representation on $\wedge^k \mathfrak{p}$. The heat operator $e^{-t\triangle_k}$ has a rapidly decreasing kernel $h_t^k$ in

$$(\mathcal{S}(G) \otimes \operatorname{End}(\wedge^k \mathfrak{p}))^{K \times K},$$

where $\mathcal{S}(G) = \cap_{p>0} \mathcal{C}^p(G)$ is Harish-Chandra's Schwartz space [BM] and $K \times K$ acts on $\mathcal{S}(G)$ via right and left translation and on $\operatorname{End}(\wedge^k \mathfrak{p})$ via the adjoint representation. We put for $r \geq 0$

$$f_t^r \stackrel{\text{def}}{=} \sum_{k=0}^{\dim X} (-1)^{k+r} \left( \begin{array}{c} k \\ r \end{array} \right) \operatorname{tr} h_t^k,$$

where tr means the trace in $\operatorname{End}(\wedge^k \mathfrak{p})$.

We want to use the inversion formula for orbital integrals as in [HC-S]. We compute the principal series contributions to the inversion formula. So let $H$ denote a $\theta$-stable Cartan subgroup, let $H_\mathbb{R}$ denote its split component and $H_I$ the compact component then we have a direct product $H = H_\mathbb{R} H_I$. Let $P$ denote a cuspidal parabolic subgroup with Langlands decomposition $P = MH_\mathbb{R}N$. Let $(\xi, W_\xi)$ be an irreducible unitary representation of $M$ and $e^\nu$ a quasicharacter of $H_\mathbb{R}$. Set $\pi_{\xi,\nu} = Ind_P^G(\xi \otimes e^\nu \otimes 1)$ the principal series representation. This is the right regular representation of G on the space

$$H_{\xi,\nu} = \{f : G \to W_\xi \mid f(m\ exp(h_\mathbb{R})ng) = \xi(m)e^{(\rho+\nu)(h_\mathbb{R})}f(g)\}.$$

Recall the notion of a pseudocuspform [MS2]: A rapidly decaying function f on G is called a **pseudocuspform** if $tr\ \pi_{\xi,\nu}(f) = 0$ whenever H is not fundamental.

**Proposition 3.2.3** *Assume the fundamental rank of G is positive.The function $f_t^r$ with $r = FR(G)$ is a pseudocuspform. More precisely we have:*

$$tr\ \pi_{\xi,\nu}(f_t) = \begin{cases} 0 & if\ H\ is\ not\ fundamental \\ -e^{t\pi_{\xi,\nu}(C)}dim(W_\xi \otimes \wedge^*\mathfrak{p}_\mathfrak{h})^{K\cap M^0} & otherwise \end{cases}.$$

*Here C is the Casimir operator of G and for a finite dimensional vector space V we write*

$$\wedge^* V = \sum_{k=0}^{dim\ V} (-1)^k \wedge^k V.$$

*Further $\mathfrak{p}_\mathfrak{h}$ denotes the orthocomplement of $\mathfrak{h}_\mathbb{R}$ in $\mathfrak{p}$, where $\mathfrak{h}_\mathbb{R} = LieH_\mathbb{R}$.*

*For $r < FR(G)$ we have*

$$tr\ \pi_{\xi,\nu}(f_t) = 0.$$

**Proof:** The decomposition $\mathfrak{p} = \mathfrak{h}_\mathbb{R} + \mathfrak{p}_\mathfrak{h}$ is $K \cap M$-stable. As $K \cap M$-modules we get

$$\wedge^k\mathfrak{p} \cong \sum_{i+j=k} \wedge^j\mathfrak{h}_\mathbb{R} \otimes \wedge^i\mathfrak{p}_\mathfrak{h}.$$

Recall that $\mathfrak{h}_\mathbb{R}$ is a trivial $K \cap M$-module. Let $r' = dim\mathfrak{h}_\mathbb{R}$. We get by Frobenius reciprocity:

$$
\begin{aligned}
tr\pi_{\xi,\nu}(f_t) &= tr\pi_{\xi,\nu}(\sum_{k=0}^{d} \binom{k}{r}(-1)^k e^{-t\tilde{\triangle}_k}) \\
&= e^{t\pi_{\xi,\nu}(C)}\sum_{k=0}^{d}\binom{k}{r}(-1)^k dim(H_{\xi,\nu} \otimes \wedge^k\mathfrak{p})^K \\
&= e^{t\pi_{\xi,\nu}(C)}\sum_{k=0}^{d}\binom{k}{r}(-1)^k dim(W_\xi \otimes \wedge^k\mathfrak{p})^{K\cap M} \\
&= e^{t\pi_{\xi,\nu}(C)}\sum_{k=0}^{d}\sum_{j=0}^{k}(-1)^k\binom{k}{r}\binom{r'}{j} dim(W_\xi \otimes \wedge^{k-j}\mathfrak{p}_\mathfrak{h})^{K\cap M} \\
&= e^{t\pi_{\xi,\nu}(C)}\sum_{i=0}^{d}\sum_{j=0}^{d}(-1)^{i+j}\binom{i+j}{r}\binom{r'}{j} dim(W_\xi \otimes \wedge^i\mathfrak{p}_\mathfrak{h})^{K\cap M}.
\end{aligned}
$$

Now let for $r, r' \in \mathbb{N}$; $i \in \mathbb{N}_0$:

$$A(i, r, r') = \sum_{j=0}^{r'} \binom{i+j}{r} \binom{r'}{j} (-1)^j.$$

The claim now follows from the fact that r=r' if and only if H is fundamental and the elementary

**Lemma. 3.2.1** *We have*

$$A(i, r, r') = (-1)^{r'} \binom{i}{r - r'}.$$

Now consider the Laplace operator $\triangle_{k,\Gamma}$ on $X_\Gamma$ and the heat operator $e^{-t\triangle_{k,\Gamma}}$. Plugging the function $f_t^r$ into the **Selberg Trace Formula** we get

$$\sum_{k=0}^{\dim X} (-1)^{k+r} \binom{k}{r} \operatorname{tr}\, e^{-t\triangle_{k,\Gamma}} = \sum_{[\gamma] \subset \Gamma} \operatorname{vol}(\Gamma_\gamma \backslash G_\gamma) \mathcal{O}_\gamma(f_t^r).$$

The sum on the right runs over all conjugacy classes in $\Gamma$ and the expression $\mathcal{O}_\gamma(f)$ stands for the **orbital integral**:

$$\mathcal{O}_\gamma(f) = \int_{G/G_\gamma} f(g\gamma g^{-1}) dg.$$

Note that the left hand side of the trace formula tends to $\chi_r(X_\Gamma)$ as $t \to \infty$.

Now by Harish-Chandra's formula for the Fourier transform of the orbital integrals and the fact that $\operatorname{tr}\pi_{\xi,\nu}(f_t^r)$ vanishes for $r < FR(G)$ we get that for $r < FR(G)$ all orbital integrals vanish. This gives

**Theorem 3.2.1** *All the Euler Characteristics* $\chi_0(X_\Gamma), \ldots, \chi_{FR(G)-1}(X_\Gamma)$ *vanish.* ∎

## 3.3   Haar Measures

The bilinear form $B$ on $\mathfrak{g}$ which we always take to be a positive multiple of the Killing form induces a Haar measure on the group $G$ the following way: Firstly $B$ defines a positive inner product on the vector space $\mathfrak{p}_0$, which is $K$-invariant. By the differential of the projection $G \to G/K = X$ the space $\mathfrak{p}_0$ is mapped isomorphically to the tangent space of $X$ at $eK$. By translation we get a Riemannian metric on $X$. Let $\mu = \mu_B$ denote the volume form of this metric. Let $dk$ denote the Haar measure on $K$ with total volume 1. Then the formula

$$\int_G f d\mu = \int_X \int_K f(xk) dk d\mu_B(x)$$

defines a Haar measure also denoted $\mu$ om $G$ the so called **standard measure** (see [HC-HA1], sec. 7).

In the case when the fundamental rank is zero there is however a different measure which we are going to use, the **Euler-Poincaré measure** $\mu_{EP}$, which is determined by the property that

$$\mu_{EP}(\Gamma\backslash G) = (-1)^{\frac{\dim X}{2}}\chi(X_\Gamma)$$

for every torsionfree uniform lattice $\Gamma$. To express the relation between these two measures let $T$ be a compact Cartan subgroup and let

$$c = (\prod_{\alpha>0}(\alpha,\rho))(2\pi)^{-|\Phi^+|}2^{(\text{rank}\,G/K-\dim G/K)/2}\frac{v(K)|W(G_{\mathbb{C}},T_{\mathbb{C}})|}{v(T)|W(K,T)|}.$$

Then we have

$$\mu_{EP} = c\mu.$$

# Chapter 4

# Geometric Zeta Functions in the Rank One Case

This chapter contains a treatment of two kinds of zeta functions in the rank one case. They are expressed in a "complete" determinant formula which also gives a determinant expression for the "factor at infinity".

## 4.1 The Zeta Functions of Selberg and Ruelle

Before preceding in a more general setting let $Y$ denote a compact Riemannian surface of genus $g \geq 2$. There is a unique hyperbolic metric on $Y$ of curvature $-1$. Consider the infinite product

$$R(s) \overset{\text{def}}{=} \prod_c (1 - e^{-sl(c)})$$

convergent for Re $(s) >> 0$, where the product is taken over all primitive closed geodesics $c$ in $Y$ and $l(c)$ is the length of the geodesic $c$.

By the uniformisation theorem $Y$ is a Shimura manifold $Y = X_\Gamma = \Gamma \backslash \mathcal{H}_2(\mathbb{R})$, where $\mathcal{H}_2(\mathbb{R}) = \{z \in \mathbb{C} | z = x + iy, \ y > 0\}$ is the upper half plane in $\mathbb{C}$ with the hyperbolic metric $\frac{dx^2 + dy^2}{y^2}$. Then the group of orientation preserving isometries of $X$ is $G = PSL_2(\mathbb{R})^0 = SL_2(\mathbb{R})/\pm 1$, the connected component of $PSL_2(\mathbb{R})$. For the sphere bundle we actually have $SX_\Gamma = \Gamma \backslash G$. The complex Lie algebra $\mathfrak{g}$ of $G$ is $sl_2(\mathbb{C}) = \{M \in \text{Mat}_2(\mathbb{C}) | \text{tr} M = 0\}$. Let

$$H = \frac{1}{2} \begin{pmatrix} 1 & \\ & -1 \end{pmatrix} \in \mathfrak{g}$$

then via the identification $SX_\Gamma = \Gamma \backslash G$ the **geodesic flow** becomes

$$\Phi_t : SX_\Gamma = \Gamma \backslash G \to \Gamma \backslash G$$

$$\Phi_t(\Gamma g) = \Gamma g \exp(tH).$$

Now let $\mathfrak{n} = \left\{ \begin{pmatrix} 0 & * \\ 0 & 0 \end{pmatrix} \right\} \subset \mathfrak{g}$ and $\bar{\mathfrak{n}} = \left\{ \begin{pmatrix} 0 & 0 \\ * & 0 \end{pmatrix} \right\} \subset \mathfrak{g}$ then with $\mathfrak{a} = \mathbb{C}H$ the decomposition

$$\mathfrak{g} = \bar{\mathfrak{n}} \oplus \mathfrak{a} \oplus \mathfrak{n}$$

induces a decomposition of the tangent bundle of $G$ which pushes down to a decomposition of the tangent bundle of $SX_\Gamma = \Gamma \backslash G$ as

$$TSX_\Gamma = T^s \oplus T^0 \oplus T^u$$

into **stable**, **central** and **unstable** part, where the geodesic flow is contracting on $T^s$ and expanding on $T^u$, i.e. for $v \in T^s$ the metric norm of $D\Phi_t(v)$ decreases exponentially as $t \to \infty$ and for $v \in T^u$ the norm of $D\Phi_t(v)$ increases exponentially. Further the norm of $D\Phi_t(v)$ for $v$ in the line bundle $T^0$ remains constant. These properties are summarized in saying that $\Phi$ is an **Anosov flow**.

The function $R(s)$ above is called the **Ruelle zeta function** since it was introduced by David Ruelle in [Rue]. Ruelle introduced it in the more general setting of arbitrary Anosov flows which satisfy a mild analyticity condition. For these flows Ruelle was able to show that $R(s)$ admits an analytic continuation to a meromorphic function of finite order. In the general setting there is not very much more one can prove concerning $R(s)$. For a geodesic flow of a Riemann surface however, one has more at hand.

Atle Selberg [Sel] introduced the zeta function

$$Z(s) \stackrel{\text{def}}{=} \prod_c \prod_{N \geq 0} (1 - e^{-(s+N)l(c)})$$

and he showed that $Z(s)$ extends to an entire function whose zero set can be split into two parts, the **spectral zeroes** and the **topological zeroes**: The spectral zeroes all lie in the set $[0,1] \cup (\frac{1}{2} + i\mathbb{R})$ and s is a spectral zero if and only if $s(1-s)$ is an eigenvalue of the Laplacian $\triangle_\Gamma$ of $X_\Gamma$, the multiplicity of the zero being equal to the multiplicity of the eigenvalue. The topological zeroes are at $s = 0, -1, -2, \ldots$ and the order of a zero at $s = -k$ is (2k+1)(2g-2). The orders at s=0 add up.

All this information can be summarized to

**Theorem 4.1.1** *[CaVo] The Selberg zeta function satisfies the determinant formula:*

$$Z(s) = \det(\triangle_\Gamma + s(1-s))(e^{-(s-\frac{1}{2})^2} \det(P + s - \frac{1}{2}))^{2g-2},$$

*where $P = \sqrt{\triangle^d + \frac{1}{4}}$ and $\triangle^d$ is the Laplacian of the 2-sphere $S^2$ with the usual metric. This means that the eigenvalues of $P$ are the numbers $j + \frac{1}{2}$ for $j = 0, 1, \ldots$ and the eigenvalue $j + \frac{1}{2}$ has multiplicity $2j + 1$.*

From this we also get that the function $\hat{Z}(s) = Z(s)(e^{-(s-\frac{1}{2})^2}\det(P + s - \frac{1}{2}))^{2-2g}$ satisfies the functional equation

$$\hat{Z}(s) = \hat{Z}(1 - s).$$

Back to the Ruelle zeta function it is easy to see that

$$R(s) = \frac{Z(s)}{Z(s+1)}$$

so that all poles and zeroes of $R(s)$ may be read off. Even more, we find that $R(s)$ satisfies a functional equation as $s$ is replaced by $-s$ and in the center of the functional equation we find the following first Taylor coefficient:

$$R(s) = \frac{T_{\mathrm{hol}}(X_\Gamma)}{T_{\mathrm{hol}}^{(2)}(X_\Gamma)} s^{2g-2} + \text{higher order terms},$$

where $T_{\mathrm{hol}}$ is the holomorphic torsion, $T_{\mathrm{hol}}^{(2)}$ the holomorphic $L^2$-torsion.

This is going to be generalized in the following chapters.

## 4.2   Rank One Spaces

In this chapter we assume $X$ to be a rank one symmetric space. By definition this means that the group of orientation preserving isometries $G$ acts transitively on the sphere bundle $SX$ of $X$. By 3.1.1 it then follows that any maximal abelian subspace $\mathfrak{a}_0$ of $\mathfrak{p}_0$ has dimension one and thus we conclude that $G$ has splitrank one. There are three series and one sporadic space (see [Helg]), we have:

- the real hyperbolic space of dimension $n = 2, 3, \ldots$

$$X = \mathcal{H}_n(\mathbb{R}) = SO(n,1)^+/SO(n),$$

- the complex hyperbolic space of dimension $2n$; $n = 1, 2, \ldots$

$$X = \mathcal{H}_{2n}(\mathbb{C}) = SU(n,1)/S(U(n) \times U(1)),$$

- the quaternionic hyperbolic space of dimension $4n$; $n = 1, 2, \ldots$

$$X = \mathcal{H}_{4n}(\mathbb{H}) = Sp(n,1)/Sp(n),$$

  and

- the Cayley hyperbolic space of dimension 16

$$X = \mathcal{H}_{16}(\mathrm{Cal}) = (\mathfrak{f}_{4(-20)}, \mathfrak{so}(9)).$$

The only coincindences here are $\mathcal{H}_2(\mathbb{R}) \cong \mathcal{H}_2(\mathbb{C})$ and $\mathcal{H}_4(\mathbb{R}) \cong \mathcal{H}_4(\mathbb{H})$.

Since the splitrank of $G$ is one, the fundamental rank can only be zero or one which is in this setting equivalent to the dimension being even or odd. So we only have fundamental rank equals 1 in the case of $\mathcal{H}_n(\mathbb{R})$ for $n$ odd, in all other cases we have $FR(G) = 0$. The only spaces occurring here that are **hermitian symmetric spaces**, i.e. that have a $G$-invariant holomorphic structure are the complex hyperbolic spaces.

We will fix a maximal abelian space $\mathfrak{a}_0 \subset \mathfrak{p}_0$ and the corresponding split torus $A = \exp(\mathfrak{a}_0)$. Further let $M$ denote the centralizer of $A$ in $K$ and fix a minimal parabolic $P = MAN$.

Since it doesn't change the geometry we will not use $G = Iso^+(X)$ but a finite cover of this, namely we assume $G = \mathbf{G}(\mathbb{R})$, where $\mathbf{G}$ is a simply connected algebraic group.

### 4.2.1   K-separability

A finite dimensional representation $(\tau, V)$ of a maximal compact subgroup $K$ of $G$ defines a vector bundle $E_\tau$ over $X$ via $E_\tau \stackrel{\text{def}}{=} (G \times V)/K$ where $K$ acts on $G \times V$ from the right by $(g, v)k \stackrel{\text{def}}{=} (gk, \tau(k^{-1})v)$. The projection $E_\tau \to X$ is the obvious one: $(g, v)K \mapsto gK$. The group $G$ then acts on $E_\tau$ via $g'((g, v)K) \stackrel{\text{def}}{=} (g'g, v)K$ and this action makes $E_\tau$ a **homogeneous vector bundle**, that is, the $G$-action is smooth, maps fibres to fibres, is fibrewise linear and the projection $E_\tau \to X$ is $G$-equivariant.

The other way round given a homogeneous vector bundle $E$ over $X$ the $K$-action on the fibre over the point $eK \in X = G/K$ defines a finite dimensional representation $(\tau, V)$ of $K$ and it is easy to see that $E \cong E_\tau$. This defines an equivalence of categories between the category of homogeneous vector bundles over $X$ and the category of finite dimensional representations of the compact group $K$.

By rank$(X) = 1$ we have for the sphere bundle $SX$:

$$SX \cong G/M.$$

The pullback of any homogeneous vector bundle $E_\tau$ to $G/M$ will be a homogeneous vector bundle over $G/M$ corresponding to the $M$-representation $\tau|_M$. Actually, what we are going to do is to represent homogeneous vector bundles over $G/M$ as linear combinations of bundles over $G/K$, hence we have to study the restriction of representations from $K$ to $M$.

For any compact group $H$ let Rep$(H)$ denote the **Representation ring** of $H$, as abelian group Rep$(H)$ is the Grothendieck group of the category of finite dimensional unitary representations of $H$ and the tensor product makes it a ring.

The restriction gives a ring homomorphism res : Rep$(K) \to$ Rep$(M)$. A representation $\sigma$ of $M$ is called $K$-**separable** if $\sigma$ lies in the image of res.

**Proposition 4.2.1** *If* dim $X$ *is even then every finite dimensional representation of $M$ is $K$-separable. If* dim $X$ *is odd and $\sigma \in \hat{M}$ is not $K$-separable then $\sigma + \sigma^*$ is, where $\sigma^*$ is the dual representation. Furthermore, the constituents of the adjoint representation of $M$ on the space $\wedge^p \mathfrak{p}$ are $K$-separable for all $p$.*

**Proof:** Fix a maximal torus $T$ in $K$ such that $T_M = T \cap M$ is a maximal torus in $M$. Denote the corresponding Lie algebras by $\mathfrak{t}$ and $\mathfrak{t}_M$. For dim $X$ odd we have $K = Spin(2n+1)$, $M = Spin(2n)$ for some $n$. The claim follows easily ([BrotD]). Now assume dim X even. This means rank$(g) = $ rank$(K)$. Now by the minimal $K$-type formula ([Knapp], p.629) we see that in case of compatible orderings every highest weight of $M$ is the restriction of a highest weight of $K$. To get the $K$-separability from this one proceeds by induction on the highest weights of $M$ with respect to the root order we have chosen. For the trivial weight the assertion is clear. Now let $\lambda \in \mathfrak{t}_M^*$ be a highest weight and assume the claim proven for all highest weights $\lambda' < \lambda$. Denote the $M$-representation with highest weight $\lambda$ by $R_M(\lambda)$. Fix a highest weight $\lambda_K$ of $K$ which restricts to $\lambda$ and write $R_K(\lambda_K)$ for the corresponding $K$-representation. Now let $\lambda'$ be a highest weight of $M$ such that $R_M(\lambda')$ occurs in the restriction of $R_K(\lambda_K)$. Since $\lambda'$ is the restriction of a highest $K$-weight which is smaller in the $K$-order than $\lambda_K$, it follows $\lambda' < \lambda$ and so we have

$$R_K(\lambda_K)|_M = R_M(\lambda) + \text{lower terms}$$

and by the induction hypothesis the lower terms already lie in Rep$(K)$. The claim now follows.

Finally to prove the assertion about the constituents of $\wedge^p \mathfrak{p}$ we only have to worry about the case $K = SO(n)$ and $M = SO(n-1)$. Let $\tau_p$ denote the representation of $K$ on $\wedge^p \mathfrak{p}$, then $\tau_1$ is just the natural representation of $K = SO(n)$ on $\mathbb{C}^n$. Restricted to M this splits into $\mathbb{C} \oplus \mathbb{C}^{n-1}$, where $M$ acts trivially on the first summand and via the natural representation $\sigma_1$ on the second. Therefore we get $\tau_p|_M = \wedge^p \sigma_1 \oplus \wedge^{p-1} \sigma_1$ and thus $\wedge^p \sigma_1 = \oplus_{k=0}^p (-1)^{k+p} \tau_p|_M$ in Rep$(M)$. ∎

## 4.3  Selberg's Zeta for Even Dimensions

### 4.3.1   The Trace Formula for M-Types

From now on we will assume dim $X$ even. Recall that the torus $\mathfrak{t}_M \oplus \mathfrak{a}$ admits a real root $\alpha$, unique up to sign. To fit conveniently into normalizations of Harish-Chandra we will work with the form $B$, a scalar multiple of the Killing form such that $B(\alpha) = 4$ From [Wall] we take the following version of the trace formula. Let $\Gamma$ denote a cocompact torsionfree lattice in $G$. Since dim $X$ is even $G$ admits a compact Cartan subgroup and the Euler-Poincaré measure $\mu_{EP}$ on G is nonzero. Thus $\mu_{EP}(\Gamma \backslash G) = (-1)^m \chi(\Gamma \backslash X) = (-1)^m \chi(\Gamma)$, the Euler-Poincaré characteristic. Let $(\varphi, V_\varphi)$ denote a finite dimensional unitary representation of the group $\Gamma$. Consider the

Hilbert space

$$L^2(\Gamma\backslash G, \varphi) = \left\{ f : G \to V_\varphi \mid f(\gamma g) = \varphi(\gamma) f(g), \int_{\Gamma\backslash G} \parallel f(x) \parallel^2 d\mu_{EP}(x) < \infty \right\}$$

(modulo nullfunctions). The group $G$ acts unitarily on $L^2(\Gamma\backslash G, \varphi)$ by right translations and it is known that we have a direct sum decomposition

$$L^2(\Gamma\backslash G, \varphi) = \bigoplus_{\pi \in \hat{G}} N_{\Gamma,\varphi}(\pi)\pi$$

with finite multiplicities $N_{\Gamma,\varphi}(\pi)$ for all $\pi \in \hat{G}$. Let $C$ denote the Casimir operator of $G$.

Since $C$ is central it acts on the representation space of $\pi \in \hat{G}$ as a scalar $\pi(C)$ times the identity. We have the well known formula

$$\pi(C) = B(\lambda_\pi) - B(\rho).$$

Let $\tau \in \hat{K}$ and denote by $\chi_\tau$ the character of $\tau$. By [BM] we may insert the heat kernel $e^{tC} *_K \chi_\tau$ into the trace formula. We introduce some notation: For an irreducible unitary representation $\pi$ of G let $[\pi : \tau] = [\pi \mid_K : \tau]$ denote the multiplicity of $\tau$ in the restriction $\pi \mid_K$. Let $d_\pi$ denote the formal degree of $\pi$ if $\pi$ is a discrete series representation and $d_\pi = 0$ otherwise. This can be summarized in saying $d_\pi = \mu_{\hat{E}P}(\{\pi\})$, where $\mu_{\hat{E}P}$ is the Plancherel measure corresponding to $\mu_{EP}$.

For $\pi \in \hat{G}$ let

$$\lambda_\pi = N_{\Gamma,\varphi}(\pi) - (-1)^m \chi(X_\Gamma) d_\pi \dim(\varphi).$$

Note that, as a consequence of the $L^2$−index theorem of Atiyah [AtSch], the number $\chi(X_\Gamma) d_\pi = \mu_{EP}(\Gamma\backslash G) d_\pi$ is an integer. For $\xi \in \hat{M}$ let $\wedge_\xi$ denote its infinitesimal character and $P_\xi$ the Plancherel density on the principal series attached to $\xi$.

Let

$$f_\xi(t) = \int_{-\infty}^{\infty} e^{-tv^2} P_\xi(v) dv.$$

The conjugacy classes $[\gamma]$ of $\Gamma$ are in bijection to the free homotopy classes of closed paths in $X_\Gamma$. Denote by $l_\gamma$ the minimal length of paths in the class $[\gamma]$ and by $[\gamma_0]$ the underlying primitive class. The letter $\rho$ denotes the half sum of positive roots and $\rho_0 = \rho \mid_\mathfrak{a}$. Let $|\rho_0|$ denote the length of $\rho_0$ under the Killing form. Since $\Gamma$ is torsionfree, any $\gamma \in \Gamma$ is conjugate in $G$ to an element $a_\gamma m_\gamma \in AM$. Now let

$$T_\gamma = \frac{l_{\gamma_0} e^{-|\rho_0| l_\gamma} tr\varphi(\gamma)}{\det(1 - \gamma^{-1} \mid \mathfrak{n})}.$$

The trace formula yields

$$\sum_{\pi \in \hat{G}} [\pi : \tau] \lambda_\pi e^{t(\pi(C)+B(\rho))}$$

$$= \sum_{\xi \in \hat{M}} [\tau : \xi] e^{tB(\wedge_\xi)} \left\{ (-1)^m \chi(X_\Gamma) \dim(\varphi) f_\xi(t) + \sum_{[\gamma] \neq 1} T_\gamma tr\xi(m_\gamma) \frac{e^{-l_\gamma^2/4t}}{\sqrt{4\pi t}} \right\}.$$

Since we have restricted to the even dimensions there is, by Proposition 4.2.1 to every $\sigma \in \hat{M}$ an element $r = \sum_\tau a_\tau \tau$ in Rep(K) such that

(2.3.1) $$\sigma = r \mid_M = \sum_{\tau \in \hat{K}} a_\tau \tau \mid_M .$$

For $\pi \in \hat{G}$ and $r \in Rep(K)$ define

$$N_r(\pi) = \sum_{\tau \in \hat{K}} a_\tau [\pi : \tau].$$

Note that if $\pi = \pi_{\xi,\nu}$ is a principal series representation and $r \in Rep(K)$ with (2.3.1) we have

$$N_r(\pi) = \left\{ \begin{array}{ll} 1 & if \ \xi = \sigma \\ 0 & else. \end{array} \right.$$

Together with the trace formula this leads to the

**Observation:** Let r be in the kernel of the restriction map $Rep(K) \to Rep(M)$ then for every $z \in \mathbb{C}$ we have

$$\sum_{\pi \in \hat{G}, \pi(C) = z} N_r(\pi) \bigg( N_\Gamma(\pi) - \mu_{EP}(\Gamma \backslash G) d_\pi \bigg) = 0.$$

Apply this to the case $G = SL_2(\mathbb{R})$. Here $K = SO(2)$ and $M = \{\pm 1\}$, so there is a lot of $r$ as above. It follows for example for the Mock representations $\pi_\pm$ (see[Lang]) that we have $N_\Gamma(\pi_+) = N_\Gamma(\pi_-)$. Further, the discrete series representations of $G = SL_2(\mathbb{R})$ are parametrized by $\mathbb{Z} - \{0\}$, so $\hat{G}_d = \{\pi_n \mid n \neq 0\}$. From the observation we get $N_\Gamma(\pi_n) = d_{\pi_n} vol(\Gamma \backslash G)$ for $n \neq \pm 1$ and $N_\Gamma(\pi_{\pm 1}) = d_{\pi_{\pm 1}} vol(\Gamma \backslash G) + 1$ (compare [Langl]).

Now for $\sigma \in \hat{M}$ fix any $r \in Rep(K)$ which restricts to $\sigma$. Set $N_\sigma(\pi) = N_r(\pi)$ then the following expression is indeed independent of $r$:

$$\sum_{\pi \in \hat{G}} N_\sigma(\pi) \lambda_\pi e^{t(\pi(C) + B(\rho) - B(\wedge_\sigma))}$$

$$= (-1)^m \chi(X_\Gamma) \dim(\varphi) f_\sigma(t) + \sum_{[\gamma] \neq 1} T_{\gamma,\sigma} \frac{e^{-l_\gamma^2/4t}}{\sqrt{4\pi t}}.$$

### 4.3.2    The zeta function

Fix a minimal parabolic subgroup $P$ of $G$, let $P = MAN$ be its Langlands decomposition, so $A$ is a one dimensional split torus, $M$ its centralizer in $K$ and $N$ is the unipotent radical of $P$. Denote the corresponding Lie algebras by $\mathfrak{p}_0, \mathfrak{m}_0, \mathfrak{a}_0, \mathfrak{n}_0$ and their complexifications by $\mathfrak{p}, \mathfrak{m}, \mathfrak{a}, \mathfrak{n}$. For $N \geq 0$ let $S^N(\mathfrak{n})$ denote the $N$-th symmetric power of $\mathfrak{n}$. Let $\Gamma \subset G$ denote as before a torsion free uniform lattice in $G$. Since $\Gamma$ is torsion free and $G$ of real rank one, every $\gamma \in \Gamma$ is in G conjugate to an element $m_\gamma a_\gamma$ of $MA$.

As before, $(\varphi, V_\varphi)$ denotes a finite dimensional unitary representation of $\Gamma$. We further fix a representation $(\sigma, V_\sigma) \in \hat{M}$ and define the **generalized Selberg Zeta Function** as

$$Z_{\sigma,\varphi}(s) = \prod_{[\gamma] \text{ prime}} \prod_{N \geq 0} \det(1 - e^{-sl_\gamma}\gamma \mid V_N),$$

where $\gamma \in \Gamma$ induces an operator, well defined up to conjugacy on $V_N = V_\varphi \otimes V_\sigma \otimes S^N(\mathfrak{n})$ via $\gamma \mapsto \varphi(\gamma) \otimes \sigma(m_\gamma) \otimes Ad^N((m_\gamma a_\gamma)^{-1})$. Recall that $l_\gamma$ is the minimal length of a closed path freely homotopic to $\gamma$ in $X_\Gamma$. Now in the rank one situation which we are in each free homotopy class of closed paths contains a unique geodesic. So the product above can also be interpreted as a product over all primitive closed geodesics.

For $\pi \in \hat{G}$ and $z \in \mathbb{C}$ let

$$m_z^\sigma = \sum_{\pi \in \hat{G}, \ \pi(C) = -z} N_\sigma(\pi)\lambda_\pi.$$

We consider the virtual operator $\triangle_{\sigma,\varphi}$ defined by the divisor $(z, m_z^\sigma)_{z \in I}$, where $I$ is the set of all $z$ with $m_z^\sigma \neq 0$. Since there are only finitely many $\pi \in \hat{G}$ with $N_\sigma(\pi)d_\pi \neq 0$, we conclude that $\triangle_{\sigma,\varphi}$ is a virtual differential operator plus a finite dimensional one. More explicitly, let $r = \sum_\tau a_\tau \tau$ be in Rep(K) such that $r$ restricts to $\sigma$. Every $\tau \in \hat{K}$ defines a homogeneous vector bundle $E_\tau$. The Casimir acts on the sections of $E_\tau$ as a generalized Laplacian $\triangle_\tau$, so we define $\triangle'_{\sigma,\varphi} = \sum_\tau a_\tau \triangle_\tau$, then $\triangle_{\sigma,\varphi} = \triangle'_{\sigma,\varphi} +$ finite operator. It follows

**Lemma. 4.3.1** *The operator $\triangle_{\sigma,\varphi}$ is zeta admissible. Moreover, we have the asymptotics*

$$tre^{-t\triangle_{\sigma,\varphi}} \sim \sum_{k=-m}^{\infty} c_k t^k. \blacksquare$$

Now set

$$M_\sigma(z,\lambda) = \int_0^\infty t^{z-1} e^{-\lambda t} tr(e^{-t\triangle_{\sigma,\varphi}}) dt.$$

By taking Mellin transforms of all terms in the trace formula we get

$$M_\sigma(z, \lambda + B(\rho) - B(\wedge_\sigma)) = M_1(z,\lambda) + M_{geo}(z,\lambda),$$

now a calculation similar to [D-Hitors, 3.12] shows

$$\exp(M_{geo}(0, s^2)) = Z_{\sigma,\varphi}(s + |\rho_0|)$$

and

$$M_1(z,\lambda) = (-1)^m \chi(X_\Gamma) \dim(\varphi) \int_0^\infty \int_{-\infty}^{\infty} e^{-t(v^2+\lambda)} P_\sigma(v) t^{z-1} dv dt.$$

This will give the analytic continuation to the entire plane. Because of

$$\exp(M_\sigma(0, \lambda)) = \det(\triangle_{\sigma, \varphi} + \lambda)$$

it only remains to consider $M_1(z, \lambda)$.

In order to compute $M_1(z, \lambda)$ we write down the Plancherel formula in our normalization of measures. For this purpose let $\hat{G}_d$ denote the set of discrete series representations of $G$ and for $\pi \in \hat{G}_d$ we get the formal degree as

$$d_\pi = \frac{\mid W(T, K) \mid}{\mid W(T_\mathbb{C}, G_\mathbb{C}) \mid} \prod_{\alpha > 0} \frac{(\alpha, \lambda)}{(\alpha, \rho)},$$

where $\lambda$ is the Harish-Chandra parameter of $\pi$ [Knapp]. To derive this, we take the formula for the formal degree provided by [AtSch]. Recall the notion of the **dual symmetric space** $X^d$ to $X$: Let $G^d$ denote a maximal compact subgroup of the complexification $G_\mathbb{C}$ of $G$ such that $K \subset G^d$. Then $X^d$ is by definition equal to $G^d/K$. The Haar-measure used by [AtSch] is normalized to give the dual space $X^d$ of X the volume 1. By means of the Hopf-Samelson formula one gets that in our normalization the dual space has the volume

$$vol(X^d) = \chi(X^d) = \frac{\mid W(G^d) \mid}{\mid W(K) \mid}.$$

From this the assertion follows. To $\sigma \in \hat{M}$ let $\wedge_\sigma$ denote the corresponding infinitesimal character, i.e. $\wedge_\sigma = w_\sigma + \rho_M$ where $w_\sigma$ is the highest weight.

There is a distinguished element $\gamma_0$ in the center of $M$ [Knapp,p.487] of order two. The set $\hat{M}$ can be devided into an even and an odd part depending on whether the eigenvalue of $\gamma_0$ is $\pm \gamma_0^\rho$. So $\sigma \in \hat{M}_{ev} \Leftrightarrow \sigma(\gamma_0) = \gamma_0^\rho$. Now, from [Knapp] we take the Plancherel formula: Suppose $f \in C_c^\infty(G)$, then we have

$$
\begin{aligned}
f(1) \quad &= \sum_{\pi \in \hat{G}} d_\pi tr(\pi(f)) + i \frac{(-1)^m}{4} \frac{\mid W(K) \mid}{\mid W(G_\mathbb{C}) \mid} \times \\
&\left\{ \sum_{\sigma \in \hat{M}_{ev}} \int_{-\infty}^{\infty} tr \pi_{\sigma i \nu}(f) \prod_{\alpha > 0} \frac{< \wedge_\sigma + i\nu\nu_0, \alpha >}{< \rho, \alpha >} \coth(\pi\nu/2) d\nu \right. , \\
&\left. \sum_{\sigma \in \hat{M}_{odd}} \int_{-\infty}^{\infty} tr \pi_{\sigma i \nu}(f) \prod_{\alpha > 0} \frac{< \wedge_\sigma + i\nu\nu_0, \alpha >}{< \rho, \alpha >} \tanh(\pi\nu/2) d\nu \right\}
\end{aligned}
$$

where $\nu_0 \in \mathfrak{a}^*$ is such that $B(\nu_0, \alpha_0) = B(\alpha_0, \alpha_0)/2$ for the unique positive real root $\alpha_0$ of $\mathfrak{a} + \mathfrak{t}_m$. This means that the Plancherel density $P_\sigma$ is given by

$$P_\sigma(\nu) = i \frac{(-1)^m}{4} \frac{\mid W(K) \mid}{\mid W(G_\mathbb{C}) \mid} \prod_{\alpha > 0} \frac{< \wedge_\sigma + i\nu\nu_0, \alpha >}{< \rho, \alpha >} \coth(\pi\nu/2)$$

if $\sigma$ is even and

$$P_\sigma(\nu) = i\frac{(-1)^m}{4}\frac{\mid W(K)\mid}{\mid W(G_\mathbb{C})\mid}\prod_{\alpha>0}\frac{<\wedge_\sigma + i\nu\nu_0,\alpha>}{<\rho,\alpha>}\tanh(\pi\nu/2)$$

if $\sigma$ is odd. Note that $P_\sigma$ is always a nonnegative function.

The term $p_{\leq k}$ will from now on stand for a polynomial of degree $\leq k$. The expressions 4.3 and 3.4 give

$$\begin{aligned}\log Z_{\sigma,\varphi}(s+\mid\rho_0\mid) &= M_{geo}(0,s^2)\\ &= (-1)^{m+1}\int^{m+1}M_\sigma(m+1,s^2+B(\rho)-B(\wedge_\sigma))d(s^2)\\ &+(-1)^m\int^{m+1}M_1(m+1,s^2)d(s^2)+p_{\leq m}(s^2),\end{aligned}$$

where $\int^k$ means integrating $k$-times (well defined up to $p_{\leq k-1}$). Thus

$$\begin{aligned}\log(Z_{\sigma,\varphi}(s+\mid\rho_0\mid)) &= \text{logdet}(\triangle_{\sigma,\varphi}+s^2-B(\rho)+B(\wedge_\sigma))\\ &-\dim\varphi\,\text{logdet}^{(2)}(\triangle_{\sigma,\varphi}+s^2-B(\rho)+B(\wedge_\sigma))\\ \\ &= \text{logdet}(\triangle_{\sigma,\varphi}+s^2-B(\rho)+B(\wedge_\sigma))\\ &+(-1)^m\int^{m+1}M_1(m+1,s^2)d(s^2)+p_{\leq m}(s^2).\end{aligned}$$

Further we have

$$M_1(m+1,s^2)=m!(-1)^m\chi(X_\Gamma)\int_{-\infty}^\infty\frac{P_\sigma(\nu)}{(\nu^2+s^2)^{m+1}}d\nu.$$

For $\lambda>0$ and $n>k\geq 0$ define

$$F_n^k(\lambda)=(-1)^{n+1}\Gamma(n)\int_{-\infty}^\infty\frac{r^{2k-1}\tanh(\frac{\pi r}{2})}{(r^2+\lambda)^n}dr$$

and

$$F_0^k(\lambda)=\int^{k+1}F_{k+1}^k(\lambda)d\lambda\quad mod\ p_{\leq k}(\lambda).$$

Further we set

$$G_n^k(t)=(-1)^{n+1}\Gamma(n)\int_{-\infty}^\infty\frac{r^{2k-1}\tanh(\frac{\pi r}{2})}{(r^2t+1)^n}dr.$$

Then

$$\begin{array}{llll}\frac{\partial}{\partial\lambda}F_n^k(\lambda) &= F_{n+1}^k(\lambda) & \frac{\partial}{\partial t}G_n^k(t &= G_{n+1}^{k+1}(t)\\ \\ F_n^k(\lambda) &= \lambda^{-n}G_n^k(\frac{1}{\lambda}) & G_n^k(t) &= t^{-n}F_n^k(\frac{1}{t})\end{array}$$

The residue calculus gives

$$F_2^1(\lambda) = -\frac{2}{\sqrt{\lambda}} \sum_{n=0}^{\infty} \frac{1}{(2n+1+\sqrt{\lambda})^2}.$$

For $j \geq 0$ let $D_j$ denote the operator with divisor

$$div D_j = (2n+1, (2n+1)^j)_{n \geq 0},$$

then its zeta function is given by

$$\zeta_{D_j+a}(s) = \sum_{n=0}^{\infty} \frac{(2n+1)^j}{(2n+1+a)^s} \quad Res >> 0.$$

We define $\hat{\zeta}_{D_j+a}(n) = (-1)^{n+1}\Gamma(n)\zeta_{D_j+a}(n)$ and get

$$\frac{\partial}{\partial a}\hat{\zeta}_{D_j+a}(n) = \hat{\zeta}_{D_j+a}(n+1).$$

We further let $\hat{\zeta}_{D_j+a}(0) = \int^{j+2} \zeta_{\hat{D_j}+a}(j+2)da + p_{\leq j+1}$.

**Proposition 4.3.1**
$$F_0^k(a^2) = 8(-1)^k\hat{\zeta}_{D_{2k-1}+a}(0).$$

**Proof:** We first derive the equation:

$$F_0^{k+1} = \int F_0^k(\lambda)d\lambda - \lambda F_0^k(\lambda),$$

this is done as follows:

$$
\begin{aligned}
F_0^{k+1}(\lambda) \quad &= \int^{k+2} F_{k+2}^{k+1}(\lambda)d\lambda \\
&= \int^{k+2} \lambda^{-k-2}G_{k+2}^{k+1}(\frac{1}{\lambda})d\lambda \\
&= -\int^{k+2} \lambda^{-k}\frac{\partial}{\partial\lambda}(G_{k+1}^k(\frac{1}{\lambda}))d\lambda \\
&= -\int^{k+2} \lambda^{-k}\frac{\partial}{\partial\lambda}(\lambda^{k+1}F_{k+1}^k(\lambda))d\lambda \\
&= -(k+1)\int^{k+2} F_{k+2}^k(\lambda)d\lambda - \int^{k+2} \lambda F_{k+2}^k(\lambda)d\lambda.
\end{aligned}
$$

Iterated partial integration gives the claim.

Now we prove the proposition by induction. The case k=1 is easy. The induction step amounts to

$$\hat{\zeta}_{D_{2k+1}+a}(1) = a^2\hat{\zeta}_{D_{2k-1}+a}(1) + p_{\leq 2k+1},$$

which is shown by differentiating $2k + 2$ times. ∎

Now let $E_j$ denote the operator attached to the divisor

$$\mathrm{div}E_j = (2n, (2n)^j)_{n\in\mathbb{N}}.$$

Define

$$H_n^k(\lambda) = (-1)^{n+1}\Gamma(n)\int_{-\infty}^{\infty}\frac{r^{2k-1}\coth(\frac{\pi r}{2})}{(r^2+\lambda)^n}dr$$

and

$$H_0^k(\lambda) = \int^{k+1}H_{k+1}^k(\lambda)d\lambda \quad \mathrm{mod}\ p_{\leq k}(\lambda).$$

We get in perfect analogy to the above

**Proposition 4.3.2**

$$H_0^k(a^2) = 8(-1)^{k+1}\hat{\zeta}_{E_{2k-1}+a}(0).$$

For the proof we only remark that residue calculus gives

$$H_2^1(a^2) = -\frac{1}{a^3} - \frac{2}{a}\sum_{n=1}^{\infty}\frac{1}{(2n+a)^2}.$$

Let $Q_\sigma$ denote the polynomial over $\mathbb{Q}$

$$Q_\sigma(x) = \prod_{\alpha>0}\frac{<\wedge_\sigma + x\nu_0, \alpha>}{<\rho, \alpha>}.$$

**Lemma. 4.3.2** *We have $Q_\sigma(\mathbb{Z}) \subset \mathbb{Z}$ and there is a $k \in \mathbb{N}$ such that $Q_\sigma(\mathbb{N} + k) \subset \mathbb{N}$.*

**Proof:** The existence of such a k follows from the Weyl dimension formula, since for n large enough $Q_\sigma(n)$ is the dimension of a representation. For the rest let d be a natural number such that $F = dQ_\sigma \in \mathbb{Z}[x]$ the we have $F(\mathbb{N}+k) \subset d\mathbb{N}$ or $F(n+k) \equiv 0(d)$ for $n \in \mathbb{N}$. Since n+k runs through all classes modulo d as n runs through the natural numbers it follows that $F(c) \equiv 0(d)$ for all $c \in \mathbb{Z}$. ∎

As we will see, this lemma gives the analytic continuation of $Z_\sigma$ itself whereas Wakayama [Wak] could only show the continuation of a power of $Z_\sigma$.

Define the operator $D_\sigma$ by its divisor according to the cases:

- if $\sigma$ is even:

$$\mathrm{div}\, D_\sigma = (2n, Q_\sigma(2n))_{n \in \mathbb{N}}$$

- if $\sigma$ is odd:

$$\mathrm{div}\, D_\sigma = (2n-1, Q_\sigma(2n-1))_{n \in \mathbb{N}}.$$

Putting things together we have proven that there is a polynomial of degree $\leq 2m$ such that

$$Z_{\sigma,\varphi}(|\rho_0| + s) = \det\left( \triangle_{\sigma,\varphi} + s^2 - B(\rho) + B(\wedge_\sigma) \right)$$

$$\times \left( \det(D_\sigma + s)\exp(P(s)) \right)^{2(-1)^m \dim(\varphi)\chi(X_\Gamma)/\chi(X^d)}.$$

We finally want to show that $P$ is even. To this end recall that $\log Z_{\sigma,\varphi}(|\rho_0| + s)$ tends to zero as $s \to +\infty$. It is thus possible to compute $P(s)$ by means of the asymptotics according to Proposition 1.1.3. Since s only enters in the $\triangle_{\sigma,\varphi}$-term via $s^2$ there is only need to consider the $D_\sigma$-term. Recall the operators $D_j$ and $E_j$. An easy computation shows

$$e^{-tD_j} = (-1)^j \frac{j!}{2} t^{-(j+1)} + O(1) \qquad as\ t \to 0$$

and the same formula holds for $D_j$ replaced by $E_j$. Since only odd j enter in our operators $D_\sigma$ there are only even powers in the asymptotic of $D_\sigma$, so we get

**Theorem 4.3.1** *Let $c_\sigma = B(\wedge_\sigma) - B(\rho)$. The Selberg Zeta Function $Z_{\sigma,\varphi} = Z$ can be expressed as a product of determinants*

$$Z(|\rho_0| + s) = \det\left( \triangle_{\sigma,\varphi} + s^2 + c_\sigma \right) \det^{(2)}\left( \triangle_{\sigma,\varphi} + s^2 + c_\sigma \right)^{-\dim\varphi}$$

$$= \det\left( \triangle_{\sigma,\varphi} + s^2 + c_\sigma \right)\left( \det(D_\sigma + s)\exp(P(s^2)) \right)^{2(-1)^m \dim(\varphi)\chi(X_\Gamma)/\chi(X^d)}.$$

*where $P$ is a polynomial of degree $\leq m$. It follows that $Z_{\sigma,\varphi}$ extends to a meromorphic function with finitely many poles and that*

$$\frac{Z(|\rho_0| + s)}{Z(|\rho_0| - s)} = \left( \frac{\det(D_\sigma + s)}{\det(D_\sigma - s)} \right)^{2(-1)^m \dim(\varphi)\chi(X_\Gamma)/\chi(X^d)}.$$

*Further $Z$ satisfies a modified Riemann hypothesis, its singularities are devided into the spectral ones and the topological ones. A complex number s is a spectral singularity if and only if the number $-c_\sigma - s^2$ is an eigenvalue of $\triangle_{\sigma,\varphi}$ and its order is its multiplicity then, i.e. spectral singularities are zeroes. The complex number s is a topological singularity if and only if $s = -2n$, $n \in \mathbb{N}$ if $\sigma$ is even or $s = -2n+1$, $n \in \mathbb{N}$ if $\sigma$ is odd. Its order is then equal to the number $(-1)^m Q_\sigma(s) \dim(\varphi)\chi(X_\Gamma)/\chi(X^d)$.*

Finally, we compare our results with that of Juhl [Ju]. Set $\varphi = 1$. We introduce some notation. Let $G^d$ denote the dual group to $G$. This is a compact form of $G$ which we may assume to contain $K$. The compact dual $A^d$ of $A$ is a torus in $G^d$. The Lie algebra $\mathfrak{n}$ defines a complex structure on the compact manifold $G^d/MA^d$ by

$$T^{(1,0)}(G^d/MA^d) \cong G^d \times_{MA^d} \mathfrak{n}.$$

For a character $\chi^d$ of A with $\chi^d(\gamma_0)\sigma(\gamma_0) = 1$, where $\gamma_0$ is defined as in 4.4 let $ind(\mathfrak{n}, \sigma, \chi^d)$ denote the index of the Dolbeault complex with values in the vector bundle $E_{\sigma \otimes \chi^d}$ defined by $\sigma \otimes \chi^d$ over $G^d/MA^d$. Let $H_0 \in \mathfrak{a}$ as before defined by $\alpha(H_0) = 1$ for the short root $\alpha$ of A. Set

$$Z^+_{(\mathfrak{n},\sigma)}(s) = \prod_{\substack{\chi^d \in \hat{A}^d \\ \chi^d(H_0) > \rho_0 \\ \chi^d(\gamma_0)\sigma(\gamma_0) = 1}} (s + \chi^d(H_0))^{\operatorname{ind}(\mathfrak{n},\sigma,\chi^d)}.$$

where the product denotes $\zeta$-regularized product.

Let $X^d$ denote the compact dual space to X. A. Juhl has shown [Ju, p. 184]:

$$\frac{Z_{\sigma,\varphi}(|\rho_0| + s)}{Z_{\sigma,\varphi}(|\rho_0| - s)} \cong \left( \frac{Z^+_{(\mathfrak{n},\sigma)}(|\rho_0| + s)}{Z^+_{(\mathfrak{n},\sigma)}(|\rho_0| - s)} \right)^{-2\chi(X_\Gamma)/\chi(X^d)}.$$

where $\cong$ means equality up to the exponential of a polynomial. Since $Z_{(\mathfrak{n},\sigma)}$ is holomorphic and nonzero for $Res > 0$ we conclude

$$Z^+_{(\mathfrak{n},\sigma)}(s) = \det(D_\sigma + s)^{(-1)^{m+1}},$$

it further follows that we have equality above instead of "$\cong$". Let $A^d$ denote the compact dual of A acting on the compact form $G^d$ then we easily get

$$\operatorname{ind}(\mathfrak{n}, \sigma, \chi^d) = (-1)^{m+1} \prod_{\alpha > 0} \frac{< \wedge_\sigma + \chi^d + \rho \mid_{\mathfrak{a}}, \alpha >}{< \rho, \alpha >}.$$

## 4.4 Ruelle's Zeta for Even Dimensions

We extend the definition of the zeta function $Z_{\sigma,\varphi}$ to finite dimensional unitary representations $\sigma$ which are not necessarily irreducible in the obvious way: If $\sigma = \oplus_{\xi \in \hat{M}} n_\xi \xi$ is the isotypical decomposition of $\sigma$ then

$$Z_{\sigma,\varphi}(s) \stackrel{\text{def}}{=} \prod_{\xi \in \hat{M}} (Z_{\xi,\varphi}(s))^{n_\xi}.$$

Let $\mathfrak{a}$ denote the complex Lie algebra of $A$, so that $M$ is the centralizer of $\mathfrak{a}$ in $K$. Let $P = MAN$ be a minimal parabolic. Since we have normalized the form $B$ such that the real root $\alpha_r$ has length 2 the possible root length of $(\mathfrak{g}, \mathfrak{a})$ are 1 and 2. Let $\mathfrak{n} = \mathfrak{n}_1 \oplus \mathfrak{n}_2$, where $A$ acts according to $\alpha_r$ on $\mathfrak{n}_2$ and according to $\frac{1}{2}\alpha_r$ on $\mathfrak{n}_1$. Let $\sigma_1, \sigma_2$ be the representations of $M$ on $\mathfrak{n}_1, \mathfrak{n}_2$ and write $\tau_p$ for the $K$-representation on $\wedge^p \mathfrak{p}$. Then we have

$$\tau_1|_M = \sigma_0 \oplus \sigma_1 \oplus \sigma_2,$$

where $\sigma_0$ is the trivial one dimensional representation.

**Theorem 4.4.1** *For the Ruelle zeta function*

$$R_\varphi(s) = \prod_{[\gamma] \text{prime}} \det(1 - e^{-sl_\gamma} \varphi(\gamma))$$

*we have*

$$R_\varphi(s) = \prod_{l=0}^{\dim \mathfrak{n}_1} \prod_{k=0}^{\dim \mathfrak{n}_2} Z_{\wedge^l \sigma_1 \otimes \wedge^k \sigma_2, \varphi}(s + l + 2k)^{(-1)^{l+k}}.$$

**Proof:** For an arbitrary representation $\sigma$ of $M$ we write

$$
\begin{aligned}
\log Z_{\sigma, \varphi}(s) &= \sum_{[\gamma] \text{prime}} \sum_{N \geq 0} \text{tr} \log(1 - e^{-sl_\gamma} \gamma | V_N) \\
&= -\sum_{[\gamma] \text{prime}} \sum_{N \geq 0} \sum_{n=1}^{\infty} \frac{1}{n} e^{-snl_\gamma} \text{tr} \varphi(\gamma^n) \text{tr} \sigma(m_\gamma^n) \text{tr} \text{Ad}^N((m_\gamma a_\gamma)^{-n}) \\
&= -\sum_{[\gamma]} \text{tr} \varphi(\gamma) \frac{e^{-sl_\gamma}}{\mu_\gamma} \frac{\text{tr} \sigma(m_\gamma)}{\det(1 - \gamma^{-1} | \mathfrak{n})}.
\end{aligned}
$$

Now we have

$$
\begin{aligned}
\det(1 - \gamma^{-1} | \mathfrak{n}) &= \det(1 - \gamma^{-1} | \mathfrak{n}_1) \det(1 - \gamma^{-1} | \mathfrak{n}_2) \\
&= \sum_{l=0}^{\dim \mathfrak{n}_1} \sum_{k=0}^{\dim \mathfrak{n}_2} (-1)^{k+l} \text{tr}(\gamma^{-1} | \wedge^l \mathfrak{n}_1 \otimes \wedge^k \mathfrak{n}_2) \\
&= \sum_{l=0}^{\dim \mathfrak{n}_1} \sum_{k=0}^{\dim \mathfrak{n}_2} (-1)^{k+l} e^{-(l+2k)l_\gamma} \text{tr}(m_\gamma | \wedge^l \mathfrak{n}_1 \otimes \wedge^k \mathfrak{n}_2),
\end{aligned}
$$

from which the proposition follows. ∎

## 4.5   Selberg's Zeta Function for Odd Dimensions

Here we assume $n = \dim X$ odd. It follows $X = \mathcal{H}_n(\mathbb{R})$ the real hyperbolic space. We have $G = SO(n,1)^+$, $K = SO(n)$ and $M = SO(n-1)$. We will normalize the form $B$ in such a way that $X$ has curvature $-1$. This means

$$B = \frac{1}{n-1} B_{\text{Killing}},$$

where $B_{\text{Killing}}$ denotes the Killing form. The Haar measure $\mu = \mu_B$ will be the standard measure attached to $B$.

We will follow the notation of 4.3. Since there is no compact Cartan subgroup there is no discrete series, so we needn't consider the numbers $\lambda_\pi$. Indeed, let $\sigma \in \hat{M}$ be $K$-separable, choose $r \in \text{Rep}(K)$ which restricts to $\sigma$ and let $N_\sigma(\pi) = N_r(\pi)$. We then get

$$\sum_{\pi \in \hat{G}} N_\sigma(\pi) N_{\Gamma,\varphi}(\pi) e^{t(\pi(C) + B(\rho) - B(\wedge_\sigma))} = \mu(\Gamma \backslash G) \dim(\varphi) f_\sigma(t) + \sum_{[\gamma] \neq 1} T_{\gamma,\sigma} \frac{e^{-l_\gamma^2/4t}}{\sqrt{4\pi t}}.$$

Here again $f_\sigma(t) = \int_{-\infty}^{\infty} e^{-tv^2} P_\sigma(v) dv$ but now the Plancherel density $P_\sigma$ is an even polynomial of degree $n-1$ ([Knapp], p.485). This makes the factor at infinity much easier to handle. Let

$$\zeta_a(s) \overset{\text{def}}{=} \frac{1}{\Gamma(s)} \int_0^\infty t^{s-1} f_\sigma(t) e^{-ta^2} dt.$$

For $a > 0$ fixed this extends to a meromorphic function, holomorphic at $s = 0$ and with $Z_1(a) \overset{\text{def}}{=} \exp(-\zeta_a'(0))$ we have

**Lemma. 4.5.1**

$$Z_1(a) = \exp(2\pi \int_0^a P_\sigma(iy) dy)$$

**Proof:** ([Fried], p.533). ∎

Denote the polynomial $2\pi \int_0^a P_\sigma(iy) dy$ by $\tilde{P}_\sigma(a)$. Along the lines of section 4.3 we get

**Theorem 4.5.1** *Let $X$ be a real hyperbolic space of odd dimension. Fix a $K$-separable representation $\sigma \in \hat{M}$. Consider the zeta function*

$$Z_{\sigma,\varphi}(s) = \prod_{[\gamma] \text{prime}} \prod_{N \geq 0} \det(1 - e^{-sl_\gamma} \gamma | V_N),$$

*where $V_N = V_\varphi \otimes V_\sigma \otimes S^N(\mathfrak{n})$ and $\gamma$ operates as $\varphi(\gamma) \otimes \sigma(M_\sigma) \otimes \mathrm{Ad}^N((m_\gamma a_\gamma)^{-1})$. Then $Z_{\sigma,\varphi}$ satisfies the equation*

$$Z_{\sigma,\varphi} = \det(\triangle_{\sigma,\varphi} + s^2 + c_\sigma) \exp(\dim \varphi \mu(\Gamma \backslash G) \tilde{P}_\sigma(s)).$$

*It follows that $Z_{\sigma,\varphi}$ extends to an entire function and that*

$$\hat{Z}_{\sigma,\varphi} = Z_{\sigma,\varphi}(s) \exp(-\dim \varphi \mu(\gamma \backslash G) \tilde{P}_\sigma(s))$$

*satisfies*

$$\hat{Z}_{\sigma,\varphi}(s) = \hat{Z}_{\sigma,\varphi}(n-1-s).$$

## 4.6   Ruelle's Zeta for Odd Dimensions

We still stick to the situation of the preceding subsection, i.e. $X = \mathcal{H}_n(\mathbb{R})$, $n$ odd. Denote by $\tau_p$ the representation of $K$ on $\wedge^p \mathfrak{p}$ the we have $\tau_p|_M = \sigma_p + \sigma_{p-1}$ for some $\sigma_j \in \hat{M}$. Define the Ruelle zeta function as

$$R_\varphi(s) = \prod_{[\gamma] \mathrm{prime}} \det(1 - e^{-sl_\gamma} \varphi(\gamma))$$

Then we get

**Proposition 4.6.1**  *([Fried],p.532)*

$$R_\varphi(s) = \prod_{p=0}^{n-1} Z_{\sigma_p,\varphi}(s+p)^{(-1)^p}.$$

From this one reads off:

**Theorem 4.6.1**  *([Fried]) The Ruelle zeta function $R_\varphi(s)$ extends to a meromorphic function on $\mathbb{C}$, its zeroes and poles are all in the set $\mathbb{R} \cup (-\frac{n-1}{2} + i\mathbb{R}) \cup (-\frac{n-1}{2} + 1 + i\mathbb{R} \cup \ldots \cup (\frac{n-1}{2} + i\mathbb{R})$ and can be expressed in terms of Laplace eigenvalues for the Laplacians $\triangle_{\varphi,p}$, for $p = 0, 1, \ldots, n$.*

*The function $R_\varphi(s)$ satisfies a functional equation*

$$R_\varphi(s) = R_\varphi(-s) \exp(\dim \varphi \mu(\Gamma \backslash G) Q(s)),$$

*where $Q$ is the polynomial*

$$Q(s) = \sum_{p=0}^{n-1} \tilde{P}_{\sigma_p}(s+p) - \tilde{P}_{\sigma_p}(n-1-s-p).$$

*In the center of the functional equation we have the special value*

$$R_\varphi(0) = \tau(\phi),$$

*the analytic torsion of $\varphi$.* ∎

# Chapter 5

# The Holomorphic Torsion Zeta Function

In this chapter we are going to consider a different generalization of the Selberg Zeta Function. From this we will construct a zeta function which has the holomorphic torsion as special value. Of particular neatness is the case of the $m$-th power $\mathcal{H}^m$, $m \in \mathbb{N}$, of the hyperbolic plane $\mathcal{H}$, where the Euler product looks very simple (sec. 5.2).

## 5.1   The trace of the heat kernel

Let $\bar{X}$ denote a compact locally symmetric space whose universal covering $X$ is hermitian globally symmetric of the noncompact type. Then $X = G/K$ where $G$ is the group of orientation preserving isometries of $X$ and $K$ is a maximal compact subgroup of $G$. Then there is a uniform lattice $\Gamma$ in $G$ such that $\bar{X} = \Gamma \backslash X = \Gamma \backslash G/K$. Further, $\Gamma \cong \pi_1(\bar{X})$, the fundamental group of $\bar{X}$. We will therefore write $X_\Gamma$ instead of $\bar{X}$. It follows that $G$ is a semisimple Lie group without center that admits a compact Cartan subgroup $T \subset K$. We denote the real Lie algebras of $G$, $K$ and $T$ by $\mathfrak{g}_0, \mathfrak{k}_0$ and $\mathfrak{t}_0$ and their complexifications by $\mathfrak{g}, \mathfrak{k}$ and $\mathfrak{t}$. We will denote the Killing form of $\mathfrak{g}$ by B. As well we will write $B$ for the diagonal of the Killing form so $B(X) = B(X, X)$. Denote by $\mathfrak{p}_0$ the orthocomplement of $\mathfrak{k}_0$ in $\mathfrak{g}_0$ with respect to $B$ then via the differential of exp the space $\mathfrak{p}_0$ is isomorphic to the real tangent space of $X = G/K$ at the point $eK$. Let $\Phi(\mathfrak{t}, \mathfrak{g})$ denote the system of roots of $(\mathfrak{t}, \mathfrak{g})$, let $\Phi_c(\mathfrak{t}, \mathfrak{g}) = \Phi(\mathfrak{t}, \mathfrak{k})$ denote the subset of compact roots and $\Phi_{nc} = \Phi - \Phi_c$ the set of noncompact roots. To any root $\alpha$ let $\mathfrak{g}_\alpha$ denote the corresponding root space. Fix an ordering $\Phi^+$ on $\Phi = \Phi(\mathfrak{t}, \mathfrak{g})$ and let $\mathfrak{p}_\pm = \bigoplus_{\alpha \in \Phi_{nc}^+} \mathfrak{g}_{\pm \alpha}$. Then the complexification $\mathfrak{p}$ of $\mathfrak{p}_0$ splits as $\mathfrak{p} = \mathfrak{p}_+ \oplus \mathfrak{p}_-$ and the ordering can be chosen such that this decomposition corresponds via exp to the decomposition of the complexified tangent space of $X$ into holomorphic and antiholomorphic part.

Let $\theta$ denote the Cartan involution on $\mathfrak{g}_0$ and on $G$ corresponding to the choice of $K$. Extend $\theta$ linearly to $\mathfrak{g}$. Let $H$ denote a $\theta$-stable Cartan subgroup of $G$ then $H = AB$ where $A$ is a split torus and $B$ compact. The use of the letter $B$ here will not cause any confusion. The dimension of $A$ is called the split rank of H. Let $\mathfrak{a}$ denote the complex Lie algebra of $A$. Then $\mathfrak{a}$ is an abelian subspace of $\mathfrak{p} = \mathfrak{p}_+ \oplus \mathfrak{p}_-$. Let $X \mapsto X^c$ denote the complex conjugation on $\mathfrak{g}$ according to the real form $\mathfrak{g}_0$. The next lemma shows that $\mathfrak{a}$ lies skew to the decomposition $\mathfrak{p} = \mathfrak{p}_+ \oplus \mathfrak{p}_-$.

**Lemma. 5.1.1** *Let $Pr_\pm$ denote the projections from $\mathfrak{p}$ to $\mathfrak{p}_\pm$ the we have* $\dim Pr_+(\mathfrak{a}) = \dim$ *$Pr_-(\mathfrak{a}) = \dim \mathfrak{h}_R$, or, what amounts to the same:* $\mathfrak{a} \cap \mathfrak{p}_\pm = 0$.

Proof: $\mathfrak{a}$ is stable under complex conjugation which interchanges $\mathfrak{p}_+$ and $\mathfrak{p}_-$. So let $X \in \mathfrak{p}_+$ be such that $X + X^c \in \mathfrak{a}$. Since $\mathfrak{a}$ is abelian, the assumption $i(X - X^c) \in \mathfrak{a}$ would lead to $0 = [X + X^c, i(X - X^c)] = 2i[X^c, X]$ and the latter only vanishes for $X = 0$. So $\mathfrak{a}$ does not contain $X$ or $X^c$ and the claim follows. ∎

Now let $\mathfrak{a}'$ denote the orthocomplement of $\mathfrak{a}$ in $Pr_+(\mathfrak{a}) \oplus Pr_-(\mathfrak{a})$. For later use we write $\mathcal{V} = \mathfrak{a} \oplus \mathfrak{a}' = \mathcal{V}_+ \oplus \mathcal{V}_-$, where $\mathcal{V}_\pm = \mathcal{V} \cap \mathfrak{p}_\pm = Pr_\pm(\mathfrak{a})$.

Let $n = 2m$ denote the real dimension of $X$ and for $0 \le p, q \le m$ let $\Omega^{p,q}(X)$ denote the space of smooth $(p,q)$-forms on $X$. Fix a finite dimensional unitary representation $(\varphi, V_\varphi)$ of the group $\Gamma$. Then $\varphi$ defines an hermitian flat holomorphic vector bundle $E_\varphi = X \times_\Gamma V_\varphi$ over $X_\Gamma$. On the space $\Omega^{p,q}(X_\Gamma, E_\varphi)$ of $E_\varphi$-valued forms we have a Laplacian $\triangle^\varphi_{p,q,\Gamma}$ and by construction it is clear that $\triangle^\varphi_{p,q,\Gamma} = \triangle_{p,q,\Gamma} \otimes 1$ where $\triangle_{p,q,\Gamma}$ is the Laplacian on $\Omega^{p,q}(X_\Gamma)$ and we have used $\Omega^{p,q}(X_\Gamma, E_\varphi) = C^\infty(\wedge^p T^*_{hol}(X_\Gamma) \otimes \wedge^q T^*_{ahol}(X_\Gamma) \otimes E_\varphi)$. Further $\triangle_{p,q,\Gamma}$ is the pushdown of the Laplacian $\triangle_{p,q}$ on $\Omega^{p,q}(X)$.

By [BM] the heat operator $e^{-t\triangle_{p,q}}$ has a smooth kernel $h_t^{p,q}$ of rapid decay in

$$(C^\infty(G) \otimes \text{End}(\wedge^p \mathfrak{p}_+ \otimes \wedge^q \mathfrak{p}_-))^{K \times K}.$$

Now fix p and set for $t > 0$

$$f_t^p = \sum_{q=0}^m q(-1)^q \text{tr } h_t^{p,q},$$

where tr means the trace in $\text{End}(\wedge^p \mathfrak{p}_+ \otimes \wedge^q \mathfrak{p}_-)$.

We want to compute the trace of $f_t^p$ on the principal series representations. To this end let $H = AB$ be a $\theta$-stable Cartan subgroup with split part $A$ and compact part $B$. Let $P$ denote a parabolic subgroup of $G$ with Langlands decomposition $P = MAN$. Note that $M$ is of inner type (see [HC-HA1], Lemma 4.9). Let $(\xi, W_\xi)$ denote an irreducible unitary representation of $M$, $e^\nu$ a quasicharacter of $A$ and set $\pi_{\xi,\nu} = Ind_P^G(\xi \otimes e^{\nu+\rho_P} \otimes 1)$, where $\rho_P$ is the half of the sum of the $P$-positive roots.

**Proposition 5.1.1** *For the trace of $f_t^p$ under $\pi_{\xi,\nu}$ we have*

$$\text{tr}\,\pi_{\xi,\nu}(f_t^p) = \begin{cases} 0 & \text{if } \dim \mathfrak{a} > 1 \\ -e^{t\pi_{\xi,\nu}(C)} \sum_{q=0}^{\dim(\mathfrak{p}_-)-1} (-1)^q \dim(W_\xi \otimes \wedge^p \mathfrak{p}_+ \otimes \wedge^q(\mathfrak{a}^\perp \cap \mathfrak{p}_-))^{K \cap M} & \text{if } \dim \mathfrak{a} = 1. \end{cases}$$

Proof: As before let $\mathfrak{a}'$ the span of all $X - X^c$, where $X + X^c \in \mathfrak{a}$, $X \in \mathfrak{p}_+$ then $\mathcal{V} = \mathfrak{a} \oplus \mathfrak{a}' = \mathcal{V}_+ \oplus \mathcal{V}_-$ where $\mathcal{V}_\pm = \mathcal{V} \cap \mathfrak{p}_\pm$. The group $K_M = K \cap M$ acts trivially on $\mathfrak{a}$ so for $x \in K_M$ we have $X + X^c = \mathrm{Ad}(x)(X + X^c) = \mathrm{Ad}(x)X + \mathrm{Ad}(x)X^c$. Since $K$ respects the decomposition $\mathfrak{p} = \mathfrak{p}_- \oplus \mathfrak{p}_+$, we conclude that $K_M$ acts trivially on $\mathcal{V}$, hence on $\mathcal{V}_-$. Let $r = \dim\mathfrak{a} = \dim\mathcal{V}_-$. As a $K_M$-module we have

$$\wedge^p \mathfrak{p}_- = \sum_{a+b=q} \wedge^a \mathcal{V}_- \otimes \wedge^b \mathcal{V}_-^\perp$$
$$= \sum_{a+b=q} \binom{r}{a} \wedge^b \mathcal{V}_-^\perp,$$

where $\mathcal{V}_-^\perp = \mathfrak{a}^\perp \cap \mathfrak{p}_-$. By Frobenius reciprocity we get

$$\mathrm{tr}\pi_{\xi,\nu}(\mathrm{f}_t^\mathrm{p}) = \mathrm{tr}\pi_{\xi,\nu}(\sum_{q=0}^m q(-1)^q h_t^{p,q})$$
$$= e^{t\pi_{\xi,\nu}(C)} \sum_{q=0}^m q(-1)^q \dim(\mathrm{V}_{\pi_{\xi,\nu}} \otimes \wedge^p \mathfrak{p}_+ \otimes \wedge^q \mathfrak{p}_-)^\mathrm{K}$$
$$= e^{t\pi_{\xi,\nu}(C)} \sum_{q=0}^m q(-1)^q \dim(\mathrm{W}_{\pi_\xi} \otimes \wedge^p \mathfrak{p}_+ \otimes \wedge^q \mathfrak{p}_-)^{\mathrm{K} \cap \mathrm{M}}$$
$$= e^{t\pi_{\xi,\nu}(C)} \sum_{q=0}^m \sum_{a=0}^q q(-1)^q \binom{r}{q-a} \dim(\mathrm{W}_{\pi_\xi} \otimes \wedge^p \mathfrak{p}_+ \otimes \wedge^a \mathcal{V}_-^\perp)^{\mathrm{K} \cap \mathrm{M}}$$
$$= e^{t\pi_{\xi,\nu}(C)} \sum_{a=0}^m \sum_{q=a}^m q(-1)^q \binom{r}{q-a} \dim(\mathrm{W}_{\pi_\xi} \otimes \wedge^p \mathfrak{p}_+ \otimes \wedge^a \mathcal{V}_-^\perp)^{\mathrm{K} \cap \mathrm{M}}$$

By taking into account that we have $a \leq m-r$ we get $\sum_{q=a}^m q(-1)^q \binom{r}{q-a} = (-1)^r \binom{a}{1-r}$ and the claim follows. ∎

In view of the preceding proposition fix a $\Theta$-stable Cartan subgroup $H = AB$ with $\dim A = 1$ and a parabolic $P = MAN$. Fix a system of positive roots compatible with $P$ and let $\rho$ denote the half sum of positive roots. For $\xi \in \hat{M}$ let $\lambda_\xi \in \mathfrak{b}^*$ denote the infinitesimal character of $\xi$. Recall that we have

$$\pi_{\xi,\nu}(C) = B(\nu) + B(\lambda_\xi) - B(\rho).$$

Let $c = c(H)$ denote the number of positive roots in $\phi(\mathfrak{a}, \mathfrak{g})$. Since $\mathfrak{a}$ is a split torus and there is a real root, it follows $c = 1$ or $c = 2$.

**Lemma. 5.1.2** *If $c = 1$ let $\mathfrak{g} = \mathfrak{g}_1 \oplus \mathfrak{g}_2$ where $\mathfrak{g}_1$ is the simple ideal of $\mathfrak{g}$ that contains $\mathfrak{a}$. Accordingly we have $G \cong G_1 \times G_2$. Let $M_1 = M \cap G_1$ then we have an isomorphism as $K_M$-modules:*

$$\mathfrak{p}_\pm \cong \mathbb{C} \oplus \mathfrak{p}_{G_2,\pm} \oplus \mathfrak{p}_{M_1}.$$

*In the case $c = 2$ the space $M/K_M$ inherits the complex structure from $G/K$, i.e. $\mathfrak{p}_M = \mathfrak{p}_{M,+} \oplus \mathfrak{p}_{M,-}$ where $\mathfrak{p}_{M,\pm} = \mathfrak{p}_M \cap \mathfrak{p}_\pm$ and we have a $K_M$-module isomorphism:*

$$\mathfrak{p}_\pm \cong \mathbb{C} \oplus \mathfrak{p}_{M,\pm} \oplus \mathfrak{n}_\pm,$$

*where*

$$n_\pm = \bigoplus_{\substack{\alpha \, \in \, \phi_{nI}(\mathfrak{h}, \mathfrak{g}) \\ \pm\alpha|_\mathfrak{b} > 0}} \mathfrak{g}_\alpha.$$

**Proof:** Let $H$ be a generator of $\mathfrak{a}_0$. Write $\alpha_r$ for the unique positive real root in $\phi(\mathfrak{h}, \mathfrak{g})$. Since for any root $\alpha$ we have $2B(\alpha, \alpha_r)/B(\alpha_r) \in \pm\{0, 1, 2, 3\}$ the only possible roots in $\phi^+(\mathfrak{a}, \mathfrak{g})$ are $\alpha_r/2, \alpha_r, 3\alpha_r/2$. Consider an embedding $\mathfrak{g} \hookrightarrow gl_n$ such that the Cartan involution becomes $\theta(X) = -X^t$ and $\mathfrak{a}$ is mapped to the diagonal. Since $[\mathfrak{g}_{\alpha_r}, \theta(\mathfrak{g}_{\alpha_r})] = \mathfrak{a}$ it is easy to see that $3\alpha_r/2$ does not occur.

Now write $H = Y + Y^c$ with $Y \in \mathfrak{p}_+$. According to the root space decomposition $\mathfrak{g} = \mathfrak{a} \oplus \mathfrak{k}_M \oplus \mathfrak{p}_M \oplus \mathfrak{n} \oplus \theta(\mathfrak{n})$ we write $Y = Y_a + Y_k + Y_p + Y_n + Y_{\theta(n)}$. Because of $\theta(Y) = -Y$ it follows $Y_k = 0$ and $Y_{\theta(n)} = -\theta(Y_n)$. For arbitrary $Z \in \mathfrak{b}$ we have $[Z, Y] = 0$ since $[Z, H] = 0$ and the projection $Pr_+$ is $k_M$-equivariant. Hence $0 = [Z, Y] = [Z, Y_a] + [Z, Y_p] + [Z, Y_n - \theta(Y_n)]$. Since $[Z, Y_a] = 0$ and $[Z, Y_p] \in \mathfrak{p}_M$ it follows $[Z, Y_n - \theta(Y_n)] \in \mathfrak{p}_M \cap (\mathfrak{n} \oplus \theta(n)) = 0$. Therefore $Y_p = 0$ and $Y_n \in \mathfrak{n}_r = \mathfrak{g}_{\alpha_r}$. Now $k_M$ acts trivially on $H$ and thus on $Y$, so on $Y_n$, i.e. $[\mathfrak{k}_M, Y_n] = 0$. Further $\mathfrak{a} \cap \mathfrak{p}_\pm$ is trivial, so $Y_n \neq 0$, so it generates the root space $\mathfrak{g}_{\alpha_r}$, so $[\mathfrak{k}_M, \mathfrak{g}_{\alpha_r}] = 0$.

**Assume** there is a root $\alpha \neq \alpha_r$ such that $\alpha|_\mathfrak{a} = \alpha_r|_\mathfrak{a}$. Then $B(\alpha, \alpha_r) > 0$ and hence $\beta = \alpha - \alpha_r$ is a root. The root $\beta$ is imaginary. Suppose $\beta$ is compact, then $\mathfrak{g}_\beta \subset \mathfrak{k}_M$ and we have $[\mathfrak{g}_{\alpha_r}, \mathfrak{g}_\beta] = \mathfrak{g}_\alpha$, which contradicts $[\mathfrak{k}_M, \mathfrak{g}_{\alpha_r}] = 0$. It follows that $\beta$ is noncompact and thus it follows

$$\mathfrak{g}_{\alpha_r|_\mathfrak{a}} = [\mathfrak{p}_M, \mathfrak{g}_{\alpha_r}] \oplus \mathfrak{g}_{\alpha_r},$$

where $\mathfrak{g}_{\alpha_r|_\mathfrak{a}} = \{X \in \mathfrak{g} | [H, X] = \alpha_r(H) X\}$ is the root space for the restricted root $\alpha_r|_\mathfrak{a}$. Let $\mathfrak{g}_1$ as above the simple ideal containing $\mathfrak{a}$, then $\mathfrak{g}_1$ also contains $\mathfrak{g}_{\alpha_r}$ and hence we may substitute $\mathfrak{p}_{M_1}$ in the above. Assume the root space $\mathfrak{g}_{\alpha_r/2|_\mathfrak{a}}$ is nonzero, then $(\mathfrak{a} \oplus \mathfrak{m} \oplus \mathfrak{n}_r \oplus \mathfrak{g}_{\alpha_r/2|_\mathfrak{a}} \oplus \mathfrak{g}_{-\alpha_r/2|_\mathfrak{a}}) \mathfrak{g}_1$ would be a nontrivial ideal of $\mathfrak{g}_1$, therefore $\alpha_r/2$ does not occur as root of $(\mathfrak{a}, \mathfrak{g})$ and we may state

$$\mathfrak{n} = [\mathfrak{p}_{M_1}, \mathfrak{g}_{\alpha_r}] \oplus \mathfrak{g}_{\alpha_r}.$$

The map $\phi : X \mapsto [X, Y_n]$ has no kernel in $\mathfrak{p}_{M_1}$ since this would similarly allow us to construct a nontrivial ideal of $\mathfrak{g}_1$. The map

$$\psi : \mathfrak{a} \oplus \mathfrak{p}_M \oplus \mathfrak{n} \quad \to \mathfrak{p}$$
$$a + p + n \quad \mapsto \quad a + p + n - \theta(n)$$

is a $K_M$-isomorphism as is the map $\phi$. So, as $K_M$-modules $\mathfrak{p} \cong \mathfrak{a} \oplus \mathfrak{n}_r \oplus \mathfrak{p}_M \oplus \mathfrak{p}_M$. To derive the assertion for the case $c = 1$ it remains to show $\mathfrak{p}_{M_1} \cap \mathfrak{p}_+ = 0$, because then $\mathfrak{p}_{M_1} \cong Pr_+(\mathfrak{p}_{M_1})$. For this assume $X \in \mathfrak{p}_{M_1} \cap \mathfrak{p}_+$, then $0 = [X, Y] = [x, Y_a] + [X, Y_n - \theta(Y_n)]$ and therefore $[X, Y_n] = 0$ which implies $X = 0$.

The remaining case is the case where there exists no such root $\alpha$ as above. This means $\mathfrak{n} = \mathfrak{g}_{\alpha_r/2|_\mathfrak{a}} \oplus n_r$. This gives $[\mathfrak{p}_M, \mathfrak{n}_r] = 0$ and thus $[\mathfrak{p}_M, Y] = 0$, which implies $\mathfrak{p}_M = \mathfrak{p}_{M,+} \oplus \mathfrak{p}_{M,-} = \mathfrak{p}_M \cap \mathfrak{p}_+ \oplus \mathfrak{p}_M \cap \mathfrak{p}_-$. So $\mathfrak{p}_M$ inherits the holomorphic structure and so does $\mathfrak{n}$. Since the maps $\phi$ and $\psi$ above are $K_M$-homomorphisms and $\mathfrak{p}_+$ is given by a choice of positive roots the claim follows. ∎

In this section we will prove a general lemma which we will apply to the group $M$ later. So here we will not assume $G$ connected, but $G$ should be of inner type (see [HC-HA1] Lemma 4.9).

Let $(\tau, V_\tau)$ denote an irreducible unitary representation of $K$. Since $G$ is of inner type, the Casimir operator $C_K$ of $K$ will acts as a scalar $\tau(C_K)$ on $V_\tau$.

**Lemma. 5.1.3** *Let $(\pi, W_\pi)$ be an irreducible unitary representation of $G$ and assume that*

$$\sum_{p=0}^{\dim \mathfrak{p}_-} (-1)^p \dim(W_\pi \otimes \wedge^p \mathfrak{p}_- \otimes V_\tau)^K \neq 0$$

*then we have*

$$\pi(C) = \tau(C_K) - B(\rho) + B(\rho_K).$$

**Proof:** Consider $\mathfrak{p}$ as a subspace of the Clifford algebra $Cl(B, \mathfrak{p})$. We will make $S = \wedge^* \mathfrak{p}_-$ a $Cl(B, \mathfrak{p})$-module. For this let $x \in \mathfrak{p}_-$ act on $S$ via

$$x.z_1 \wedge \ldots \wedge z_n = \sqrt{2} \, x \wedge z_1 \wedge \ldots \wedge z_n,$$

and $y \in \mathfrak{p}_+$ via

$$y.z_1 \wedge \ldots \wedge z_n = \sqrt{2} \sum_{i=1}^{n} (-1)^{i+1} B(y, z_i) z_1 \wedge \ldots \hat{z}_i \ldots \wedge z_n.$$

This prescription turns $S$ into a nontrivial $Cl(B, \mathfrak{p})$-module. Since there is only one such of the dimension of $S$, we conclude that $S$ is the nontrivial irreducible $Cl(B, \mathfrak{p})$-module. Therefore the argumentations of [AtSch] apply to $S$ and especially the formula of Parasarathy ([AtSch], (A13)). Note that in [AtSch] everything was done under the assumption of $G$ being connected. Since $G$ is of inner type, however, $\pi(C)$ will be a scalar. Writing $G^0$ for the connected component, the representation $\pi|_{G^0}$ will decompose as a finite sum of irreducibles on each of which the formula of Parasarathy holds. Thus it holds globally. Let $d_\pm$ be defined as in loc. cit. then our assumption leads to $ker(d_+ d_-) \cap \pi \otimes S(\tau) \neq 0$, and therefore $0 = \tau(C_K) - \pi(C) - B(\rho) + B(\rho_K)$. ∎

**Lemma. 5.1.4** *Assume $c = c(H) = 2$ and let $(\sigma, W_\sigma)$ be an irreducible $M$-subrepresentation of $\wedge^l \mathfrak{n}$ such that $W_\sigma \cap \wedge^l \mathfrak{n}_- \neq 0$. Since $M$ is of inner type there is an infinitesimal character $\lambda_\sigma$ of $\sigma$ and we have*

$$B(\lambda_\sigma) - B(\rho) = -(1 + 2\dim(\mathfrak{n}_-) - l)^2 B(\frac{\alpha_r}{2}),$$

*where $\alpha_r$ is the unique positive real root of $\mathfrak{t}_M \oplus \mathfrak{a}$.*

**Proof:** The assertion follows from Corollary 5.7 of [Kost], here one uses the fact that the boundary map $\partial$ of the Lie algebra homology $H_*(\mathfrak{n})$ maps $\bigwedge^*(\mathfrak{n}_+ \oplus \mathfrak{n}_-)$ to $\mathfrak{n}_r \wedge \bigwedge^*(\mathfrak{n}_+ \oplus \mathfrak{n}_-)$ and maps $\bigwedge^* \mathfrak{n}_-$ to zero. ∎

**Lemma. 5.1.5** *Assume $c = c(H) = 2$ and let $(\tau, V_\tau)$ be an irreducible $K_M$-subrepresentation of $\wedge^l \mathfrak{n}_-$, then the Casimir operator $C_{K_M}$ of $K_M$ acts on $V_\tau$ by the scalar*

$$\tau(C_{K_M}) = B(\rho) - (1 + 2\dim(\mathfrak{n}_-) - l)^2 B(\frac{\alpha_r}{2}) - B(\rho_{K_M}) - \frac{2l B(\rho_{M,n})}{\dim(\mathfrak{p}_{M,-})},$$

*where $\rho_{M,n}$ is the half of the sum of the noncompact roots of $M$, so $\rho_{M,n} = \rho_M - \rho_K$. In the case that $\mathrm{rank}(X) = 1$, the last summand does not occur.*

**Convention.** In the case of rank one the last summand is undefined. To keep the formulas unified we will agree to consider it as zero then.

**Proof:** Let $x_1, \ldots, x_n$ be a basis of $\mathfrak{p}_{M,+}$ consisting of root vectors with respect to $\mathfrak{t}_M$ and such that $B(x_j, x_j^c) = 1$, where $x^c$ is the complex conjugate of $x$. Let $L$ be the element of the Lie algebra of $\mathfrak{k}_M$ defined by

$$L = \sum_j [x_j, x_j^c].$$

Since $(x_j^c)$ is the dual basis to $(x_j)$ it follows that $L$ lies in the center of $\mathfrak{k}_M$. Furthermore $L$ is imaginary so that weights will take real values on $L$.

Consider the $K_M$-module $\mathfrak{n}_-$. We will show now that the center of $\mathfrak{k}_M$ acts on $n_-$ by a character. To this end recall that $\mathfrak{t}$ is spanned by $\mathfrak{t}_M$ and $X + \theta(X)$, where $0 \neq X \in \mathfrak{n}_r$. Let $Z$ be a generator of the center of $\mathfrak{k}$ then

$$Z = a(X + \theta(X)) + bL + R$$

where $R$ lies in $L^\perp \cap \mathfrak{t}_M$. Since the latter is just the dual space of the span of the compact roots in $\mathfrak{t}_M^*$ it follows $R = 0$.

Now $Z$ acts as a scalar $\mu$ on $\psi(\mathfrak{n}_-)$ and it is easy to see that $ab \neq 0$. Let $Y \in \mathfrak{g}_\alpha \subset \mathfrak{n}_-$. A computation shows

$$\mu Y = a[X, \theta(Y)] + b[L, Y].$$

So it remains to show that $Y \mapsto [X, \theta(Y)]$ acts as a scalar on $\mathfrak{n}_-$. To this end let $H = [X, \theta(X)] \in \mathfrak{a}$ then

$$(\alpha_r(H)/2)Y = [H, Y] = [X, [\theta(X), Y]].$$

Since root spaces are one dimensional there exists a number $c = c_Y$ such that $[X, \theta(Y)] = cY$. Therefore $[\theta(X), Y] = c\theta(Y)$ and $(\alpha_r(H)/2)Y = c^2 Y$. By this $c$ is determined up to sign. So $L$ acts on $\mathfrak{n}_-$ by at most two eigenvalues. Suppose there are two eigenvalues and an according decomposition

$$\mathfrak{n}_- = \mathfrak{n}_-^0 \oplus \mathfrak{n}_-^1,$$

where $\mathfrak{n}_-^0$ corresponds to the character of less absolute value. We will show that $\mathfrak{g}_1 \cap (\mathfrak{a} \oplus \mathfrak{m} \oplus \mathfrak{n}_-^0 \oplus \mathfrak{n}_-^{0\,c} \oplus \theta(\mathfrak{n}_-^0) \oplus \theta(\mathfrak{n}_-^{0\,c}))$ is an ideal of $\mathfrak{g}_1$ which cannot be since $\mathfrak{g}_1$ is simple, thus proving that there is only one eigenvalue of $L$ on $\mathfrak{n}_-$. For this it suffices to show $[\mathfrak{p}_{M,-}, \mathfrak{n}_-] = 0$. Let $X \in \mathfrak{n}_-$ then $X - \theta(X)$ is in $\mathfrak{p}_-$ and therefore $[\mathfrak{p}_{M,-}, X - \theta(X)] = 0$. But $\mathfrak{p}_M$ also leaves stable the $\mathfrak{a}$-root spaces, therefore $[\mathfrak{p}_{M,-}, X] = 0$. We have shown that the center of $\mathfrak{k}_M$ acts be a character $\chi_{\mathfrak{n}_-}$ on $\mathfrak{n}_-$. As for the value of this character recall that the center of $\mathfrak{k}_M$ acts by a character $\chi$ on $\mathfrak{p}_{M,-}$. The assumption $\chi \neq \chi_{\mathfrak{n}_-}$ would similarly allow us to construct a nontrivial ideal, hence it follows $\chi = \chi_{\mathfrak{n}_-}$.

Recall that each $x_j$ is a root vector, say $x_j \in \mathfrak{m}_\alpha$ then $x_j^c \in \mathfrak{m}_{-\alpha}$ and $[x_j, x_j^c] = H_\alpha$ where $H_\alpha \in \mathfrak{t}_M$ is defined by $\alpha(H) = B(H, H_\alpha)$. Therefore

$$L = \sum_{\alpha \in \phi_{noncompact}^+(\mathfrak{t}_M, \mathfrak{m})} H_\alpha.$$

So let $\alpha$ be a noncompact negative root in $\phi(\mathfrak{t}_M, \mathfrak{m})$. Since $L$ acts on $\mathfrak{n}_-$ by the same scalar as on $\mathfrak{p}_{M,-}$ we get

$$\chi_{\mathfrak{n}_-}(L) = 2B(\rho_{M,n}, \alpha) = -\frac{2B(\rho_{M,n})}{\dim(\mathfrak{p}_{M,-})}.$$

Now let $\mu$ be the lowest weight of $\tau$ and let $\sigma$ be the $M$-representation generated by $\tau$. We may assume that $\sigma$ is irreducible. We claim that $\mu$ is the lowest weight vector of $\sigma$. To see this it suffices to see that $\sigma(\mathfrak{p}_{M,-})V_\tau = 0$, which in turn follows from $\sigma(\mathfrak{p}_{M,-})\mathfrak{n}_- = [\mathfrak{p}_{M,-}, \mathfrak{n}_-] = 0$. We thus have shown that $\mu$ is the lowest weight vector of $\sigma$. So by Lemma 5.1.4 we get

$$B(\mu - \rho_M) - B(\rho) = -(1 + 2\dim(\mathfrak{n}_-) - l)^2 B(\frac{\alpha_r}{2}).$$

On the other hand we know $\tau(C_{K_M}) = B(\mu - \rho_{K_M}) - B(\rho_{K_M})$. We have

$$B(\mu - \rho_M) = B(\mu - \rho_{K_M} - \rho_{M,n}).$$

Note that $2\rho_{M,n}$ is the dual of $L$. Since $L$ is in the center of $\mathfrak{k}_M$ it follows

$$B(\rho_{K_M}, \rho_{M,n}) = 0.$$

Therefore

$$\begin{aligned}
B(\mu - \rho_M) &= B(\mu - \rho_{K_M}) - 2B(\mu - \rho_{M,n}) \\
&= B(\mu - \rho_{K_M}) - \mu(L) \\
&= B(\mu - \rho_{K_M}) + \frac{2lB(\rho_{M,n})}{\dim(\mathfrak{p}_{M,-})}.
\end{aligned}$$

From this the claim follows. ∎

Define $\hat{f}^0_{t,H}$ as in [HC-S]. Letting $b^* \in B^*$ we therefore get in the case that $c(H) = 2$:

$$\hat{f}^0_{t,H}(\nu, b^*) = -e^{tB(\nu)}$$

$$\times \sum_{l=0}^{\dim(\mathfrak{n}_-)} \sum_{p=0}^{\dim \mathfrak{p}_{M,-}} (-1)^{p+l} e^{-t((1+2\dim(\mathfrak{n}_-)-l)^2 B(\frac{\alpha_r}{2}) + \frac{2lB(\rho_{M,n})}{\dim(\mathfrak{p}_{M,+})})} \dim(V_{b^*} \otimes \wedge^p \mathfrak{p}_{M,-} \otimes \wedge^l \mathfrak{n}_-)^{K_M}.$$

In the other case, $c(H) = 1$ we get

**Lemma. 5.1.6** *Assume $c(H) = 1$ then for $\pi_{\xi,\nu}$ as in Lemma 5.1.1 we get*

$$\mathrm{tr}(\pi_{\xi,\nu}(f^0_t)) = -e^{t(B(\nu)+B(\rho_{K\cap G_2})+B(\rho_{M_1})-B(\rho))}$$

$$\times (\sum_{j=0}^{\dim(\mathfrak{p}_{M_1})} (-1)^j \dim(W_{\xi_1} \otimes \wedge^j \mathfrak{p}_{M_1})^{K\cap M_1})(\sum_{k+0}^{\dim(\mathfrak{p}_{G_2,-})} (-1)^k \dim(W_{\xi_2} \otimes \wedge^k \mathfrak{p}_{G_2,-})^{K\cap G_2})$$

**Proof:** We have to show that $\mathrm{tr}(\pi_{\xi,\nu}(f^0_t)) \neq 0$ implies $B(\lambda_\xi) = B(\rho_{M_1} + \rho_{K\cap G_2})$. We have $M = M_1 \times G_2$ and thus $\xi = \xi_1 \otimes \xi_2$. Therefore

$$\mathrm{tr}(\pi_{\xi,\nu}(f^0_t)) = -e^{t\pi_{\xi,\nu}(C)} \sum_{q=0}^{\dim(\mathfrak{p}_-)-1} (-1)^q \sum_{j+k=q} \dim(W_{\xi_1} \otimes \wedge^j \mathfrak{p}_{M_1})^{K\cap M_1} \dim(W_{\xi_2} \otimes \wedge^k \mathfrak{p}_{G_2,-})^{K\cap G_2}$$

$$= -e^{t\pi_{\xi,\nu}(C)} (\sum_{j=0}^{\dim(\mathfrak{p}_{M_1})} (-1)^j \dim(W_{\xi_1} \otimes \wedge^j \mathfrak{p}_{M_1})^{K\cap M_1})(\sum_{k+0}^{\dim(\mathfrak{p}_{G_2,-})} (-1)^k \dim(W_{\xi_2} \otimes \wedge^k \mathfrak{p}_{G_2,-})^{K\cap G_2})$$

Now assume $\mathrm{tr}(\pi_{\xi,\nu}(f^0_t)) \neq 0$ then Lemma 5.1.3 implies $B(\lambda_{\xi_2}) = B(\rho_{K\cap G})$ whereas Lemma 2.4 in [MS-2] gives $B(\lambda_{\xi_1}) = B(\rho_{M_1})$. ∎

Therefore in the case $c(H) = 1$ we get

$$\hat{f}^0_{t,H}(\nu, b^*) = \hat{f}^0_{t,H}(\nu, b^*_1 + b^*_2) = -e^{t(B(\nu)+B(\rho_{K\cap G_2})+B(\rho_{M_1})-B(\rho))}$$

$$(\sum_{j=0}^{\dim(\mathfrak{p}_{M_1})} (-1)^j \dim(V_{b^*_1} \otimes \wedge^j \mathfrak{p}_{M_1})^{K\cap M_1})(\sum_{k+0}^{\dim(\mathfrak{p}_{G_2,-})} (-1)^k \dim(V_{b^*_2} \otimes \wedge^k \mathfrak{p}_{G_2,-})^{K\cap G_2}).$$

Since $\Gamma$ is the fundamental group of $X_\Gamma$, every conjugacy class $[\gamma]$ in $\Gamma$ defines a free homotopy class of closed paths in $X_\Gamma$. It is known [DKV] that the union $X_\gamma$ of all closed geodesics which are homotopic to $[\gamma]$ is a smooth submanifold of $X_\Gamma$. Let $\chi_1(X_\gamma)$ denote the first higher Euler number of $X_\gamma$ (see [D-Hitors]), i.e.

$$\chi_1(X_\gamma) = \sum_{p=0}^{\dim(X_\gamma)} p(-1)^p b_p(X_\gamma),$$

where $b_p(X_\gamma)$ is the $p$-th Betti number of $X_\gamma$.

Now let $\mathcal{E}_H(\Gamma)$ denote the set of nontrivial $\Gamma$-conjugacy classes, which are in $G$ conjugate to an element of H. For $l \geq 0$ define

$$b_l(H) = (\frac{c(H)}{2} + \dim(\mathfrak{n}) - 1 - l)|\frac{\alpha_r}{c(H)}|.$$

If $c(h) = 2$ let

$$d_l(H) = \sqrt{b_l(H)^2 + \frac{2lB(\rho_{M,n})}{\dim(\mathfrak{p}_{M,-})}}.$$

In the case $c(H) = 1$ we finally set

$$d(H) = \sqrt{B(\rho) - B(\rho_{K \cap G_2}(H)) - B(\rho_{M_1}(H))}.$$

With an argumentation as in [D-Hitors] we get

**Theorem 5.1.1** *Let $X_\Gamma$ be a compact locally hermitian space with fundamental group $\Gamma$ and such that the universal covering is globally symmetric of the noncompact type. Let $\triangle_{p,q,\varphi}$ be the Hodge Laplacian on (p,q)-forms with values in a flat bundle $E_\varphi$, then*

$$\Theta(t) = \sum_{q=0}^{\dim_{\mathbb{C}} X_\Gamma} q(-1)^q Tre^{-t\triangle_{0,q,\varphi}}$$

$$= \sum_{\substack{H/conj. \\ c(H) = 2}} \sum_{[\gamma] \in \mathcal{E}_H(\Gamma)} \frac{l_{\gamma_0} \chi_1(X_\gamma) \mathrm{tr}\varphi(\gamma)}{|W(\mathfrak{h}, \mathfrak{g}_\gamma)| \prod_{\alpha \in \phi_\gamma^+} B(\rho_\gamma, \alpha) \det(1 - \gamma^{-1}|\mathfrak{n})} \frac{e^{-l_\gamma^2/4t}}{\sqrt{4\pi t}}$$

$$\sum_{l=0}^{\dim(\mathfrak{n}_-)} (-1)^l e^{-td_l(H)^2} \frac{\tilde{\omega}_\gamma(\gamma^{\rho_M} \det(1 - \gamma^{-1}|(\mathfrak{k}/\mathfrak{t})^+) \mathrm{tr}(\gamma_I| \wedge^l \mathfrak{n}_-))}{\gamma^{\rho_M} \prod_{\alpha \in \phi_I^+ - \phi_\gamma^+} (1 - \gamma^{-\alpha})}$$

$$+ \sum_{\substack{H/conj. \\ c(H) = 1}} \sum_{[\gamma] \in \mathcal{E}_H(\Gamma)} l_{\gamma_0} \chi_1(X_\gamma) \frac{e^{-l_\gamma^2/4t}}{\sqrt{4\pi t}} e^{-b_0(H)l_\gamma} e^{td(H)^2}$$

$$+ f_t^0(e) \dim\varphi \ \mathrm{vol}(X_\Gamma),$$

*where $\tilde{\omega}_\gamma$ is the usual differential operator [HC-DS] p. 33. For unexplained notation we refer to [D-Hitors].*

The reader should keep in mind that by its definition we have for the term of the identity:

$$f_t^0(e) \, \dim\varphi \, \text{vol}(X_\Gamma) = \sum_{q=0}^{\dim_\mathbb{C} X} q(-1)^q \text{tr}_\Gamma(e^{-t\triangle_{0,q,\varphi}}),$$

where $\text{tr}_\Gamma$ is the $\Gamma$-trace (1.2.1). Further note that by the Plancherel theorem the Novikov-Shubin invariants of all operators $\triangle_{0,q}$ are positive.

Now assume we are given a unitary irreducible representation $(\tau, V_\tau)$ of $K$. The representation $\tau$ defines a $G$-homogeneous vector bundle $E_\tau$ over $X = G/K$ and we will consider the space of smooth sections of $E_\tau$:

$$C^\infty(E_\sigma) = (C^\infty(G) \otimes V_\tau)^K.$$

The Casimir operator $C$ of $G$ acts on this space and defines a second order elliptic differential operator $C_\tau$ on $E_\tau$. The methods of [BM] indicate that the heat operator $e^{tC_\tau}$, $t > 0$, acts by a smooth convolution kernel $h_t^\tau$ which is in the Harish-Chandra Schwartz space $\mathcal{S}(G)$. For any $\pi \in \hat{G}$ it follows

$$\text{tr}\pi(h_t^\tau) = e^{t\pi(C)} \dim(V_\pi \otimes V_\tau)^K.$$

By Lemma 5.1.3 there is a linear combination $g_\tau$ of the functions $h_t^{\tau'}$ such that $g_\tau \in \mathcal{S}(G)$ and

$$\text{tr}\pi(g_\tau) = \sum_{p=0}^{\dim(\mathfrak{p}_-)} \dim(V_\pi \otimes \wedge^p \mathfrak{p}_- \otimes V_\tau)^K.$$

Similar argumentations as in 5.1.1 lead to $\text{tr}\pi_{\xi,\nu}(g_\tau) = 0$ for any properly induced representation $\pi_{\xi,\nu}$.

Now let $(\sigma, W_\sigma)$ be a finite dimensional irreducible representation of $G$. Assume $W_\sigma$ is furnished with a scalar product such that $\sigma|_K$ is unitary. Lemma 2.4 of [MS-2] says that if

$$\sum_{p=0}^{\dim \mathfrak{p}} (-1)^p \dim(V_\pi \otimes \wedge^p \mathfrak{p} \otimes W_\sigma)^K \neq 0,$$

then

$$B(\lambda_\pi) = B(\lambda_\sigma).$$

So there is a function $f_\sigma \in \mathcal{S}(G)$ such that for any $\pi \in \hat{G}$:

$$\text{tr}\pi(f_\sigma) = \sum_{p=0}^{\dim \mathfrak{p}} (-1)^p \dim(V_\pi \otimes \wedge^p \mathfrak{p} \otimes W_\sigma)^K.$$

Let $f_\sigma, g_\tau$ as in the previous section. We want to compute orbital integrals of these functions. Recall that an element $g$ of $G$ is called **elliptic** if it lies in a compact Cartan subgroup.

**Proposition 5.1.2** *Let $g$ be a semisimple element of the group $G$. If $g$ is not elliptic, the orbital integrals $\mathcal{O}_g(f_\sigma)$ and $\mathcal{O}_g(g_\tau)$ vanish. If $g$ is elliptic we may assume $g \in T$, where $T$ is a Cartan in $K$ and then we have*

$$\mathcal{O}_g(f_\sigma) = \frac{\overline{\mathrm{tr}\sigma(g)}|W(\mathfrak{t}, \mathfrak{g}_g)| \prod_{\alpha \in \Phi_g^+}(\rho_g, \alpha)}{[G_g : G_g^0]c_g},$$

*for all elliptic $g$ and*

$$\mathcal{O}_g(g_\tau) = \frac{\overline{\mathrm{tr}\tau(g)}}{\det(1 - g^{-1}|\mathfrak{p}_+)},$$

*if $g$ is regular elliptic. For general elliptic $g$ we have*

$$\mathcal{O}_g(g_\tau) = \frac{\sum_{s \in W(T,K)} \det(s)\tilde{\omega}_g g^{s\lambda_{\tau_*} + \rho - \rho_K}}{[G_g : G_g^0]c_g g^\rho \prod_{\alpha \in \phi^+ - \phi_g^+}(1 - g^{-\alpha})},$$

*where $c_g$ is Harish-Chandra's constant, it does only depend on the centralizer $G_g$ of $g$. Its value is given in [D-Hitors], further $\tilde{\omega}_\gamma$ is the differential operator as in [HC-DS] p.33.*

**Proof:** The vanishing of $\mathcal{O}_g(f_\sigma)$ for nonelliptic $g$ is immediate by Harish-Chandra's formula for the Fourier transform of orbital integrals [HC-S]. Now let $g \in K \cap G'$, where $G'$ denotes the set of regular elements. We apply Lemma 4.3 of [MS-1] to get

$$\mathcal{O}_g(f_\sigma) = \overline{\mathrm{tr}\sigma(g)}.$$

This proves the proposition in the regular case. The general case is derived from this by standard considerations (see [HC-DS], p32 ff.). A word of comment is in order here. In [MS-1] everything is formulated in terms of a certain subgroup $M^+$ of $M$. One could either do the same here or see that Lemma 4.3 in [MS-1] actually holds for $M$ as well since the latter still is in the Harish-Chandra class (Lemma 4.9 in [HC-HA1]). One sees by Frobenius reciprocity that the right hand side of Lemma 4.2 in [MS-1] actually gives the same for $M$ as for $M^+$.

The case $g_\tau$ is treated similarly. ∎

## 5.2 A Simple Example

Let m be a natural number and consider the group $G = SL_2(\mathbb{R})^m$. A maximal compact subgroup is given by $K = SO(2)^m$ and we have $G/K \cong \mathcal{H}_2(\mathbb{R})^m$. On the latter there is a canonical metric giving it the curvature $-1$. To induce this metric we normalize the form $B$ such that $B(\alpha) = 1$ for every root $\alpha$. For Haar measure we take the Euler-Poincaré measure.

The Cartan subgroups of splitrank one are modulo conjugation given by

$$A \times T \times \ldots \times T, \; T \times A \times T \times \ldots \times T, \; \ldots T \times \ldots \times T \times A,$$

where

$$A = \left\{ \begin{pmatrix} a & \\ & a^{-1} \end{pmatrix}, a \neq 0 \right\}, \quad T = SO(2).$$

Hence it follows that all the constants $c(H)$ in Theorem 6.1.1 are equal to $\frac{1}{2}$, all spaces $\mathfrak{n}_\pm$ are trivial and thus the summation index $l$ only takes the value zero. As to $f_t^0(e)$ we have

$$\begin{aligned} f_t^0(e) &= \sum_{q=0}^m q(-1)^q \mathrm{tr} < x|e^{-t\triangle_{0,q,\varphi}}|x > \\ &= \dim(\varphi) \sum_{q=0}^m q(-1)^q \mathrm{tr} < x|e^{-t\triangle_{0,q}}|x > \end{aligned}$$

Now let $\triangle^0$ denote the laplacian on $\mathcal{H}_2(\mathbb{R})$ and let $\triangle^1$ twice the Laplacian on $(0,1)$-forms on $\mathcal{H}_2(\mathbb{R})$. Then we have

**Lemma. 5.2.1**

$$\mathrm{tr} < x|e^{-\triangle^1}|x > -\mathrm{tr} < x|e^{-\triangle^0}|x > = \frac{1}{2}$$

**Proof:** Write for the moment $G = SL_2(\mathbb{R})$ and $K = SO(2)$. Let $h_t^0$ denote the convolution kernel of $e^{-t\triangle^0}$ in

$$(C^\infty(G) \otimes \mathrm{End}(\mathfrak{p}_-))^{K \times K}$$

and write $h_t^1$ for the $\triangle^1$-analogue. Let $g_t^0(x) = \mathrm{tr} h_t^0(x)$ and $g_t^1(x) = \mathrm{tr} h_t^1(x)$, where the trace is the trace of $\mathrm{End}(\mathfrak{p}_-)$. The lemma amounts to

$$g_t^1(e) - g_t^0(e) = \frac{1}{2}.$$

We apply the Plancherel theorem to $g_t^0$ and $g_t^1$. For any $\pi \in \hat{G}$ we get

$$\mathrm{tr}\pi(G_t^0) = e^{t\pi(C)} \dim(V_\pi)^K,$$

$$\mathrm{tr}\pi(G_t^1) = e^{t\pi(C)} \dim(V_\pi \otimes \mathfrak{p}_-)^K.$$

There is, however, only one $\pi \in \hat{G}$ for which these two dimensions differ, namely $\pi$ equals the discrete series representation with Harish-Chandra parameter $\lambda_\pi = \rho$ and for this representation we have with our normalizations:

$$\pi(C) = 0, \ \dim(V_\pi)^K = 0, \ \dim(V_\pi \otimes \mathfrak{p}_-)^K = 1$$

and $d_\pi = \frac{1}{2}$, whence the lemma. ∎

Now since $X = \mathcal{H}_2(\mathbb{R})^m$ is a product space, the heat kernels can be written as products, so

$$e^{-t\triangle_{0,q}} = \sum_{I \subset \{1,\dots;m\}|I|=q} e^{-t\triangle_{1_I}(1)} \otimes \dots \otimes e^{-t\triangle_{1_I}(m)},$$

where $1_I$ is the characteristic function of $I$. From this we take

$$
\begin{aligned}
f_t^e(e) &= \dim(\varphi) \sum_{q=0}^{m} q(-1)^q \sum_{\substack{I \subset \{1,\dots,m\} \\ |I| = q}} \prod_{j=1}^{m} \text{tr}(, x|e^{-t\triangle_{1_I}(j)}|x> \\
&= \dim(\varphi) \sum_{q=0}^{m} q(-1)^q \binom{m}{q} g_t^0(e)^{m-q} g_t^1(e)^q \\
&= \dim(\varphi) m \sum_{q=1}^{m} \binom{m-1}{q-1} g_t^0(e)^{m-q}(-g_t^1(e))^q \\
&= \dim(\varphi) m(-g_t^1(e))(g_t^0(e) - g_t^1(e))^{m-1} \\
&= \dim(\varphi) m(-1)^m (\tfrac{1}{2} + g_t^0(e)) 2^{1-m}
\end{aligned}
$$

These calculations now easily lead to

**Theorem 5.2.1** *Let $m$ be a natural number and $\Gamma$ a nice cocompact lattice in $G = SL_2(\mathbb{R})^m$. Let $\varphi$ be a finite dimensional unitary representation of $\Gamma$. The zeta function*

$$Z_\varphi(s) \stackrel{\text{def}}{=} \prod_{\substack{[\gamma]\text{prime} \\ \mathbb{R}-\text{rank}[\gamma]=1}} \prod_{k,n \geq 0} \det(1 - e^{(s+n)l_\gamma} \varphi(\gamma) \otimes \text{Ad}(\gamma_I | S^k(\mathfrak{p}_{M,-})))^{\chi_1(X_\gamma)}$$

*extends to a meromorphic function on $\mathbb{C}$. One has*

$$Z_\varphi(\tfrac{1}{2} + s) = \prod_{q=0}^{m} \det(\triangle_{0,q,\varphi} + s^2 - \tfrac{1}{4})^{q(-1)^q} \left( \exp(s^2) \det(P+s) \right)^{-2m\dim\varphi \text{Ar}(X_\Gamma)}$$

*where $\text{Ar}(X_\Gamma)$ is the arithmetic genus of $X_\Gamma$ and $P = \sqrt{\triangle^d + \tfrac{1}{4}}$ as in Theorem 4.1.1*∎

**Corollary 5.2.1** *The function $\hat{Z}_\varphi(s) = Z_\varphi(s)\det(P+s-\tfrac{1}{2})^{2m\dim\varphi\text{Ar}(X_\Gamma)}$ satisfies*

$$\hat{Z}_\varphi(s) = \hat{Z}_\varphi(1-s).\blacksquare$$

**Corollary 5.2.2** *The "factor at infinity" is for* Re $(s) >> 0$ *representable as*

$$\left(\exp(s^2)\det(P+s)\right)^{-2m\dim\varphi\mathrm{Ar}(X_\Gamma)} = \prod_{q=0}^{m}\left(\det{}^{(2)}(\triangle_{0,q,\varphi}+s^2-\frac{1}{4})\right)^{q(-1)^{q+1}}.\blacksquare$$

**Corollary 5.2.3** *The function $Z_\varphi$ vanishes at $s = \frac{1}{2}$ up to order*

$$n_0 = -2m\dim\varphi\mathrm{Ar}(X_\Gamma) + \sum_{q=0}^{m}q(-1)^q h_{0,q}(\varphi),$$

*where $h_{p,q}(\varphi)$ is the $(p,q)$-th Hodge number. With $R_\varphi(s) = Z_\varphi(s)s^{-n_0}$ we have*

$$R_\varphi(\frac{1}{2}) = \frac{T_{\mathrm{hol}}(X_\Gamma,\varphi)}{T_{\mathrm{hol}}^{(2)}(X_\Gamma)^{\dim\varphi}}.$$

*Further we have*

$$T_{\mathrm{hol}}^{(2)}(X_\Gamma) = \left(2^{12}\exp(\zeta'(1))\right)^{-2m\dim\varphi\mathrm{Ar}(X_\Gamma)},$$

*where $\zeta$ is the zeta function of Riemann.* $\blacksquare$

In the theorem we assumed $\Gamma$ to be nice in order to have a simple Euler product. Note that in the case $G = SL_2(\mathbb{R})$ every torsionfree $\Gamma$ automatically is nice. But the theorem also extends to the case of nonnice groups as follows:

Let $\gamma \in \Gamma$ then $\gamma$ is the power of some prime $\gamma_0$, but since $\Gamma$ s not nice, it may happen that $G_\gamma \neq G_{\gamma_0}$. But then there is a least power $\gamma_0^n$ of $\gamma_0$ such that $G_{\gamma_0^n} = G_\gamma$. Let $\lambda_\gamma$ denote the length of the geodesic attached to $\gamma_0^n$, the so called **generic geodesic** in $X_\gamma$. Then let $\mu_\gamma = l_\gamma/\lambda_\gamma$ the **generic multiplicity** of $\gamma$. For example, assume that the only root of unity of $\gamma_0$ is $-1$, then we have

$$\mu_{\gamma_0^{2n}} = n \quad \mu_{\gamma_0^{2n+1}} = 2n+1.$$

Further for $\gamma \in \Gamma$ let $Per(\gamma)$ denote the set of orders of nontrivial roots of unity occurring as eigenvalues of the adjoint representation $\mathrm{Ad}(\gamma)$, so $m$ is in $Per(\gamma)$ if and only if a primitive $m$-th root of unity occurs as an eigenvalue of $\mathrm{Ad}(\gamma)$. We have: $\Gamma$ is nice if and only if $Per(\gamma) = \emptyset$ for all $\gamma \in \Gamma$.

For any subset $I \subset Per(\gamma)$ let $n_I$ denote the least common multiple of the elements of $I$. Set $n_\emptyset = 1$. Consider

$$Z_\varphi(s) \overset{\mathrm{def}}{=} \prod_{\substack{[\gamma]\mathrm{prime} \\ \mathbb{R}-\mathrm{rank}[\gamma]=1}} \prod_{I\subset Per(\gamma)}\prod_{k,n\geq 0}\det(1-e^{n_I(s+n)l_\gamma}\varphi(\gamma)\otimes\mathrm{Ad}(\gamma_I|S^k(\mathfrak{p}_{M,-})))^{m_I(\gamma)},$$

where

$$m_I(\gamma) = \sum_{J \subset I} (-1)^{|I|+|J|+1} n_J \chi(\hat{X}_{\gamma^{n_J}})/\mu_{\gamma^{n_J}}.$$

Here $\chi(\hat{X}_{\gamma^{n_J}})$ denotes the orbifold Euler characteristic of $\hat{X}_{\gamma^{n_J}} = X_{\gamma^{n_J}}/(\text{geodesic flow})$. Note that $\chi(\hat{X}_{\gamma^{n_J}})$ is no longer an integer but a rational number. This cannot cause any trouble since Euler factors are all close to one and the powers are defined by the principal branch of the logarithm.

**Theorem 5.2.2** *With* $\Gamma \subset SL_2(\mathbb{R})$ *torsionfree and cocompact and* $Z_\varphi$ *defined as above the preceding theorem and all corollaries hold as well.*∎

## 5.3   The Rank One Case

Now assume that the rank of $X$ is one. From the classification of globally symmetric spaces it follows that $X = SU(m,1)/S(U(m) \times U(1))$. So that $G = SU(m,1)/center$ and $K = S(U(m) \times U(1))/center$. Since it doesn't change the geometry we will instead work with $G = SU(m,1)$ and $K = S(U(m) \times U(1))$. The Lie algebra $\mathfrak{g}_0$ can be written as

$$\mathfrak{g}_0 = \left\{ \begin{pmatrix} Z_1 & Z_2 \\ Z_2^* & -\text{tr}Z_1 \end{pmatrix} \begin{array}{l} Z_1 \in \text{Mat}_m(\mathbb{C}), \ Z_1^* = -Z_1 \\ Z_2 \in \mathbb{C}^m \text{ arbitrary} \end{array} \right\}$$

where for a complex matrix $A$ we write $A^* = \bar{A}^t$. The Lie algebra $\mathfrak{g}_0$ is considered as a real Lie algebra but it is formed out of complex matrices. To avoid conflicts with the complexification we will write $Z \in \mathfrak{g}_0$ as $Z = X + IY$ with $X, Y \in \text{Mat}_{m+1}(\mathbb{R})$ further we write $\bar{Z} = X - IY$ whereas $Z^c$ will denote the complex conjugation on $\mathfrak{g}$ as before.

Instead of the Killing form on $\mathfrak{g}$ we will rather use the form

$$B = \frac{1}{4(m+1)} \times \text{Killing form}.$$

Then it follows that $B(X,Y) = \frac{1}{2}\text{Re}(\text{tr}(XY))$ where Re stands for the real part on $\mathfrak{g}_0$.

A maximal split torus of $\mathfrak{g}_0$ is given by

$$\mathfrak{a}_0 = \mathbb{R}H_0 \quad \text{where} \quad H_0 = \begin{pmatrix} \mathbf{0} & & \\ & 0 & 1 \\ & 1 & 0 \end{pmatrix}.$$

The Lie algebra $\mathfrak{k}_0$ is given by

$$\mathfrak{k}_0 = \left\{ \begin{pmatrix} Z & 0 \\ 0 & -\text{tr}Z \end{pmatrix} \mid Z^* = -Z \right\},$$

and the centralizer of $\mathfrak{a}_0$ is $\mathfrak{a}_0 + \mathfrak{m}_0$, where

$$m_0 = \left\{ \begin{pmatrix} X & 0 & 0 \\ 0 & -\frac{1}{2}\mathrm{tr}X & 0 \\ 0 & 0 & -\frac{1}{2}\mathrm{tr}X \end{pmatrix} \begin{array}{c} X \in \mathrm{Mat}_{\mathrm{m}-1}(\mathbb{C}) \\ X^* = -X \end{array} \right\}.$$

Let $\mathfrak{t}_0$ denote the space of diagonal matrices in $\mathfrak{g}_0$, then $T = \exp(\mathfrak{t}_0)$ is a compact Cartan subgroup of $G$ lying in $K$. Let $\mathfrak{t}_{\mathrm{m},o} = \mathfrak{t}_0 \cap \mathfrak{m}_0$. This is a Cartan subalgebra of $m_0$ and $\mathfrak{t}_{\mathrm{m}} \oplus \mathfrak{a}$ is a Cartan of $\mathfrak{g}$. We are going to consider the root system $\Phi(\mathfrak{t}_{\mathrm{m}} \oplus \mathfrak{a}, \mathfrak{g})$ of this torus. We have $\Phi(\mathfrak{t}_{\mathrm{m}} \oplus \mathfrak{a}, \mathfrak{g}) = \Phi(\mathfrak{t}_{\mathrm{m}}, \mathfrak{m}) \cup \Phi_{nI}(\mathfrak{t}_{\mathrm{m}} \oplus \mathfrak{a}, \mathfrak{g})$, here nI means nonimaginary roots. There is up to sign a unique real root $\alpha_r$. We have $\alpha_r(H_0) = 2$.

Choose a real parabolic subalgebra $\mathfrak{a}_0 \oplus \mathfrak{m}_0 \oplus \mathfrak{n}_0$ and fix an ordering on the roots compatible with the choice of the parabolic. Then the unipotent radical $\mathfrak{n}$ splits as $n = \mathfrak{n}_r \oplus \mathfrak{n}_- \oplus \mathfrak{n}_+$ where $n_r = \mathfrak{g}_{\alpha_r}$ and

$$\mathfrak{n}_{\pm} = \bigoplus_{\substack{\alpha \in \Phi_{nI}^+ \\ \pm\alpha \mid \mathfrak{t}_{\mathrm{m}} > 0}} \mathfrak{g}_{\alpha}.$$

We will make use of the results in chapter 4. In chapter 4 we considered the **generalized Selberg Zeta Function**:

$$Z_{\sigma,\varphi}(s) = \prod_{\substack{1 \neq [\gamma] \\ \text{prime}}} \prod_{N \geq 0} \det(1 - \mathrm{e}^{-\mathrm{sl}_\gamma} \gamma \mid \mathrm{V}).$$

The product is taken over all nontrivial primitive conjugacy classes in $\Gamma$ and any $\gamma \in \Gamma$ induces an operator, well defined up to conjugacy, on $V = V_\varphi \otimes V_\sigma \otimes S^N(\mathfrak{n})$ via $\gamma \mapsto \varphi(\gamma) \otimes \sigma(m_\gamma) \otimes \mathrm{Ad}^N((m_\gamma a_\gamma)^{-1})$. Let $\lambda_\sigma$ denote the infinitesimal character of $\sigma$ and let $c_\sigma = B(\lambda_\sigma) - B(\rho)$. The main result of chapter 4 is the **determinant formula**

$$Z_{\sigma,\varphi}(m+s) = \det(\triangle_{\sigma,\varphi} + \mathrm{s}^2 + \mathrm{c}_\sigma)\left(\det(\mathrm{D}_\sigma + \mathrm{s})\exp(\mathrm{P}_\sigma(\mathrm{s}^2))\right)^{2(-1)^{\mathrm{m}}\dim\varphi\chi(\mathrm{X}_\Gamma)/\chi(\mathrm{X}^{\mathrm{d}})}$$

where $\triangle_{\sigma,\varphi}$ is a virtual differential operator. $P_\sigma$ is a polynomial which can be calculated via 1.2 and the fact that $\log Z_{\sigma,\varphi}$ tends to zero as $s \to \infty$. This will be needed later.

Now fix any $K$-type $\tau \in \hat{K}$ then it is known [Zel] that $\tau$ reduces to $M$ with multiplicity one. We define $Z_{\tau,\varphi} = \prod_{\sigma \in \hat{M}} \sigma_{|\tau} Z_{\sigma,\varphi}$, where the symbol $\sigma \mid \tau$ means that $\sigma$ occurs in hermitian homogeneous vector bundle $E_\tau$ over $X$ which pushes down to a bundle $E_{\tau,\Gamma}$ over $X_\Gamma$. Let $\triangle_{\tau,\varphi}$ denote the Laplacian corresponding to $E_{\tau,\Gamma} \otimes \varphi$. Let $\tau_l$ denote the representation of $K$ on $\wedge^l \mathfrak{p}_+$ then it is easy to see that $\tau \mid_M = \sigma_l \oplus \sigma_{l-1}$ where $\sigma_l$ is the representation of $M$ on $\wedge^l \mathfrak{n}_+$. Let $\triangle_l = \sum_{n \geq 0}(-1)^{l+n}\triangle_{\tau_n}$. Where $\triangle_{\sigma_l}$ is defined as in chapter 4. Note that we have $\triangle_l = \triangle_{\sigma_l} +$ finite operator. More precisely we get

$$det(\triangle_l + \lambda) = det(\triangle_{\sigma_l} + \lambda)\prod_{n=0}^{l}\prod_{\pi \in \check{G}_d}(\pi(C) + \lambda)^{[\pi:\tau_n]\dim_{\Gamma}(\pi)(-1)^{l+n}}$$

a calculation now shows that only for $n = l = m$ we have a nontrivial contribution. Namely for the $\pi$ in the discrete series with Harish-Chandra parameter

$$\lambda_\pi = \rho_G$$

and therefore $\pi(C) = 0$. Further we have $\dim_\gamma(\pi) = (-1)^m \mathrm{Ar}(X_\Gamma)$, where $\mathrm{Ar}(X_\Gamma)$ denotes the arithmetic genus of $X_\Gamma$. From this we get

$$\triangle_l = \triangle_{\sigma_l} \quad \text{for } l < m \text{ and}$$

$$\det(\triangle_m + \lambda) = \det(\triangle_{\sigma_l} + \lambda)\lambda^{(-1)^m \mathrm{Ar}(X_\Gamma)}.$$

Further let $D_l = D_{\sigma_l}$, we are going to give this operator a geometric meaning relating to the dual space later.

Let $P = MAN$ the Langlands decomposition of a minimal parabolic $P$, fix an irreducible unitary representation $\xi$ of $M$ with infinitesimal character $\lambda_\xi$. Recall that we have for the Casimir eigenvalue on a principal series representation $\pi_{\xi,\nu}$ that $\pi_{\xi,\nu}(C) = B(\nu) + B(\lambda_\xi) - B(\rho)$. Now for $l \geq 0$ let $Z_\varphi^l$ denote the zeta function attached to the representation of $M$ on $\wedge^l \mathfrak{n}_-$. Let

$$Z_\varphi(s) = \prod_{l=0}^{m-1} Z_\varphi^l(s+l)^{(-1)^l}.$$

Then we easily see that

$$Z_\varphi(s) = \prod_{\substack{1 \neq [\gamma] \\ \text{prime}}} \prod_{N \geq 0} \det(1 - e^{-sl_\gamma}\varphi(\gamma) \otimes \mathrm{Ad}^N(\gamma^{-1} \mid \mathfrak{n}_r \oplus \mathfrak{n}_-))$$

**Proposition 5.3.1** *We have*

$$Z_\varphi(s) = \prod_{l=0}^{m-1} \left( \frac{\det(\triangle_l + (s+1-m)^2 - (l-m)^2)}{\det^{(2)}(\triangle_l + (s+1-m)^2 - (l-m)^2)^{\dim\varphi}} \right)^{(-1)^l}$$

$$= \prod_{l=0}^{m-1} \det(\triangle_l + (s+1-m)^2 - (l-m)^2)^{(-1)^l}$$

$$\left( \det(D_l + s + 1 - m)\exp(P_l((s+1-m)^2)) \right)^{(-1)^{m+1}\dim\varphi\chi(X_\Gamma)/\chi(X^d)}$$

*By the second identity $Z_\varphi(s)$ can be continued to a meromorphic function on the entire plane. Note that $\chi(X_\Gamma)/\chi(X^d) = \mathrm{Ar}(X_\Gamma)$, the arithmetic genus of $X_\Gamma$.*

We further consider

$$\tilde{Z}_\varphi(s) = \prod_{l=0}^{m-1} Z_\varphi^l(s+2m-l)^{(-1)^l}.$$

It follows that

$$\tilde{Z}_\varphi(s) = \prod_{\substack{1 \neq [\gamma] \\ \text{prime}}} \prod_{N \geq 0} \det(1 - e^{-(s+2)l_\gamma}\gamma \mid V)$$

where $V = V_\varphi \otimes \wedge^{m-1}\mathfrak{n}_+ \otimes (\mathfrak{n}_r \oplus \mathfrak{n}_-)$ and $\gamma$ acts on $V$ as $\varphi(\gamma) \otimes \mathrm{Ad}(\gamma^{-1}) \otimes \mathrm{Ad}^N(\gamma^{-1})$

**Proposition 5.3.2** *We have*

$$\tilde{Z}_\varphi(s) = \prod_{l=0}^{m-1} \left( \frac{\det(\triangle_l + (s-1+m)^2 - (l-m)^2)}{\det^{(2)}(\triangle_l + (s-1+m)^2 - (l-m)^2)^{\dim\varphi}} \right)^{(-1)^l}$$

$$= \prod_{l=0}^{m-1} \det(\triangle_l + (s-1+m)^2 - (l-m)^2)^{(-1)^l}$$

$$\left( \det(D_l + s - 1 + m)\exp(P_1((s-1+m)^2)) \right)^{(-1)^{m+1}\dim\varphi\chi(X_\Gamma)/\chi(X^d)}$$

We want to understand the determinant of the operator $D_l$ from above. It will turn out that it has a geometric interpretation as a virtual differential operator on the dual space.

Let $X^d$ denote the dual space to $X$. The decomposition $\mathfrak{p} = \mathfrak{p}_+ \oplus \mathfrak{p}_-$ establishes a complex structure on $X^d$ as well. Let $\triangle_{p,q}^d$ denote the Hodge Laplace operator on $(p,q)$-forms. Further define $\triangle_l^d$ by $\triangle_l^d = D_l^2 - (m-l)^2$.

**Theorem 5.3.1** *We have* $\det(\triangle_l^d + s)\det(\triangle_{l-1}^d + s) = \det(\triangle_{0,l}^d + s)R_l(s)$, *where* $R_l(s)$ *is the polynomial given by* $R_0(s) = 1$ *and*

$$R_l(s) = \prod_{\frac{m}{2} \leq n \leq \frac{m+l}{2}} (4n^2 - (m-l)^2 + s)^{Q_l(2n)}$$

*for* $l > 0$ *and* $l+m$ *even and*

$$R_l(s) = \prod_{m \leq 2n-1 \leq m+l} ((2n-1)^2 - (m-l)^2 + s)^{Q_l(2n-1)}$$

*for* $l > 0$, $l+m$ *odd.*

Proof of the theorem: Let $\tau_l \in \hat{K}$ be the representation of $K$ on $\wedge^l\mathfrak{p}_-$. Then the space of $(0,l)$-forms on $X^d$ is as $G$-module equal to $(C^\infty(G^d) \otimes \tau_l)^K$. Further we have $\tau_l \mid_M = \sigma_l \oplus \sigma_{l-1}$ where $\sigma_{-1} = 0$, and the Laplacian is given by the Casimir operator $C$. Now let $\lambda \in (\mathfrak{a} + \mathfrak{t_m})^*$ be a dominant integral weight and let $V(\lambda) \in \hat{G}^d$ denote the corresponding irreducible representation of the compact group $G^d$. The representations $\tau_l$ have the following remarkable property extending Theorem1.7 of [CaWo].

**Lemma. 5.3.1** *For $\lambda$ as above we have that $\tau_l$ occurs in the restriction of $V(\lambda)$ if and only if $\lambda \mid_{\mathfrak{t_m}}$ is the highest weight of $\sigma_l$ or $\sigma_{l-1}$.*

Proof: Lets show the only if part first: Assume $\tau_l$ occurs in $V(\lambda)$ and let $W \subset V(\lambda)$ denote the space of $\tau_l$. Write $U(\mathfrak{g}), U(\mathfrak{k}), U(\bar{\mathfrak{n}}), U(\mathfrak{a})$ for the universal enveloping algebras of $\mathfrak{g}, \mathfrak{k}, \bar{\mathfrak{n}}, \mathfrak{a}$, where $\bar{\mathfrak{n}}$ is the nilpotent Lie algebra corresponding to the negative roots. The Iwasawa decomposition $\mathfrak{g} = \bar{\mathfrak{n}} \oplus \mathfrak{a} \oplus \mathfrak{k}$ and the PBW theorem show $U(\mathfrak{g}) = U(\bar{\mathfrak{n}})U(\mathfrak{a})U(\mathfrak{k})$. Let $Pr$ denote the projection on $V(\lambda)$ to the highest weight space $V(\lambda)_\lambda$. Then $U(\bar{\mathfrak{n}})$ and $U(\mathfrak{a})$ commute with $Pr$ and we get $V(\lambda)_\lambda = PrV(\lambda) = Pr(U(\mathfrak{g})W) = Pr(U(\bar{\mathfrak{n}})U(\mathfrak{a})W) = U(\bar{\mathfrak{n}})U(\mathfrak{a})PrW$. Which yields $PrW \neq 0$ and so $\lambda \mid_{T_{\mathfrak{m}}}$ must occur as a weight of $\sigma_l$ or $\sigma_{l-1}$. But since $\lambda$ was dominant, $\lambda \mid_{\mathfrak{t_m}}$ will be M-dominant and an inspection shows that the only dominant weight of $\sigma_l$ is its highest weight. This gives the only if direction.

For the if part we will consider the underlying algebraic groups. We also change from the torus $A$ to a Cayley transform of $A$ lying in $K$. This doesn't change highest weight theory. So let for the moment $G = GL_{m+1}, K = GL_m \times GL_1$ embedded in $G$ in the obvious manner, let $M = GL_{m-1} \times GL_1$ embedded in $G$ as follows: $(A, a) \in M$ is mapped to $\mathrm{diag}(a, A, a)$, Further let $A = GL_1$ embedded as $a \mapsto \mathrm{diag}(a, 1, a^{-1})$. Under these embedding the natural ordering on the roots of $G$ will give our chosen ordering. Note that we also have added a center since we had to consider $SL_{m+1}$ otherwise. This simplifies argumentation and doesn't change the subject. Now $\tau_l$ and $\sigma_l$ are given by $\tau_l(k, c) = c^{-l}\lambda_l k$, the representation space is $\wedge^l \mathbb{C}^m$ and $\sigma_l(a, A, a) = a^{-l} \wedge^l A$, the representation space is $\wedge^l \mathbb{C}^{m-1}$. Let T denote the set of diagonal matrices in $G$. Any character on T will be of the form $\chi : \mathrm{diag}(a_1, \ldots, a_{m+1}) \mapsto a_1^{\chi_1} \ldots a_{m+1}^{\chi_{m+1}}$ with $\chi_i \in \mathbb{Z}$. We will write $\chi = (\chi_1, \ldots \chi_{m+1})$. The condition of dominance is $\chi_1 \geq \chi_2 \geq \ldots \geq \chi_{m+1}$. Now let $\lambda = (\frac{n-l}{2}, 1, \ldots, 1, 0, \ldots, 0, \frac{-n-l}{2})$ with l-ones after the $\frac{n-l}{2}$. The integrality gives $n \equiv l(2)$ and the dominance $n > l$ if $l > 0$ and $n \geq 0$ if l=0. We assume $n > 0$ and we want to show that $V(\lambda)$ contains $\tau_l$ and $\tau_{l+1}$. Here we want to make use of the Weyl principle which says that any irreducible representation of $GL_{m+n}$ splits, when reduced to $GL_m \times GL_n$ as a sum $\bigoplus_{\pi \in GL_n^{irr}} \pi' \otimes \pi$ where $\pi'$ is an irreducible representation of $GL_m$, each $\pi$ occurring maximally once. Now it is easy to see that $V(\lambda) \mid_{GL_1}$ contains the weight (-l) and that the corresponding highest weight of $K$ is $(1, \ldots, 1, 0, \ldots, 0, -l)$ with l-ones, i.e. the highest weight of $\tau_l$. In the same manner one exhibits $\tau_{l+1}$ as counterpart of (-l-1). The proof of the lemma is finished. ∎

Proof of the theorem (continued): The Peter-Weyl theorem shows $C^\infty(G^d) = \bigoplus_{\pi \in \hat{G}^a} \pi \otimes \pi^*$, as a $G^d \times G^d$ module. It is known that any $\pi \in \hat{G}^d$ restricts to $K$ with multiplicity one [Zel]. From this we see that, as a $G$-module

$$(C^\infty(G^d) \otimes \tau_l)^K = \bigoplus_{\substack{\pi \in \hat{G}^d \\ \tau_l \mid \pi}} \pi.$$

So that

$$\det(\triangle_{0,1}^d + s) = \prod_{\substack{\pi \in \hat{G}^d \\ \tau_l \mid \pi}} (\pi(C) + s)^{\dim \pi},$$

where the product means the zeta regularized product.

Now let $\lambda = (\frac{n-l}{2}, 1, \ldots, 1, 0, \ldots, 0, \frac{-n-l}{2})$ dominant integral and $\pi = \pi_\lambda$ the irreducible representation of $G^d$ with highest weight $\lambda$ then $\pi(C) = B(\lambda + \rho) - B(\rho) = B(\frac{n}{2}\alpha_r + \sigma_l + \rho) - B(\rho) = (m+n)^2 - (m-l)^2$ and for the dimension of $\pi$ we have by the Weyl character formula: $\dim \pi = Q_l(m+n)$. Note now that the element $\gamma_0$ is, in the algebraic setting, given by the matrix $\mathrm{diag}(-1, 1, -1)$ and so we get: $\sigma_l$ even $\Leftrightarrow$ $l + m$ even. Recall that for $\lambda = \sigma_l + \frac{n}{2}\alpha_r$ the dominance condition was $n > l$ so the right hand side of the equality in the theorem does contain $\det(\triangle_{0,1}^d + s)$ and all factors with $\lambda$ not dominant. This gives the assertion. ∎

We assemble the results of this section to

**Theorem 5.3.2** *Let $X$ denote a hermitian symmetric space of real rank one and let $\Gamma$ a torsionfree uniform lattice in the group $G$ of orientation preserving isometries. Let $X_\Gamma$ denote the quotient $\Gamma \backslash X$. Let $\varphi$ denote a finite dimensional unitary representation of $\Gamma$ and define the zeta function $Z_\varphi$ as:*

$$Z_\varphi(s) = \prod_{\substack{1 \neq [\gamma] \\ prime}} \prod_{N \geq 0} \det(1 - e^{-sl_\gamma}\varphi(\gamma) \otimes \mathrm{Ad}^N(\gamma^{-1} \mid \mathfrak{n}_r \oplus \mathfrak{n}_-)),$$

*then $Z_\varphi$ admits analytic continuation to a meromorphic function. The function $Z_\varphi(s)$ has vanishing order at zero:*

$$n_0 = \mathrm{ord}_{s=0}Z_\varphi(s) = (-1)^m \dim\varphi \mathrm{Ar}(X_\Gamma) + \sum_{q=0}^m q(-1)^q h_{0,q}(\varphi)$$

*where $h_{p,q}(\varphi)$ is the $(p,q)$-Hodge number of $X_\Gamma$ with $\varphi$-coefficients. The function $R_\varphi(s) = Z_\varphi(s)s^{-n_0}$ has the value:*

$$R_\varphi(0) = \frac{T(X_\Gamma, \varphi)}{T^{(2)}(X_\Gamma)^{\dim \varphi}}$$

*where $T$ is the Dolbeault torsion and $T^{(2)}$ the $L^2$-Dolbeault torsion. This can also be written*

$$R_\varphi(0) = T(X_\Gamma, \varphi)(T(X^d)E(m))^{(-1)^m \dim \varphi \mathrm{Ar}(X_\Gamma)}$$

*where $\mathrm{Ar}(X_\Gamma) = \sum_{q=0}^m (-1)^q h_{0,q}(X_\Gamma)$ denotes the arithmetic genus of $X_\Gamma$ and $E(m)$ is an elementary factor to be computed in section 4.9.*

**Corollary 5.3.1** *From the theorem we get a formula joining the $L^2$-torsion to the torsion of the dual space:*

$$T^{(2)}(X_\Gamma) = (T(X^d)E(m))^{(-1)^{m+1}\mathrm{Ar}(X_\Gamma)}.$$

For the proof we note that the first identities in Prop. 4.4 and 4.5 show that $Z_\varphi(0) = \tilde{Z}_\varphi(0)$. The second identities give

$$Z_\varphi(s)\tilde{Z}_\varphi(-s) = \prod_{l=0}^{m-1} \det(\triangle_l + (s+l-m)^2 - (l-m)^2)^{2(-1)^l}$$

$$\left(\det(D_1^2 - (s+1-m)^2)\exp(2P_1((s+1-m)^2))\right)^{(-1)^{m+1}\dim_\varphi\chi(X_\Gamma)/\chi(X^d)}.$$

Now Theorem 4.6 gives the assertion. ∎

We will give the recipe to compute the constant $E(m)$. We have

$$E(m) = \prod_{l=1}^{m-1} R_l(0)^{l(-1)^l}\exp(N_m),$$

where the $R_l$ are the polynomials of Theorem 3.6 and $N_m$ is given as follows: Consider the polynomial

$$F_l(x) = \binom{m-1}{l} \prod_{\substack{k=0 \\ k\neq l}}^{m-1} (x^2 - (m-2k-l)^2)$$

we get

$$Q_l(x) = \frac{xF_l(x)}{mF_0(m)} = \sum_{j=1}^{m} b_j^l x^{2j-1}.$$

Let $P_l(x) = \sum_{j=1}^{m} b_j^l c_j x^j$, where $c_j = -\frac{1}{j}\sum_{r=1}^{j} \frac{1}{2r-1}$, then

$$N_m = 2\sum_{l=0}^{m-1} (-1)^l P_l((m-l)^2).$$

The $P_l$'s are the polynomials occurring in th determinant formula in Proposition 3.5. The first values of $E$ are:

$$E(1) = \exp(-2),$$

$$E(2) = \exp(-1/2),$$

$$E(3) = 2^{-24}3^{-12}\exp(-\frac{739}{3\cdot 2^5}),$$

$$E(4) = 2^{-895}3^{-435}\exp(-\frac{11\cdot 277}{2\cdot 3^3\cdot 7}).$$

## 5.4   The Generalized Selberg Zeta Function

In this section we may assume $X$ to be an arbitrary globally symmetric space of the noncompact type. Fix a $\theta$-stable Cartan subgroup $H = AB$ of splitrank 1 and a parabolic subgroup $P$ with Langlands decomposition $P = MAN$. Let $G$ act on itself by conjugation, write $g.x = gxg^{-1}$, write $G.x$ for the orbit, so $G.x = \{gxg^{-1}|g \in G\}$ as well as $G.S = \{gsg^{-1}|s \in S, g \in G\}$ for any subset $S$ of $G$. We are going to consider functions that are supported on the set $G.(MA)$. Let $\sigma$ be a finite dimensional representation of $M$. Then let $f_\sigma$ be as in 5.1 but with $M$ taking the place of $G$. Now consider for $s \in \mathbb{C}$ and $n \in \mathbb{N}$ the function $g_s^n$ on $\bar{A}^+$ given by

$$g_s^n(a) = l_a^{n+1} e^{-sl_a},$$

where $l_a$ is the length of $a$, so $l_a = |\log(a)|$. Choose any $\eta : N \to G$ which has compact support, is positive and such that

$$\int_N \eta(n)dn = 1.$$

Given these data, let $\Phi = \Phi_{\eta,\sigma,n,s} : G \to \mathbb{C}$ be defined by

$$\Phi(knma(kn)^{-1}) = \eta(n)f_\sigma(m)\frac{g_s^n(a)}{\det(1 - (ma)^{-1}|\mathfrak{n})},$$

for $k \in K$, $n \in N$, $m \in M$ and $a \in A$. Further $\phi(g) = 0$ if $g$ is not in $G.(MA)$.

To fix notations let $(\varphi, V_\varphi)$ be a finite dimensional unitary representation of $\Gamma$ and recall that the group $G$ acts by right shifts unitarily on the Hilbert space $L^2(\Gamma\backslash G, \varphi)$ consisting of all measurable functions $f : G \to V_\varphi$ such that $f(\gamma x) = \varphi(\gamma)f(x)$ and $f$ is square integrable over $\Gamma\backslash G$ (modulo null functions). It is known that as a $G$-representation this space splits as

$$L^2(\Gamma\backslash G) = \bigoplus_{\pi \in \hat{G}} N_{\Gamma,\varphi}(\pi)\pi$$

with finite multiplicities $N_{\Gamma,\varphi} < \infty$.

**Proposition 5.4.1** *The function $\Phi$ is $(n - \dim(\mathfrak{n}))$-times continuously differentiable. For $n$ and $\mathrm{Re}\,(s)$ large enough it goes into the trace formula and for any finite dimensional unitary representation $\varphi$ of $\Gamma$ we have*

$$\sum_{\pi \in \hat{G}} N_{\Gamma,\varphi}(\pi)\mathrm{tr}(\Phi) = \sum_{[\gamma]} \mathrm{vol}(\Gamma_\gamma\backslash G_\gamma)\mathcal{O}_{m_\gamma}^M(f_\sigma)\frac{g_s^n(a)\mathrm{tr}(\varphi(\gamma))}{\det(1 - (ma)^{-1}|\mathfrak{n})},$$

*where the sum runs over all classes $[\gamma]$ such that $\gamma$ is conjugate to an element $m_\gamma a_\gamma$ of $MA^+$.*

**Proof:** The function $f_\sigma$ is rapidly decreasing [BM], replace for the moment $g_s^n$ by $g_s^n\chi_T$, where $\chi_T : \bar{A}^+ \to [0,1]$ is a $C^\infty$-function which is constant 1 for $l_a \leq T$ and constant 0 for $l_a \geq T+1$. Then $\Phi$ is rapidly decreasing and it thus goes into the trace formula. Considering each side of the trace formula separately it is easy to see that, as $T$ tends to infinity, both converge dominatedly. ∎

### 5.4.1   The generalized Selberg zeta function

Besides the parabolic $P = MAN$ we also consider the opposite parabolic $\bar{P} = MA\bar{N}$. The Lie algebra of $\bar{N}$ is written $\bar{\mathfrak{n}}$. Let $V$ denote a Harish-Chandra module of $G$ then we consider the Lie algebra homology $H_*(\bar{\mathfrak{n}}, V)$ and cohomology $H^*(\bar{\mathfrak{n}}, V)$. It is shown in [HeSchm] that these are Harish-Chandra modules of the group $MA$.

We will say that a discrete subgroup $\Gamma \subset G$ is **nice** if for every $\gamma \in \Gamma$ the adjoint $\mathrm{Ad}(\gamma)$ has no roots of unity as eigenvalues. Every arithmetic $\Gamma$ has a nice subgroup of finite index [Bor]. Let $H_1 \in A^+$ be the unique element with $B(H_1) = 1$.

We will denote by $\mathcal{E}_H(\Gamma)$ the set of nontrivial $\Gamma$-conjugacy classes $[\gamma]$, which are such that $\gamma$ is in $G$ conjugate to an element of H. Such an element will then be written $a_\gamma b_\gamma$ or $a_\gamma m_\gamma$.

**Theorem 5.4.1** *Let $\Gamma$ be nice and $(\varphi, V_\varphi)$ a finite dimensional unitary representation of $\Gamma$. Choose a $\theta$-stable Cartan $H$ of splitrank one. For* Re $(s) >> 0$ *define the zeta function*

$$Z_{H,\sigma,\varphi}(s) = \prod_{[\gamma] \in \mathcal{E}_H(\Gamma)} \prod_{N \geq 0} \det(1 - e^{-sl_\gamma}\gamma | V_\varphi \otimes W_\sigma \otimes S^N(\mathfrak{n}))^{\chi_1(X_\gamma)},$$

*where $\gamma$ acts on $V_\varphi \otimes W_\sigma \otimes S^N(\mathfrak{n})$ via $\varphi(\gamma) \otimes \sigma(m_\gamma) \otimes Ad^N((m_\gamma a_\gamma)^{-1})$. Then $Z_{H,\sigma,\varphi}$ has a meromorphic continuation to the entire plane. If $s_0$ is a pole or zero of $Z_{H,\sigma,\varphi}$ then $s_0 = \lambda(H_1)$, where $\lambda \in \mathfrak{a}^*$ is a weight of the $\mathfrak{a}$ action on the finite dimensional virtual space*

$$\sum_{p,q} (-1)^{p+q}(H^q(\bar{\mathfrak{n}}, \pi) \otimes \wedge^p\mathfrak{p}_M \otimes W_{\sigma^*})^{K_M},$$

*for some $\pi \in \hat{G}$ with $N_{\Gamma,\varphi}(\pi) \neq 0$. Here $\mathfrak{a}$ only acts on the first tensor factor. Let $m_{s_0,\pi,\sigma} \in \mathbb{Z}$ denote the dimension of the generalized weight space of $\lambda$ in this virtual module. For the order of $Z_\sigma$ at $s = s_0$ we have*

$$\mathrm{ord}_{s=s_0}(Z_{H,\sigma,\varphi}(s)) = -\sum_{\pi \in \hat{G}} N_{\Gamma,\varphi}(\pi) m_{s_0,\pi,\sigma}$$

*Remark.* To give the reader a better feeling of what this theorem says about the singularities of $Z_{H,\sigma,\varphi}$, we consider a zero or pole $s_0$ of $Z_{H,\sigma,\varphi}$. It then follows from Harish-Chandra's theory of global characters, that there is an integer $k \in \{0, \ldots, 2\dim(\mathfrak{n}_+ + 2)\}$ and a Weyl group element $w \in W(\mathfrak{h}, \mathfrak{g})$ such that

$$s_0 = w\lambda_\pi(H_1) - 2|\rho_0| - k|\alpha_r|/2,$$

where $\lambda_\pi \in \mathfrak{h}^*$ is the infinitesimal character of $\pi$.

Furthermore, there has to exist a finite subset $E \subset \Phi_{nI}^+(\mathfrak{h}, \mathfrak{g})$ of the set of nonimaginary positive roots such that

$$B(w\lambda_\pi|_\mathfrak{b} - \sum_{\alpha \in E} \alpha|_\mathfrak{b}) = B(\lambda_{\hat{\sigma}}).$$

**Proof:** From [D-Hitors] we take for $\Gamma$ nice

$$\chi_1(X_\gamma) = \frac{\mid W(\mathfrak{g}_\gamma, \mathfrak{h}) \mid \prod_{\alpha \in \Phi_\gamma^+}(\rho_\gamma, \alpha)}{\lambda_\gamma c_\gamma [G_\gamma : G_\gamma^0]} \mathrm{vol}(\Gamma_\gamma \backslash G_\gamma),$$

so that the geometric side of our trace formula will be

$$\sum_{[\gamma]} \frac{l_{\gamma_0} \chi_1(X_\gamma) \overline{\mathrm{tr}\sigma(m_\gamma)} l_\gamma^{n+1} e^{-sl_\gamma \mathrm{tr}(\varphi(\gamma))} \mathrm{tr}(\varphi(\gamma))}{\det(1 - \gamma^{-1}|\mathfrak{n})} = (-1)^{n+1} (\frac{\partial}{\partial s})^{n+2} \log Z_{H,\sigma,\varphi}(s).$$

To understand the spectral side let

$$(MA)^+ = \text{interior in } MA \text{ of the set } \{g \in MA | \det(1 - ga|\bar{\mathfrak{n}}) \geq 0 \,\, \forall a \in A^+\}.$$

Now $B$ is a compact Cartan in $M$ and $(M.B)A^+$ is a subset of $(MA)^+$. It was proven in [HeSchm] that for any $\pi \in \hat{G}$, writing the character of $\pi$ as $\Theta_\pi^G$ we have for $ma \in (MA)^+ \cap G^{\mathrm{reg}}$

$$\Theta_\pi^G(ma) = \frac{\Theta_{H_*(\bar{\mathfrak{n}},\pi)}^{MA}(ma)}{\det(1 - (ma)^{-1}|\mathfrak{n})}.$$

Let $f$ be supported on $G.(MA)$ then applying the Weyl integration formula twice gives

$$\int_G f(g)dg = \int_{G/MA} \int_{MA^+} f(gmag^{-1})|det(1 - ma|\mathfrak{n} \oplus \bar{\mathfrak{n}})dadmdg.$$

So that for $\pi \in \hat{G}$

$$\begin{aligned}
\mathrm{tr}\pi(\Phi) &= \int_G \Theta_\pi^G(x)\Phi(x)dx \\
&= -\int_{MA^+} \Theta_\pi^G(ma)f_\sigma(m)g_s^n(a)\det(1 - ma|\mathfrak{n})dadm \\
&= -\int_{A^+} \sum_L \frac{1}{|W(L,M)|} \int_L \int_{M/L} f_\sigma(mlm^{-1})dm\Theta_\pi^G(la)\det(1 - la|\mathfrak{n})dmdlg_s^n(a)da,
\end{aligned}$$

where the sum runs over all $M$-conjugacy classes of Cartan subgroups $L$ of $M$. The integral over $M/L$ is just an orbital integral of $f_\sigma$ so that in the sum only the compact Cartan $B$ will survive. Hence

$$\begin{aligned}
\mathrm{tr}\pi(\Phi) &= -\int_{MA^+} \Phi(ma)\det(1 - ma|\mathfrak{n})\Theta_{H_*(\bar{\mathfrak{n}},\pi)}^{MA}(ma)dmda \\
&= \int_{MA^+} \Phi(ma)a^{2\rho}\det(1 - (ma)^{-1}|\mathfrak{n})\Theta_{H_*(\bar{\mathfrak{n}},\pi)}^{MA}(ma)dmda \\
&= -\int_{MA^+} \Phi(ma)\det(1 - (ma)^{-1}|\mathfrak{n})\Theta_{H^*(\bar{\mathfrak{n}},\pi)}^{MA}(ma)dmda \\
&= -\int_{MA^+} f\sigma^*(m)\Theta_{H^*(\bar{\mathfrak{n}},\pi)}^{MA}(ma)dmg_s^n(a)da,
\end{aligned}$$

where we have used the isomorphism

$$H_p(\bar{\mathfrak{n}}, V) \cong H^p(\bar{\mathfrak{n}}, V) \otimes \wedge^{max}\bar{\mathfrak{n}},$$

([HeSchm], 2.18). This gives

$$\text{tr}\pi(\Phi) = -\int_{A^+} \text{tr}(a|(H^*(\bar{\mathfrak{n}}, \pi) \otimes \wedge^*\mathfrak{p}_M \otimes W_{\sigma^*})^{K_M})l_a^{n+1}e^{-sl_a}da.$$

Now let the $A$-decomposition of $(H^*(\bar{\mathfrak{n}}, \pi) \otimes \wedge^*\mathfrak{p}_M \otimes W_{\sigma^*})^{K_M}$ be

$$\bigoplus_{\lambda} E_\lambda,$$

where $E_\lambda$ is the biggest subspace on which $a - e^{\lambda(\log a)}$ acts nilpotently. Let $m_\lambda$ denote the dimension of $E_\lambda$. Integrating over $\mathfrak{a}^+$ instead over $A^+$ we get

$$
\begin{aligned}
\text{tr}\pi(\Phi) &= -\int_0^\infty \sum_\lambda m_\lambda e^{(\lambda(H_1)-s)t} t^{n+1} dt \\
&= (-1)^{n+1}(\frac{\partial}{\partial s})^{n+1} \sum_\lambda m_\lambda \frac{1}{\lambda(H_1) - s} \\
&= -\sum_\lambda m_\lambda \frac{1}{(s - \lambda(H_1))^{n+2}}.
\end{aligned}
$$

From this the theorem follows. ∎

The remark following the theorem will follow form the Wigner Lemma [BorWall], Cor 4.2 together with the

**Proposition 5.4.2** *Let $\sigma$ be a finite dimensional representation of $M$ then the order of $Z_{H,\sigma,\varphi}(s)$ at $s = s_0$ is*

$$(-1)^{\dim(N)} \sum_{\pi \in \hat{G}} N_{\Gamma,\varphi}(\pi) \sum_{q=0}^{\dim(\mathfrak{m}\oplus\mathfrak{n})} (-1)^q \dim(H^q(\mathfrak{m} \oplus \mathfrak{n}, K_M, \pi \otimes W_{\sigma^*})_\lambda$$

*with $\lambda = s_0\lambda_1$.*

**Proof:** Extend $W_{\sigma^*}$ to a $\mathfrak{m} \oplus \mathfrak{n}$-module by letting $\mathfrak{n}$ act trivially. We then get

$$H^p(\mathfrak{n}, \pi) \otimes W_{\sigma^*} \cong H^p(\mathfrak{n}, \pi \otimes W_{\sigma^*}).$$

The $(\mathfrak{m}, K_M)$-cohomology of the module $H^p(\mathfrak{n}, \pi \otimes W_{\sigma^*})$ is the cohomology of the complex $(C^*)$ with

$$
\begin{aligned}
C^q &= \text{Hom}_{K_M}(\wedge^q\mathfrak{p}_M, H^p(\mathfrak{n}, \pi) \otimes W_{\sigma^*}) \\
&= (\wedge^q\mathfrak{p}_M \otimes H^p(\mathfrak{n}, \pi) \otimes W_{\sigma^*})^{K_M},
\end{aligned}
$$

since $\wedge^p \mathfrak{p}_M$ is a self-dual $K_M$-module. Therefore we have an isomorphism of virtual $A$-modules:

$$\sum_q (-1)^q (H^p(\mathfrak{n}, \pi) \otimes \wedge^q \mathfrak{p}_M \otimes W_{\sigma^*})^{K_M} \cong \sum_q (-1)^q H^q(\mathfrak{m}, K_M, H^p(\mathfrak{n}, \pi \otimes W_{\sigma^*})).$$

Now one considers the Hochschild-Serre spectral sequence in the relative case for the exact sequence of Lie algebras

$$0 \to \mathfrak{n} \to \mathfrak{m} \oplus \mathfrak{n} \to \mathfrak{m} \to 0$$

and the $(\mathfrak{m} \oplus \mathfrak{n}, K_M)$-module $\pi \otimes W_{\sigma^*}$. We have

$$E_2^{p,q} = H^q(\mathfrak{m}, K_M, H^p(\mathfrak{n}, \pi \otimes W_{\sigma^*}))$$

and

$$E_\infty^{p,q} = \mathrm{Gr}^q(H^{p+q}(\mathfrak{m} \oplus \mathfrak{n}, K_M, \pi \otimes W_{\sigma^*})).$$

Now the module in question is just

$$\chi(E_2) = \sum_{p,q} (-1)^{p+q} E_2^{p,q}.$$

Since the differentials in the Hochschild-Serre spectral sequence are $A$-homomorphisms this equals $\chi(E_\infty)$. So we get an $A$-module isomorphism of virtual $A$-modules

$$\sum_{p,q} (-1)^{p+q} (H^p(\mathfrak{n}, \pi) \otimes \wedge^q \mathfrak{p}_M \otimes W_{\sigma^*})^{K_M} \cong \sum_j (-1)^j H^j(\mathfrak{m} \oplus \mathfrak{n}, K_M, \pi \otimes W_{\sigma^*}). \blacksquare$$

### 5.4.2   The generalized Ruelle Zeta Function

**Theorem 5.4.2** *Let $\Gamma$ be nice and choose a $\theta$-stable Cartan subgroup $H$ of splitrank one. For* Re $(s) >> 0$ *define the zeta function*

$$Z_{H,\sigma,\varphi}^R(s) = \prod_{[\gamma] \in \mathcal{E}_H(\Gamma)} \det(1 - e^{-sl_\gamma} \varphi(\gamma)),$$

*then $Z_{H,\sigma,\varphi}^R(s)$ extends to a meromorphic function on $\mathbb{C}$.*

**Proof:** We consider the logarithm of the Selberg function for $\varphi = 1$

$$\log Z_{H,\sigma,\varphi}^R(s) = -\sum_{[\gamma]} \frac{\chi_1(X_\gamma) e^{-sl_\gamma} \mathrm{tr}\sigma(m_\gamma) \mathrm{tr}(\varphi(\gamma))}{\mu_\gamma \det(1 - \gamma^{-1}|\mathfrak{n})}.$$

Consider first the case when there is only one positive root $\alpha_r$ in the root system $\Phi(\mathfrak{a}, \mathfrak{g})$, we normalize $B$ such that $|\alpha_r| = 2$. Then $a$ acts on $\wedge^l \bar{\mathfrak{n}}$ by $a^{-l\alpha_r}$. Let $\sigma_l$ denote the representation of $K_M$ on $\wedge^l \bar{\mathfrak{n}}$ the we get

$$Z_{H,\sigma,\varphi}^R(s) = \prod_{l=0}^{\dim N} Z_{H,\sigma_l,\varphi}(s + 2l)^{(-1)^l}$$

and hence the claim.

If we have two positive roots, say $\alpha_r$ and $\frac{\alpha_r}{2}$, with $|\alpha_r| = 2$ then we have $\bar{\mathfrak{n}} = \bar{\mathfrak{n}}_r \oplus \bar{\mathfrak{n}}_I$ where $\bar{\mathfrak{n}}_r$ has dimension one and $a$ acts on $\bar{\mathfrak{n}}_r$ by $a^{-\alpha_r}$ and on $\bar{\mathfrak{n}}_I$ by $a^{-\frac{\alpha_r}{2}}$. Now it follows

$$Z^R_{H,\sigma,\varphi}(s) = \prod_{l=0}^{\dim N - 1} \left( \frac{Z_{H,\wedge^l \bar{\mathfrak{n}}_I,\varphi}(s+l)}{Z_{H,\wedge^l \bar{\mathfrak{n}}_I \otimes \bar{\mathfrak{n}}_r,\varphi}(s+l+2)} \right)^{(-1)^l}. \blacksquare$$

## 5.5   The holomorphic torsion zeta function

Now let $X$ be hermitian again and let $(\tau, V_\tau)$ be an irreducible unitary representation of $K_M$.

**Theorem 5.5.1** *Let $\Gamma$ be nice and $(\varphi, V_\varphi)$ a finite dimensional unitary representation of $\Gamma$. Choose a $\theta$-stable Cartan $H$ of splitrank one with $c(H) = 2$. For Re $(s) >> 0$ define the zeta function*

$$Z_{H,\tau,\varphi}(s) = \prod_{[\gamma] \in \mathcal{E}_H(\Gamma)} \prod_{k \geq 0} \det(1 - e^{-sl_\gamma} \gamma | V^\tau_k)^{\chi_1(X_\gamma)},$$

*where $V^\tau_k = V_\varphi \otimes V_\tau \otimes S^k(\mathfrak{n})$ and $\gamma$ acts via $\varphi(\gamma) \otimes \tau(m_\gamma) \otimes Ad^k(\gamma^{-1})$. Then $Z_{H,\tau,\varphi}$ has a meromorphic continuation to the entire plane. If $s_0$ is a pole or zero of $Z_{H,\tau,\varphi}$ then $s_0 = \lambda(H_1)$, where $\lambda \in \mathfrak{a}^*$ is a weight of the $\mathfrak{a}$ action on the finite dimensional virtual space*

$$\sum_{p,q} (-1)^{p+q} (H^q(\bar{\mathfrak{n}}, \pi) \otimes \wedge^p \mathfrak{p}_{M,-} \otimes W_{\tau^*})^{K_M},$$

*for some $\pi \in \hat{G}$ with $N_{\Gamma,\varphi}(\pi) \neq 0$. Here $\mathfrak{a}$ only acts on the first tensor factor. Let $m_{s_0,\pi,\sigma} \in \mathbb{Z}$ denote the dimension of the generalized weight space of $\lambda$ in this virtual module. For the order of $Z_\sigma$ at $s = s_0$ we have*

$$\text{ord}_{s=s_0}(Z_\sigma(s)) = -\sum_{\pi \in \hat{G}} N_{\Gamma,\varphi}(\pi) m_{s_0,\pi,\sigma}$$

**Proof:** The proof proceeds as the proof of Theorem 5.4.1 with $g_{\tau^1}$ instead of $f_\sigma$, where $\tau^1$ denotes the virtual representation $\tau \otimes \wedge^* \mathfrak{p}_{M,-}$. One makes use of an analogue of Lemma 4.9 in [MS-2]. $\blacksquare$

Extend the definition of $Z_{H,\tau,\varphi}(s)$ to arbitrary virtual representations in the following way. Consider a finite dimensional virtual representation $\xi = \oplus_i a_i \tau_i$ with $a_i \in \mathbb{Z}$ and $\tau_i \in \hat{K}_M$. Then let $Z_{H,\xi,\varphi}(s) = \prod_i Z_{H,\tau_i,\varphi}(s)^{a_i}$.

Now for $c(H) = 2$ let $Z_{H,l,\varphi}(s) = Z_{H,\wedge^l \mathfrak{n}_-,\varphi}(s)$ and

$$Z_{H,\varphi}(s) = \prod_{l=0}^{\dim(\mathfrak{n}_-)} Z_{H,l,\varphi}(s + d_l(H)),$$

In the case $c(H) = 1$ let

$$Z_{H,\varphi}(s) = Z^R_{H,\varphi}(s + d(H)),$$

where $Z^R_{H,\varphi}$ denotes the Ruelle zeta function. See 5.1 for the constants.

**Proposition 5.5.1** *Assume $\Gamma$ is nice, then for $\lambda > 0$ we have the identity*

$$\prod_{q=0}^{\dim_{\mathbb{C}} X} \left( \frac{\det(\triangle_{0,q,\varphi} + \lambda)}{\det^{(2)}(\triangle_{0,q,\varphi} + \lambda)} \right)^{q(-1)^q} =$$

$$\prod_{\substack{H/\text{conj.} \\ \text{splitrank} = 1 \\ c(H) = 2}} \prod_{l=0}^{\dim \mathfrak{n}_+} \left( Z_{H,l,\varphi}(b_0(H) + \sqrt{\lambda + d_l(H)^2}) \right)^{(-1)^l}$$

$$\prod_{\substack{H/\text{conj.} \\ \text{splitrank} = 1 \\ c(H) = 1}} Z^R_{H,\varphi}(b_0(H) + \sqrt{\lambda + d(H)^2}).$$

    **Proof:** The equality is gotten by taking the Mellin transform of the expressions in Theorem 5.1.1. Then one uses the fact that the identity contribution to the trace formula equals the $\Gamma$-trace. ■

    Now the results in this section assemble to

**Theorem 5.5.2** *The zeta function $Z_{H,\varphi}$ extends to a meromorphic function on the entire plane. Let*

$$Z_\varphi(s) = \prod_{\substack{H/\text{conj.} \\ \text{splitrank} = 1}} Z_{H,\varphi}(s + b_0(H)).$$

*Let $n_0$ be the order of $Z_\varphi$ at zero then Proposition 5.5.1 shows that*

$$n_0 = \sum_{q=0}^{\dim_{\mathbb{C}}(X)} q(-1)^q (h_{0,q}(X_\Gamma) - h_{0,q}^{(2)}(X_\Gamma)),$$

*where $h_{p,q}(X_\Gamma)$ is the $(p,q)$-th Hodge number of $X_\Gamma$ and $h_{p,q}^{(2)}(X_\Gamma)$ is the $(p,q)$-th $L^2$-Hodge number of $X_\Gamma$. Let $R_\varphi(s) = Z_\varphi(s)s^{-n_0}$ then*

$$R_\varphi(0) = \frac{T(X_\Gamma, \varphi)}{T^{(2)}(X_\Gamma)^{\dim \varphi}}. ■$$

# Chapter 6

# Higher Fundamental Rank

In this chapter we consider Shimura manifolds $X_\Gamma = \Gamma \backslash G/K$, where $G$ is supposed to have a fundamental rank $> 0$. Examples are $G = SL_n(\mathbb{R})$, where $n$ is $> 2$.

This is the situation where our definitions of higher torsion and higher Euler characteristics start to unfold their meaning. It turns out that all torsion numbers of index smaller than the fundamental rank are trivial. This is paralleled on the geometrical side by the vanishing of the corresponding Euler numbers.

To define the zeta functions we have to impose a condition on the group $\Gamma$ which is automatically satisfied in the cases when the fundamental rank is either zero or one.

## 6.1 "Good" Groups

In this chapter we will assume $\Gamma$ to be **good** in a sense to be explained now: A real reductive group H is called **cuspidal**, if $H/$center has a compact Cartan subgroup. An element $g$ of $G$ is called cuspidal, if the centralizer $G_g$ of $g$ in $G$ is cuspidal. Note that every semisimple regular element of $G$ is cuspidal. Further an arbitrary semisimple element $g$ lies in a Cartan subgroup $H = AB$ where $A$ is the split part and $B$ is compact. According to this decomposition we write $g = a_g b_g$. Now suppose $H$ is chosen such that the dimension of $A$ is minimal then $g$ is cuspidal if and only if $a_g$ is regular in $A$, which means that the centralizer $G_{a_g}$ of $a_g$ equals the centralizer $G_A$ of $A$.

An element $g$ of $G$ is called **fundamental**, if $g \neq 1$ and $g$ lies in a fundamental Cartan subgroup. Now the lattice $\Gamma$ is called **good** if every nontrivial fundamental element of $\Gamma$ is cuspidal.

Note that goodness only depends on the commensurability class of $\Gamma$, i.e. when $\Gamma_1$ and $\Gamma_2$ are commensurable ($\Leftrightarrow \Gamma_1 \cap \Gamma_2$ has finite index in $\Gamma_1$ and in $\Gamma_1$) and torsionfree and $\Gamma_1$ is good then so is $\Gamma_2$.

If $G$ has a compact Cartan subgroup, i.e. $G$ has fundamental rank zero then every lattice is good. If $G$ has fundamental rank one then every torsionfree lattice is good.

If $G = G_1 \times G_2$, $\Gamma = \Gamma_1 \times \Gamma_2$ and $G_1$, $G_2$ both have no compact Cartan subgroup then $\Gamma$ is not good. In this case $\Gamma$ is reducible. There are however examples examples of irreducible lattices which fail to be good.

Let $p$ denote an odd prime and $G = SL_p(\mathbb{R})$ then any lattice in $G$ is arithmetic [Marg]. So let $A$ denote a central $\mathbb{Q}$-division algebra of degree $p$ [Pierce]. The $A \otimes \mathbb{R} \cong Mat_p(\mathbb{R})$ and the norm-1-subgroup $A^1$ of $A^*$ maps into $G$. Any arithmetic subgroup of $A^1$ gives a uniform lattice in $G$. Let $\alpha$ denote a nontrivial element of $A^1$. The subfield $\mathbb{Q}(\alpha)$ of $A$ must be a strictly maximal subfield of $A$. The centralizer of $\alpha$ is a simple algebra over $\mathbb{Q}(\alpha)$ inside $A$, hence it must coincide with $\mathbb{Q}(\alpha)$. Thus the centralizer of $\alpha$ is abelian, so $\alpha$ is regular and $\Gamma$ is good.

Let $G = SL_2(\mathbb{C})^{d_1} \times SL_2(\mathbb{R})^{d_2}$, $d_1 + d_2 > 1$ and $\Gamma$ an irreducible lattice in $G$. Then there is a number field $K$ with $d_1$ pairs of conjugate complex places and $d_2$ real places such that $\Gamma$ is commensurable to $A^1$, where $A$ is a quaternion algebra over $K$. We may assume $\Gamma$ nice. Since on the one hand $SL_2(\mathbb{R})$ has a compact Cartan and on the other every noncentral semisimple element of $SL_2(\mathbb{C})$ is regular we conclude that $\Gamma$ is good.

An example of a lattice which has nonregular elements but is good can be constructed similarly in $SL_4(\mathbb{R})$ and an example of a nongood irreducible lattice in $SL_9(\mathbb{R})$.

## 6.2 The Heat Trace

Now fix a finite dimensional unitary representation $\varphi$ of $\Gamma$. Then $\varphi$ defines a flat hermitian vector bundle $E_\varphi$ over $X_\Gamma$. Let $E$ denote the pullback of $E_\varphi$ to $X$. Let $\triangle_{\varphi,k,\Gamma}$ be the Laplace operator on $E_\varphi$-valued differential forms and $\triangle_{\varphi,k}$ the Laplace operator on $E$-valued k-forms on $X$. The heat operator $e^{-t\triangle_{\varphi,k}}$ has a smooth kernel $h_t^k$, which is of rapid decay [BM]. We put

$$f_t = \sum_{k=0}^{dim\, X} (-1)^k \left( \begin{array}{c} k \\ r \end{array} \right) \mathrm{tr}(h_t^k),$$

where $r = \mathrm{rank}(G) - \mathrm{rank}(k)$ denotes the fundamental rank of $G$ and the trace means pointwise trace. Since $\triangle_{\varphi,k,\Gamma}$ is the pushdown of $\triangle_{\varphi,k}$ we get by the Selberg trace formula with $d = \dim X$:

$$\sum_{k=0}^{d} (-1)^k \left( \begin{array}{c} k \\ r \end{array} \right) \mathrm{tr}(e^{-t\triangle_{\varphi,k}}) = \sum_{[\gamma]} \mathrm{tr}\varphi(\gamma)\, \mathrm{vol}(\Gamma_\gamma \backslash G_\gamma)\, \mathcal{O}_\gamma(f_t),$$

where the sum on the right hand side is taken over the conjugacy classes $[\gamma]$ in $\Gamma$ and $\mathcal{O}_\gamma(f_t)$ denotes the orbital integral:

$$\mathcal{O}_\gamma(f_t) = \int_{G/G_\gamma} f_t(g\gamma g^{-1})dg.$$

We want to use the Fourier transform formula for orbital integrals as in [HC-S]. Since the fundamental rank of G is $> 0$ there is no compact Cartan subgroup in $G$ and hence it suffices to compute the principal series contributions to the Fourier transform. So let H denote a $\Theta$-stable Cartan subgroup, let $A$ denote its split component and $B$ the compact component then we have a product $H = AB$. Let P denote a cuspidal parabolic subgroup with Langlands decomposition $P = MAN$. Let $(\xi, W_\xi)$ be an irreducible unitary representation of M and $e^\nu$ a quasicharacter of $A$. Set $\pi_{\xi,\nu} = Ind_P^G(\xi \otimes e^\nu \otimes 1)$ the principal series representation.

Similar to Proposition 5.1.1 we get:

**Proposition 6.2.1** *For the trace of $f_t$ under $\pi_{\xi,\nu}$ we have:*

$$\text{tr } \pi_{\xi,\nu}(f_t) = \begin{cases} -e^{t\pi\xi,\nu(C)} \dim(W_\xi \otimes \wedge^* \mathfrak{p}_\mathfrak{h})^{K \cap M} & \text{if } H \text{ is fundamental} \\ 0 & \text{otherwise.} \end{cases}$$

*Here C is the Casimir operator of G and for a finite dimensional vector space V we write*

$$\wedge^* V = \sum_{k=0}^{dim\ V} (-1)^k \wedge^k V.$$

*Further $\mathfrak{p}_\mathfrak{h}$ denotes the orthocomplement of $\mathfrak{a}$ in $\mathfrak{p}$, where $\mathfrak{a} = Lie A$.*

From now on let $H = AB$ be a fundamental $\Theta$-stable Cartan subgroup of G. Since $\exp : \mathfrak{a} \to A$ is an isomorphism we will frequently identify these groups to simplify notation.

As in 5.1 it follows that for an irreducible unitary representation $\xi$ of $M$ with infinitesimal character $\wedge_\xi$ we have

$$(W_\xi \otimes \wedge^* \mathfrak{p}_M \otimes V_\mu)^{K \cap M} \neq 0 \Rightarrow B(\wedge_\xi) = B(\mu).$$

Choose a system of positive roots according to the chosen parabolic $P = MAN$. Write $\rho$ for the half sum of the positive roots. Recall that we have

$$\pi_{\xi,\nu}(C) = B(\nu) + B(\wedge_\xi) - B(\rho)$$

Define $'F_f^H$ as in [HC-HA1,sec.17]. By Theorem 15 in [HC-S] and the preceding proposition we get with $A^*, B^*$ denoting the unitary duals:

$$'F_{f_t}^H(h) = -|W(g \mid H)|^{-1} \sum_{s \in W(G|H)} \sum_{B^* \in B^*} \epsilon_I(s)$$

$$\int_{A^*} B^*(s(B)) A^*(s(A)) e^{t(B(B^*) + B(A^*) - B(\rho))} dim(H(B^*) \otimes \wedge^* \mathfrak{p}_M \otimes \wedge^* \mathfrak{n})^{K \cap M^0} dh_R^*.$$

Here $H(B^*)$ means the discrete series representation of $M^0$ given by $B^*$. The normalizations of measures follow [HC-S]. We have

$$\int_{A^*} A^*(s(A)) e^{tB(A^*)} dA^* = \frac{e^{-B(A)/4t}}{\sqrt{4\pi t}}.$$

Let $C_M$ be the Casimir operator of M with respect to $B\mid_M$. Let $\rho_M$ denote the half sum of positive imaginary roots. Further let

$$'\triangle_M = \prod_{\alpha \in \Phi_M^+} (1 - \xi_{-\alpha}) = {}'\triangle_{M,n}\ {}'\triangle_{M,c}$$

the c and the n denoting the product over the compact and noncompact imaginary roots and $\xi_\alpha(h) = e^{\alpha(log\ h)}$ for $h \in H$. Now by Lemma 4.2 in [MS2] and similar computations as in the proof of Lemma 4.3 in [MS1] we get for $h = ab$ regular in H, where $\mathfrak{n} = Lie\ N \otimes \mathbb{C}$:

$$'F_{f_t}^H(h) = \frac{e^{-\frac{B(a)}{4t} - t(B(\rho) - B(\rho_M))}}{\sqrt{4\pi t}'} \frac{'\triangle_M(b)}{|'\triangle_{M,n}(b)|^2} tr(e^{tC_M}b \mid \wedge^*\mathfrak{n}).$$

Now let $h \in H$ and write $\Phi_h$ for the set of roots $\alpha$ such that $\xi_\alpha(h) = 1$. Write $\Phi_{h,c}$ and $\Phi_{h,n}$ for the compact and noncompact imaginary roots in $\Phi_h$. Let $c_h$ have the same meaning as in 3.2. Write as an $M^0$-module:

$$\wedge^l\mathfrak{n} = \sum_{i \in I_l} a_{i,l} V_{\lambda_{i,l}},$$

where $V_{\lambda_{i,l}}$ is an irreducible $M^0$-module with infinitesimal character $\lambda_{i,l}$ and $a_{i,l}$ denotes the multiplicity of $V_{\lambda_{i,l}}$ in $\wedge^l\mathfrak{n}$. Let $\tilde{\omega}$ denote the following element of the symmetric algebra $S(\mathfrak{h})$ over $\mathfrak{h}$:

$$\tilde{\omega} = \prod_{\alpha \in \Phi_h^+} H_\alpha.$$

Similar to [MS2] we get for the orbital integrals for $h \in H$, $h \neq 1$:

$$\int_{G/G_\gamma} f_t(xhx^{-1})dx = \int_{G/G_h^0} f_t(xhx^{-1})dx\ [G_h : G_h^0]^{-1}$$

$$= [G_h : G_h^0]^{-1} c_h^{-1} \frac{e^{-(B(A)/(4t)) - tB(\rho)}}{\xi_\rho(h) \prod_{\alpha \in \Phi_{h,c}^+}(1 - \xi_{-\alpha}(h))}$$

$$\sum_{l=0}^{dimN} (-1)^l \sum_{i \in I_l} a_{i,l} e^{tB(\lambda_{i,l})} \sum_{w \in W_{\mathbb{C}}(M)} \epsilon(w)\tilde{\omega}_h(w(\lambda_{i,l}))\xi_{w(\lambda_{i,l})}(B).$$

We write $\mathcal{E}_2(M^0)$ for the set of discrete series representations of $M^0$ and $\mathcal{E}_f(\Gamma)$ for the set of conjugacy classes $[\gamma]$ in $\Gamma$ such that $\gamma$ is fundamental in G. By the fact that $f_t$ is a pseudocuspform we conclude that in the trace formula only conjugacy classes in $\mathcal{E}_f$ contribute. We write $l_\gamma$ for the common length of the geodesics in the homotopy class of $\gamma$.

We end up with

**Lemma. 6.2.1**

$$\Theta(t) = \sum_{k=0}^{d} (-1)^k \left( \begin{array}{c} k \\ r \end{array} \right) tr e^{-t\triangle_{\varphi,k}}$$

$$= \sum_{[1]\neq[\gamma]\in\mathcal{E}_f(\gamma)} tr\varphi(\gamma) \frac{vol(\Gamma_\gamma\backslash G_\gamma)[G_h : G_h^0]^{-1}}{\xi_\rho(\gamma)\prod_{\alpha\in\Phi^+-\Phi_{\gamma,c}^+}(1-\xi_{-\alpha}(\gamma))} c_\gamma^{-1}$$

$$\frac{e^{-(l_\gamma^2/4t)-tB(\rho)}}{\sqrt{4\pi t}^r} \sum_{l=0}^{dim\ N} (-1)^l \sum_{i\in I_l} a_{i,t} e^{tB(\lambda_{i,l})} \sum_{w\in W_\mathbb{C}(M)} \epsilon(w)\tilde{\omega}_h(w(\lambda_{i,l}))\xi_{w(\lambda_{i,l})}(B)$$

$$+ dim\ F\ vol(X) \sum_{l=0}^{dim\ N} (-1)^l \sum_{i\in I_l} a_{i,l} \sum_{\xi\in\mathcal{E}_2(M^0)} \dim(W_\xi\otimes\wedge^*\mathfrak{p}_M\otimes V_{\lambda_{i,l}})^{K\cap M^0} F_\xi(t).$$

Here we wrote

$$F_\xi(t) = \int_{A^*} e^{t(B(A^*)+B(\lambda_{i,l})+B(\rho))}$$

and $P_\xi(\nu)$ is the Plancherel density of the principal series representations $\pi_{\xi,\nu}$. By [HC-HAIII] we know that $P_\xi$ is a polynomial of degree $\dim N$, more precisely

$$P_\xi(\nu) = c \prod_{\alpha\in\Phi^+} <\lambda_\xi + \nu, \alpha>,$$

where $c > 0$ is given by Theorem 24.1 in [HC-HAIII].

Choosing an orthonormal basis of $A^*$ we are lead to consider

$$\int_{\mathbb{R}^r} r^{-t\|x\|} P_\xi(x)dx$$

for a polynomial $P_\xi$. Writing $P_\xi$ as linear combination of monomials $x^\beta = x_1^{\beta_1}\ldots x_r^{\beta_r}$ it is clear that only those monomials contribute, for which $\beta \in 2\mathbb{Z}$. But for $\beta = 2\alpha$ we get

$$\int_{\mathbb{R}^r} e^{-t\|x\|} x^{2\alpha} dx = \sqrt{\pi}^r 2^{-|\alpha|} \prod_{j=1}^{r}(2\alpha_j-1)!! t^{-|\alpha|-\frac{r}{2}},$$

where $(2n-1)!! = (2m-1)(2m-3)\ldots 1$ and $|\alpha| = \alpha_1 + \ldots + \alpha_r$. We conclude, that $F_\xi(t)$ is $t^{-\frac{r}{2}}$ times a polynomial in $(\frac{1}{t})$.

## 6.3 Decomposition of The Lifted Heat Operator

Recall that for the space $\Omega^k$ of smooth k-forms on X with values in E we have

$$\Omega^k \cong (C^\infty(\Gamma\backslash G, \varphi) \otimes \wedge^k\mathfrak{p})^K.$$

Let $\mathfrak{p}_M$ denote the positive part of the Cartan decomposition of M. Write $\mathfrak{c} = \mathfrak{a} + \mathfrak{p}_M$. As $K \cap M^0$-module we have $\mathfrak{p} \cong \mathfrak{c} + \mathfrak{n}$. Let $V(\lambda_{i,l})$ denote the $V_{\lambda_{i,l}}$-isotypic subspace of $\wedge^l\mathfrak{n}$. Let $Pr^{i,j,l} : \wedge^k\mathfrak{p} \to \wedge^j\mathfrak{c} \otimes V(\lambda_{i,l})$ denote the projection and $I^{i,j,l} : \wedge^j\mathfrak{c} \otimes V(\lambda_{i,l}) \to \wedge^k\mathfrak{p}$ the injection. Let

$$\Omega^{i,j,l} = 1 \otimes Pr^{i,j,l}(\Omega^{i+j}) \subset (C^\infty(\Gamma\backslash G, \varphi) \otimes \wedge^j\mathfrak{c} \otimes V(\lambda_{i,l}))^{K\cap M^0}.$$

The space on the left hand side is only in the $K \cap M^0$-invariants and not in the K-invariants, but $e^{-t\triangle_{\varphi,k}}$ is a convolution operator and may be regarded as an operator from $K \cap M^0$-to K-invariants. So we may set

$$h_t^{i,j,l} = (1 \otimes Pr^{i,j,l}) \circ e^{-t\triangle_{\varphi,k}} \circ (1 \otimes I^{i,j,l}).$$

This defines an integral operator $e^{-t\triangle_\varphi^{i,j,l}}$.

**Lemma. 6.3.1** *As $K \cap M^0$-modules we have*

$$\sum_{j=0}^{dim\ \mathfrak{p}_M+r} (-1)^{j+l} \begin{pmatrix} j+l \\ r \end{pmatrix} (\wedge^j\mathfrak{c} \otimes \wedge^l\mathfrak{n}) \cong (-1)^{l+r}(\wedge^*\mathfrak{p}_M \otimes \wedge^l\mathfrak{n}).$$

The proof is straightforward.

Set $f_t^{i,j,l} = tr\ h_t^{i,j,l}$ and

$$f_t^{i,l} = \sum_{j=0}^{dim\ \mathfrak{p}_M} (-1)^j f_t^{i,j,l}.$$

Similar to 3.5 we find that $f_t^{i,l}$ is a pseudocuspform and for the fundamental principal series we have

$$tr\ \pi_{\xi,\nu}(f_t^{i,l}) = -e^{-t\pi_{\xi,\nu}(C)} dim(W_\xi \otimes \wedge^*\mathfrak{p}_M \otimes V(\lambda_{i,l}))^{K\cap M^0}.$$

So put $\Theta_{i,l}(t) = \sum_j tr\ e^{-t\triangle_\varphi^{i,j,l}}(-1)^{j+1}$.

**Proposition 6.3.1** *We have*

$$\Theta(t) = \sum_{l=0}^{dim\ N} \sum_{i\in I_l} (-1)^l \Theta_{i,l}(t).$$

This follows from the trace formula since the orbital integrals of both sides coincide.

Note that there is some cancellation in the preceding proposition. To describe this, consider the $M^0$-decomposition

$$\wedge^l\mathfrak{n} = \bigoplus_{i\in I_l} a_{i,l} V_{\lambda_{i,l}}$$

Now there is a differential of Lie algebra homology $\partial : \wedge^l \mathfrak{n} \to \wedge^{l-1} \mathfrak{n}$, making $\wedge^* \mathfrak{n}$ a complex. Let $H_*(\mathfrak{n})$ denote its homology, then, since $\partial$ is $M$-equivariant, we have as virtual $M$-modules

$$\wedge^* \mathfrak{n} \cong H_*(\mathfrak{n}).$$

Now let $J_l \subset I_l$ be the subset of isotypes which also occur in $H_l(\mathfrak{n})$, so

$$H_l(\mathfrak{n}) = \bigoplus_{i \in J_l} b_{i,l} V_{\lambda_{i,l}}$$

for some coefficients $0 < b_{i,l} \leq a_{i,l}$, but we actually have

**Lemma. 6.3.2** *Suppose $i \in J_l$, then $a_{i,l} = b_{i,l}$.*

**Proof:** This follows immediately from Corollary 5.7 in [Kost]. ∎

**Proposition 6.3.2**

$$\Theta(t) = \sum_{l=0}^{dim} \sum_{i \in J_l}^{N} (-1)^l \Theta_{i,l}(t).$$

∎

**Lemma. 6.3.3** *Let $i \in J_l$, then*

$$B(\lambda_{i,l}) - B(\rho) = -B(\rho_0 + \chi_{i,l}),$$

*where $\rho_0 = \rho | \mathfrak{a} \in \mathfrak{a}^*$ and $\chi_{i,l} \in \mathfrak{a}^*$ is the character of $\mathfrak{a}$ on $V_{\lambda_{i,l}} \subset \wedge^l \mathfrak{n}$.*

Analogous to the above we get

**Lemma. 6.3.4**

$$\Theta_{i,l}(t) = \sum_{[1] \neq [\gamma] \in \mathcal{E}_f(\gamma)} tr\varphi(\gamma) \frac{vol(\Gamma_\gamma \backslash G_\gamma)[G_h : G_h^0]^{-1}}{\xi_\rho(\gamma) \prod_{\alpha \in \Phi^+ - \Phi_{\gamma,c}^+} (1 - \xi_{-\alpha}(\gamma))} c_\gamma^{-1}$$

$$\frac{e^{-l_\gamma^2/4t - tB(\rho)}}{\sqrt{4\pi t}^r} a_{i,l} e^{tB(\lambda_{i,l})} \sum_{w \in W_{\mathbb{C}}(M)} \epsilon(w) \tilde{\omega}_h(w(\lambda_{i,l})) \xi_{w(\lambda_{i,l})}(B)$$

$$+ dim \ F \ vol(X) a_{i,l} \sum_{\xi \in \mathcal{E}_2(M^0)} dim(W_\xi \otimes \wedge^* \mathfrak{p}_M \otimes V_{\lambda_{i,l}})^{K \cap M^0} F_\xi(t).$$

We want to compute the Mellin transform of $\Theta_{i,l}$. So we define

$$M_{\lambda}^{i,l} = \int_0^\infty t^{z-1} e^{-\lambda t} \Theta_{i,l}(t) dt.$$

To compute the geodesic terms we need to consider

$$g_b^a(z) = \int_0^\infty t^{z-1} e^{-(\frac{a}{t} + bt)} dt, \; P \qquad a, b > 0.$$

We define $g_a(z) = g_a^a(z)$ and find

$$g_b^a(z) = \sqrt{a/b}^z \, g_{\sqrt{ab}}(z).$$

We further see

$$\begin{aligned}
g_a(-z) &= g_a(z), \\
z g_a(z) &= a(g_a(z+1) - g_a(z-1)) \\
\tfrac{\partial}{\partial a} g_a(z) &= -g_a(z-1) - g_a(z+1), \\
g_a(-\tfrac{1}{2}) &= \pi^{\frac{1}{2}} e^{-2a} a^{-\frac{1}{2}}.
\end{aligned}$$

For $n \geq 0$ define a polynomial $P_n(x)$ by

$$P_0(x) = 1, \quad P_1(x) = 1,$$

$$P_{n+1}(x) = P_{n-1}(x) + (2n-1)x P_n(x).$$

Note that $deg\ P_n = n-1$, the constant term is 1 and the highest coefficient equals $(2n-3)!! = (2n-3)(2n-5)\dots 1$. By the above we get

$$g_a(-\frac{2n+1}{2}) = \sqrt{\pi} \; \frac{e^{-2a}}{\sqrt{a}} \; P_{n+1}(\frac{1}{2a}).$$

## 6.4   The Case of Odd Fundamental Rank

In this section we assume the fundamental rank $r$ of $G$ to be odd. We write

$$\sum_{\xi \in \mathcal{E}_2(M^0)} F_\xi(t) \dim(W_\xi \otimes \wedge^* \mathfrak{p}_M \otimes V_{\lambda_{i,l}})^{K \cap M^0} = \sum_{n=0}^{\dim(N)} c_n^{i,l} t^{-n-\frac{r}{2}}.$$

We further write $\mu_\gamma$ for the quotient $l_\gamma / \lambda_\gamma$. If we denote by $\triangle_\varphi^{i,l}$ the operator $\triangle_\varphi^{*,i,l}$ on $\Omega^{*,i,l}$ then by 3.11 we get

$$\begin{aligned}
&log\ det(\triangle_\varphi^{i,l} + s^2 - B(\rho) + B(\lambda_{i,l})) \\
&= \sum_{1 \neq [\gamma] \in \mathcal{E}_f(\gamma)} tr\ \varphi(\gamma)\ \chi_r(X_\gamma)\ D_{i,l}(\gamma) \\
&\quad (\frac{2s}{l_\gamma})^{\frac{1}{2}(r-1)} P_{[\frac{1}{2}r]}(\frac{1}{sl_\gamma}) \frac{\sqrt{\pi}}{s} \frac{e^{-sl_\gamma}}{\mu_\gamma} \\
&\quad + \dim\ F\ vol(X_\Gamma) \sum_{n=0}^{\dim(N)} c_n^{i,l} \Gamma(-n - \frac{r}{2}) s^{2n+r}
\end{aligned}$$

with

$$D_{i,l} = \frac{|W_{\mathbb{C}}(G_\gamma)|^{-1} \prod_{\alpha \in \Phi_\gamma^+}(\rho_\gamma, \alpha)^{-1}}{\xi_{\rho_m} \prod_{\Phi^+ - \Phi_\gamma^+}(1 - \xi_{-\alpha}(\gamma))} \sum_{w \in W_{\mathbb{C}}(M)} \epsilon(w)\tilde{\omega}_\gamma(w(\lambda_{i,l}))\xi_{w(\lambda_{i,l})}(B).$$

Now let

$$Z_\varphi^{i,l} = exp(\sum_{1 \neq [\gamma] \in \mathcal{E}_f(\gamma)} tr \; \varphi(\gamma) \; \chi_r(X_\gamma) \; D_{i,l}(\gamma)(\frac{2s}{l_\gamma})^{\frac{1}{2}(r-1)}P_{[\frac{1}{2}r]}(\frac{1}{sl_\gamma})\frac{\sqrt{\pi}}{s}\frac{e^{-sl_\gamma}}{\mu_\gamma}).$$

Let $\sqrt{\cdot}$ denote the standard branch of the square root and let

$$Z_\varphi^\pm = \prod_{i,l} Z_\varphi^{i,l}(s \pm \sqrt{B(\rho) - B(\lambda_{i,l})})^{a_{i,l}(-1)^l}.$$

Let $R_r$ denote the rational function

$$R_r(x,y) = (\frac{2x}{y})^{\frac{1}{2}(r-1)}P_{[\frac{1}{2}r]}(\frac{1}{xy})x^{-1}.$$

**Proposition 6.4.1** *We get*

$$Z_\varphi^\pm(s) = exp(\sum_{1 \neq [\gamma] \in \mathcal{E}_f(\gamma)} tr \; \varphi(\gamma) \; \chi_r(X_\gamma) \; L^\pm(s,\gamma)\frac{e^{-l_\gamma s}}{\mu_\gamma}),$$

*where*

$$L^\pm(s,\gamma) = \frac{tr(R_r(s \pm |\chi + \rho|, l_\gamma)a_\gamma^{\pm|\chi+\rho|}\gamma_I \mid H^*(\mathfrak{n}))}{tr(\gamma \mid H^*(\mathfrak{n}))},$$

*here $\chi$ stands for the weight of the $A$-action and $|.|$ means absolute value on $\mathfrak{a}^*$. As usual $H^*(\mathfrak{n})$ means Lie algebra cohomology.*

For the proof let $T_{s,\gamma}^\pm$ denote the operator on $\wedge^*\mathfrak{n}$ acting on $V(\lambda_{i,l})$ by the scalar

$$R_r(s \pm \sqrt{B(\rho) - B(\lambda_{i,l})}, l_\gamma)e^{\mp l_\gamma \sqrt{B(\rho) - B(\lambda_{i,l})}}.$$

Analogous to [MS2] we have

$$\sum a_{i,l}(-1)^l R_r(s \pm \sqrt{B(\rho) - B(\lambda_{i,l})}, l_\gamma)D_{i,l}(\gamma) = \frac{tr(T_{s,\gamma}^\pm \gamma_I \mid \wedge^*\mathfrak{n}^*)}{det(1 - \gamma \mid \mathfrak{n})}$$

Since T interchanges with the differential on $\wedge\mathfrak{n}^*$ we have

$$tr(T_{s,\gamma}^\pm \gamma_I \mid \wedge^*\mathfrak{n}^*) = tr(T_{s,\gamma}^\pm \gamma_I \mid H^*(\mathfrak{n})).$$

By corollary 5.7 in [Kost] we get the claim.

Putting things together we get

**Theorem 6.4.1** *Let the fundamental rank $r$ of $G$ be odd. The zeta functions*

$$Z_\varphi^\pm(s) = exp\Big(\sum_{1 \neq [\gamma] \in \mathcal{E}_f(\gamma)} tr \ \varphi(\gamma) \ \chi_r(X_\gamma) \ L^\pm(s,\gamma) \frac{e^{-l_\gamma s}}{\mu_\gamma}\Big)$$

*admit analytic continuation to meromorphic functions on the complex plane. They satisfy a functional equation*

$$Z_\varphi^-(-s) = Z_\varphi^+(s) e^{2P(s)}$$

*where $P$ denotes a polynomial of degree $2 \dim(N)$. Around zero we find*

$$Z_\varphi^\pm(s) = \tau_r(E_\varphi) s^{\chi_r(X_\Gamma, \varphi)} + \text{higher order terms.}$$

*Furthermore we have $\tau_1(E_\varphi) = \ldots = \tau_{r-1}(E_\varphi) = 1$.*

## 6.5   The Generalized Selberg Zeta Function

Let $r$ denote the fundamental rank of $G$ which we assume to be $> 0$. For $\gamma \in \Gamma$ let $\lambda_\gamma$ denote the volume of a maximal compact flat containing the geodesic $[\gamma]$ and let $\mu_\gamma = \frac{l_\gamma}{\lambda_\gamma}$ denote the **generic multiplicity** .

**Theorem 6.5.1** *Let the lattice $\Gamma$ be good and nice and consider the zeta function*

$$Z_\varphi(s) = \prod_{\substack{[\gamma] \text{ prime,} \\ \text{fundamental}}} \prod_{N \geq 0} \det(1 - \xi_{-\rho_0}(\gamma) e^{-sl_\gamma} | V_N)^{\frac{\chi_r(X_\gamma)}{\mu_\gamma}},$$

*where $V_N = V_\varphi \otimes H_*(\mathfrak{n}) \otimes S^N(\mathfrak{n})$ and $\gamma$ acts as $\varphi(\gamma) \otimes \mathrm{Ad}(\gamma_I) \otimes \mathrm{Ad}^N(\gamma^{-1})$.*

*a) If $r = 1$ then $Z_\varphi(s)$ extends to a meromorphic function. If the complex number $s_0$ is a singularity (pole or zero) then $s_0^2 + B(\rho_0 + \chi_{i,l})$ is an eigenvalue of $\triangle_{i,l}$ for some $l$ and some $i \in J_l$. The (vanishing-) order of $Z_\varphi$ at $s_0$ is*

$$\sum_{l=0}^{\dim(N)} (-1)^l \sum_{i \in J_l} \dim(\mathrm{Eig}(\triangle_{i,l}, s_0^2 + B(\rho_0 + \chi_{i,l})).$$

*b) If $r > 1$ is odd, the logarithm $\log(Z_\varphi(s))$ extends to a meromorphic function. The poles of $\log(Z_\varphi(s))$ have order $\frac{r-1}{2}$ and if $s_0$ is a pole then $s_0^2 + B(\rho_0 + \chi_{i,l})$ is an eigenvalue of $\triangle_{i,l}$ for some $l$ and some $i \in J_l$. The principal of the Laurent series around $s_0$ then is*

$$(\frac{1}{s^2 - s_0^2})^{\frac{r-1}{2}} \sqrt{4\pi}^{r-1} \sum_{l=0}^{\dim(N)} (-1)^l \sum_{i \in J_l} \dim(\mathrm{Eig}\triangle_{i,l}, s_0^2 + B(\rho_0 + \chi_{i,l})).$$

*c) If $r > 0$ is even, the square of the logarithm, $(\log(Z_\varphi(s)))^2$ extends to a meromorphic function. The poles of $(\log(Z_\varphi(s)))^2$ have order $r-1$ and if $s_0$ is a pole then $s_0^2 + B(\rho_0 + \chi_{i,l})$ is an eigenvalue of $\triangle_{i,l}$ for some $l$ and some $i \in J_l$. The principal part of the Laurent series around $s_0$ then is*

$$(\frac{1}{s^2 - s_0^2})^{r-1}(4\pi)^{r-1}(\sum_{l=0}^{\dim(N)}(-1)^l\sum_{i\in J_l}\dim(\text{Eig}\triangle_{i,l}, s_0^2 + B(\rho_0 + \chi_{i,l})))^2.$$

**Proof:** We take the formula of Lemma 6.3.4 and multiply both sides with $\sqrt{4\pi t}^{r-1}$. Taking the Mellin transform now gives on the spectral side

$$\frac{1}{\Gamma(s)}\int_0^\infty t^{s-1}\Theta_{i,l}(t)\sqrt{4\pi t}^{r-1}e^{t\lambda}dt = \sqrt{4\pi}^{r-1}\frac{\Gamma(s+\frac{r-1}{2})}{\Gamma(s)}\sum_{n=1}^\infty(\lambda_n + \lambda)^{-s-\frac{r-1}{2}}$$

where $(\lambda_n)$ runs through the eigenvalues of $\triangle_{i,l}$. From this the theorem follows easily. ∎

# Chapter 7

# Dynamical Systems

In this chapter we are going to present a different approach to the zeta function which gives rise to some interesting speculations concerning possible generalizations.

Let $P = MAN$ be a minimal parabolic of $G$ with $M \subset K$ and $\mathrm{Lie}(A) \perp \mathrm{Lie}(K)$. Then $M$ is the centralizer of $A$ in $K$, in fact, the group $AM$ will be the total centralizer of $A$. Note that the dimension of $A$ is just the $\mathbb{R}$-rank of $G$ or the rank of the symmentric space $X = G/K$. In the rank one case the sphere bundle writes

$$SX = G/M \to G/K.$$

The geodesic flow then is

$$\Phi_t(gM) = g \exp(tH)M,$$

where $H \in \mathrm{Lie}(A)$ is an element of norm 1. Now let $\Gamma$ be a torsionfree uniform lattice in $G$ and let $X_\Gamma = \Gamma \backslash X$. Among other things we considered the **Ruelle Zeta Function** of the geodesic flow on $SX_\Gamma = \Gamma \backslash G/M$:

$$R(s) = \prod_c (1 - e^{-sl_\gamma}),$$

and we have seen that $R$ extends to a meromorphic function which satisfies a determinant formula, involving Laplace operators on $X_\Gamma$.

On the other hand we know that the geodesic flow belongs to the very general class of Anosov flows and that the Ruelle zeta function of an Anosov flow already extends to a meromorphic function under much more general conditions then ours [Rue]. However, there's no determinant formula in general and only very little is known about the analysis of a Ruelle zeta function.

To generalize to arbitrary Anosov flows we would have to get rid of the space $X_\Gamma$ and have an approach purely in terms of the sphere bundle $SX_\Gamma = \Gamma \backslash G/M$. Such an approach indeed exists, it was discovered by A. Juhl [Ju]. The reason, why it works is the following: In our approach we had to resort to $X_\Gamma$ because the irreducible representations of $G$ have finite $K$-type multiplicities, but not finite $M$-type multiplicities. Viewing the $G$-representation space $L^2(\Gamma \backslash G)$ as a sheaf over $\hat{G}$ this means that the restriction to any $K$-type $\tau$ (for example to p-forms on

$X_\Gamma$):

$$(L^2(\Gamma\backslash G)\otimes\tau)^K$$

gives a sheaf with finite dimensional stalks over $\hat{G}$.

On the other hand the group $M$ is the maximal compact subgroup of the Levi group $AM$ so the same would count for appropriate sheaves on $\hat{AM}$. There are natural functors joining $\hat{G}$ to $\hat{AM}$. Of special interest to us is one that preserves characters of representations.

## 7.1  Juhl's Approach

### 7.1.1  The Hecht-Schmid Character Formula

For any Lie algebra $\mathcal{L}$ and any $\mathcal{L}$-module $V$ we write $H_q(\mathcal{L},V)$ for the $q$-th Lie algebra homology of $\mathcal{L}$ with coefficients in $V$ and $H^q(\mathcal{L},V)$ for the $q$-the Lie algebra cohomology.

Let $G$ be a real reductive group with maximal compact $K$ and $V$ be a **Harish-Chandra** module of $G$, that is, $V$ is an admissible $(\mathfrak{g},K)$-module, where $\mathfrak{g}$ is the Lie algebra of $G$. Recall that $V$ has a character $\Theta_G(.|V)$ defined by Harish-Chandra. $\Theta_G(.|V)$ is a distribution on $G$ given by a locally integrable and conjugation invariant function also denoted by $\Theta_G(.|V)$, which s real analytic on the regular set of $G$.

Now let $P=MAN$ be a parabolic subgroup of $G$. Let $\Phi^+(\mathfrak{g},\mathfrak{a})$ denote the set of positive roots according to the choice of $P$ and let $A^+$ denote the positive Weyl chamber in $A$. Let $\mathfrak{n}$ be the Lie algebra of $N$ and let $\mathfrak{n}^-=\theta(\mathfrak{n})$ where $\theta$ is the Cartan involution. Set

$$(MA)^+ \stackrel{\text{def}}{=} \text{ interior in } MA \text{ of the set}$$

$$\{g\in MA|\det(1-ga|\mathfrak{n}^-)\geq 0 \text{ for all } a\in A^+\}.$$

Note that for $\dim(A)$ maximal we have $(MA)^+=MA^+$. Hecht and Schmid proved in [HeSchm]:

**Theorem 7.1.1** *For every Harish-Chandra module $V$ we have that*

- $H_p(\mathfrak{n}^-,V)$ *is a Harish-Chandra module for the reductive group $MA$.*

- *For all $x\in(MA)^+\cap G^{\text{reg}}$, the set of regular elements in $MA$ we have*

$$\Theta_G(x|V)=\frac{\sum_{p=0}^{\dim(\mathfrak{n}^-)}(-1)^p\Theta_{MA}(x|H_p(\mathfrak{n}^-,V))}{\det(1-x|\mathfrak{n}^-)}.\blacksquare$$

### 7.1.2  Nice Test Functions

Now assume $G$ semisimple with finite center and $K$ a maximal compact subgroup. Let $P=MAN$ be the Langlands decomposition of a **minimal** parabolic subgroup. Then $A$ is a maximal split torus and we may assume $M\subset K$ and $\text{Lie}(A)\perp\text{Lie}(K)$. Let $\mathfrak{a}_0$ denote the Lie algebra of $A$, $\mathfrak{a}_0^+$ a positive Weyl chamber according to the choice of $P$, $A^+ = \exp(\mathfrak{a}_0^+)$. Let $\mathfrak{n}_0^+$ denote the Lie algebra of $N$ and $\mathfrak{n}_0^-$ its opposite, i.e.: $N_0^- = \theta(\mathfrak{n}_0^+)$, where $\theta$ is the Cartan involution. Let

$$D_\pm(ma) \overset{\text{def}}{=} \det(1 - ma|\mathfrak{n}_0^\pm)$$

for $ma \in MA$. We have the integration formula

$$\int_G f(g)dg = \int_{MA\backslash G}\int_M\int_{A^+} f(g^{-1}mag)|D_+(ma)D_-(ma)|dadmdg$$

for $f \in C_c^\infty(U)$, where $U = \{gxg^{-1} | g \in G, x \in (MA) \cap G^{\text{reg}}\}$. All Haar measures in this formula are supposed to be normalized in the sense of Harish-Chandra ([HC-HAI],sec. 7).

To create our test functions let $\varphi \in C_c^\infty(A^+)$ and choose $\eta \in C_c^\infty(MA\backslash G)$ such that $\eta \geq 0$ and $\int_{MA\backslash G} \eta(g)dg = 1$. Define

$$\Phi(g^{-1}mag) = \eta(g)\varphi(a)|D_+(ma)|^{-1}.$$

Then $\Phi$ is well defined, smooth and compactly supported. Further the support of $\Phi$ lies completely in the set of conjugates of elements of $MA$ with regular $A$-part.

Let $\pi$ be an irreducible unitary representation of $G$. By the Hecht-Schmid character formula we get

$$\begin{aligned}
\text{tr}\pi(\Phi) &= \int_G \Theta_G(x|\pi)\Phi(x)dx \\
&= \int_{MA\backslash G}\int_M\int_{A^+} \eta(g)|D_-(ma)|\Theta_G(ma|\pi)dadmdg \\
&= \int_{A^+} \varphi(a)\int_M \sum_{p=0}^{\dim(N)}(-1)^p\Theta_{MA}(ma|H_p(\mathfrak{n}^-,V))dmda \\
&= \int_{A^+} \varphi(a)\text{tr}(a|H_*(\mathfrak{n}^-,V)^M)da.
\end{aligned}$$

Note that the space $H_p(\mathfrak{n}^-,V)^M$ is finite dimensional.

### 7.1.3 The Lefschetz Formula

For any $\gamma \in \Gamma$ such that $\gamma$ is conjugate to an element of $MA$ write $m_\gamma a_\gamma$ for such an element. Whenever this notation is used, one has to make clear that it is independent from the particular $m_\gamma a_\gamma \in MA$ chosen.

**Theorem 7.1.2** *(Juhl) Let $\varphi \in C_c^\infty(A^+)$ then we have the identity*

$$\sum_{\pi\in\hat{G}} N_\Gamma(\pi)\int_{A^+}\varphi(a)\text{tr}(a|H_*(\mathfrak{n}^-,V)^M)da = \sum_{[\gamma]\subset\Gamma}\lambda_\gamma\frac{\varphi(a_\gamma)}{|D_+(m_\gamma a_\gamma)|},$$

*where the sum runs over all $\gamma \in \Gamma$ that are in $G$ conjugate to an element of $MA$ and $\lambda_\gamma = \mathrm{vol}(\Gamma_\gamma \backslash MA)$ is the volume of a maximal flat containing the geodesic $[\gamma]$.*

**Proof:** One simply plugs the function $\Phi$ into the trace formula. It is clear by the above calculation that the spectral side gets the above shape. The expression for the geometric side follows directly from the definition of $\Phi$. ∎

This Theorem can easily be generalized allowing test functions $\varphi$ which are defined on $A^+$, vanish on the boundary to a sufficiently high order and are decaying at infinity fast enough.

**Corollary 7.1.1** *Let $\varphi \in C_c^\infty(A^+)$ then we have the identity*

$$\sum_{\pi \in \hat{G}} N_\Gamma(\pi) \int_{A^+} \varphi(a) \mathrm{tr}(a | H^*(\mathfrak{n}^-, V_\pi)^M) da = (-1)^{\dim N} \sum_{[\gamma] \subset \Gamma} \lambda_\gamma \frac{\varphi(a_\gamma)}{\det(1 - m_\gamma a_\gamma | \mathfrak{n}^-)},$$

*where the sum runs over all $\gamma \in \Gamma$ that are in $G$ conjugate to an element of $MA$.*

**Proof:** We use the natural isomorphism:

$$H_p(\mathfrak{n}^-, V) \cong H^{d-p}(\mathfrak{n}^-, V) \otimes \wedge^d \mathfrak{n}^-,$$

([HeSchm] (2.18)) where $d$ denotes the dimension on $\mathfrak{n}$. Since $M$ acts trivially on $\wedge^d \mathfrak{n}^-$ we get

$$\mathrm{tr}(a | H_*(\mathfrak{n}^-, V)^M) = (-1)^d \mathrm{tr}(a | H^*(\mathfrak{n}^-, V)^M \xi_{-2\rho}(a).$$

So one simply replaces $\varphi(a)$ by $\varphi(a) \xi_{2\rho}(a)$ to achieve the assertion of the corollary out of the theorem. ∎

By plugging in suitable test functions $\varphi$, Juhl arrives at expressions involving the Selberg zeta function. This works well for the rank one situation. For higher rank we may put in modified test functions as in 5.4

Then we get

**Theorem 7.1.3** *(Higher rank Lefschetz Theorem) Suppose $\Gamma$ is nice. Let $H = AB$ any $\Theta$-stable Cartan subgroup of $G$ with $r = \dim(A) > 0$. Choose a parabolic $P = MAN$. Then we have an equality of distributions on the positive Weyl chamber $A^+$:*

$$\sum_{\pi \in \hat{G}} N_\Gamma(\pi) \mathrm{tr}(a | H^*(\bar{\mathfrak{n}}, \pi) \otimes \wedge^* \mathfrak{p}_M)^{K_M}) = (-1)^{\dim(N)} \sum_{[\gamma] \in \mathcal{E}_H(\Gamma)} \lambda_\gamma \chi_r(X_\gamma) \frac{\delta_{a_\gamma}}{\det(1 - m_\gamma a_\gamma | \bar{\mathfrak{n}})}. \blacksquare$$

## 7.2 The Guillemin Conjecture

The spectral side of the corollary in the last section can be viewed as the *A*-trace on the complex

$$0 \to (L^2(\Gamma\backslash G) \otimes (\wedge^0 \mathfrak{n}^-)^*)^M \to (L^2(\Gamma\backslash G) \otimes (\wedge^1 \mathfrak{n}^-)^*)^M \to \dots$$

$$\dots \to (L^2(\Gamma\backslash G) \otimes (\wedge^{\text{top}} \mathfrak{n}^-)^*)^M \to 0,$$

with the Lie algebra differential. But this complex also has a different interpretation. Recall that the stable bundle on $SX_\Gamma = \Gamma\backslash G/M$ could be described as

$$T^s = \Gamma\backslash G \times_M \mathfrak{n}^-.$$

so the above complex becomes the **stable complex** :

$$o \to L^2(\wedge^0 T^s) \to L^2(\wedge^1 T^s) \to \dots \to L^2(\wedge^{\text{top}} T^s) \to 0,$$

and a short computation shows that the differential is just the exterior differential of $SX_\Gamma$ followed by a projection $\wedge^p T \to \wedge^p T^s$, where $T$ denotes the tangent bundle of $SX_\Gamma$. Note that the stable complex is not elliptic, a fact which is remedied in our situation by the cohomology having finite multiplicities as a *G*-representation. Now one might present the following

**Conjecture 7.2.1** *Let $\Phi_t$ denote a smooth Anosov flow on a compact manifold then the flow has a distributional trace for $t > 0$ on the cohomology of the stable complex such that for $t > 0$*

$$\text{tr}(\Phi_t | H^*(L^2(\wedge^\cdot T^s))) = (-1)^{\dim T^s} \sum_c \frac{l(c_0)\delta_{l(c)}}{\det(1 - P_c^s)}.$$

*Here the sum runs over all closed orbits $c$ of $\Phi$ and $l(c)$ denotes the length of the orbit $c$. To an orbit $c$ let $c_0$ denote the underlying primitive orbit and $P_c^s$ the stable part of the Poincaré map around $c$ and $\delta_x$ denotes the delta distribution at the point $x$.*

Compare this to ([Gu], p.504).

Perhaps one should weaken this conjecture in insisting that there should exist, as in our case, a smooth perfect pairing

$$B_x : T_x^s \times T_x^u \to \mathbb{C}$$

which is invariant under the flow, i.e. $B_{\Phi_t x}(D\Phi_t v, D\Phi_t w) = B_x(v, w)$.

Once having established this conjecture in a particular setting one can proceed as in [Ju] to finally get an analysis of the Ruelle zeta function in terms of the spectrum of the flow on the stable cohomology. This would possibly clarify the relation of the Ruelle zeta function to topological invariants as torsion and Euler characteristics. Further it would show whether the Ruelle zeta function in general satisfies an analogue of the Riemann hypothesis.

# Chapter 8

# The Functional Equation of Theta Series

Geometric zeta functions are defined by geometrical data (closed geodesics) and turn out to have singularities (zeroes and poles) which are related to spectral data (Laplace eigenvalues). The theta series to be discussed in this section are dual to the zeta functions in that they are defined by the very same spectral data and have their singularities at the lengths of the closed orbits.

## 8.1 The Rank One Case

For the contents of this subsection compare [CaVo] and [BuOl]. We consider the situation of chapter 4. So $G$ is the group of orientation preserving isometries of the even-dimensional rank one symmetric space $X$, which we identify with $G/K$ where $K$ is a maximal compact subgroup of $G$. We fix a torsionfree uniform lattice $\Gamma$ in $G$ and consider the compact Shimura manifold $X_\Gamma = \Gamma \backslash X$.

In chapter 4 we considered the zeta function

$$Z_{\sigma,\varphi}(s) = \prod_{[\gamma] \text{ prime}} \prod_{N \geq 0} \det(1 - e^{-sl_\gamma}\gamma \mid V_N),$$

where $V_N = V_\varphi \otimes V_\sigma \otimes S^N(\mathfrak{n})$ and $\gamma$ acts on $V_N$ as $\gamma \mapsto \varphi(\gamma) \otimes \sigma(m_\gamma) \otimes Ad^N((m_\gamma a_\gamma)^{-1})$. We showed (4.3.1) that with $c_\sigma = B(\lambda_\sigma) - B(\rho)$ we have

$$Z(\rho_0 + s) = \det(\triangle_{\sigma,\varphi} + s^2 + c_\sigma)(\det(D_\sigma + s)\exp(P(s^2)))^{2(-1)^m \dim(\varphi)\chi(X_\Gamma)/\chi(X^d)}.$$

Now we are going to consider the **theta series** :

$$\Theta_{\sigma,\varphi}(\tau) = \text{tr}(e^{-\tau\sqrt{\triangle_{\sigma,\varphi} + c_\sigma}}),$$

where the root is positive on $\mathbb{R}_+$ and negative imaginary on $\mathbb{R}_-$. Let

$$
\begin{aligned}
F(s) &= \det(\triangle_{\sigma,\varphi} + s^2 + c_\sigma) \\
&= \prod_{j \geq 0}{}' (a_j^2 + s^2)
\end{aligned}
$$

where $a_j \geq 0$ or $ia_j \geq 0$ and the product means regularized product. So, if $(\lambda_j)_{j \in \mathbb{N}}$ is the sequence of eigenvalues of the operator $\triangle_{\sigma,\varphi}$, then $a_j^2 = \lambda_j + c_\sigma$. Assume that the $a_j$ are numerated in such a way that the sequence $a_j^2$ is ascending. There are only finitely many $a_j$ on the imaginary axis and there is a $p > 0$ such that $\sum_j' |a_j|^{-p} < \infty$. Let $a > |a_0|$ be arbitrary. Then the theta series is

$$
\Theta_{\sigma,\varphi}(\tau) = \sum_{j=0}^{\infty} e^{-\tau a_j} \quad \mathrm{Re}\ (\tau) > 0.
$$

We are further going to consider the **dual space theta function** :

$$
\Theta_{\sigma,d}(\tau) = \mathrm{tr}(e^{-\tau D_\sigma}), \quad \mathrm{Re}\ (\tau) > 0.
$$

Since we fix the $M$-type $\sigma$ and the $\Gamma$-representation $\varphi$ we are going to leave these out in our notation so we will write $\Theta = \Theta_{\sigma,\varphi}$ and $\Theta_d = \Theta_{\sigma,d}$. Our first result is

**Proposition 8.1.1** *The function $\Theta_d = \Theta_{\sigma,d}$ extends to a meromorphic function on the entire plane. It is even, i.e.:*

$$
\Theta_d(-\tau) = \Theta_d(\tau),
$$

*its poles are all of order $\dim(X)$ and are located at $\tau = \pi i k,\ k \in \mathbb{Z}$.*

**Proof:** We first consider the case when $\sigma$ is even in the sense of 4.3.2:

$$
\begin{aligned}
\Theta_d(\tau) &= \sum_{n=1}^{\infty} Q_\sigma(2n) e^{-2\tau n} \\
&= Q_\sigma(-\frac{\partial}{\partial \tau}) \sum_{n=1}^{\infty} e^{-2\tau n} \\
&= Q_\sigma(-\frac{\partial}{\partial \tau})(\frac{1}{e^{2\tau} - 1})
\end{aligned}
$$

This gives the meromorphic continuation and the location and order of the poles (recall that the polynomial $Q_\sigma$ has order $\dim(X) - 1$. To see that $\Theta_d$ is even, recall that

$$
Q_\sigma(x) = x\tilde{Q}_\sigma(x^2)
$$

for a polynomial $\tilde{Q}_\sigma$. It follows

$$
\Theta_d(\tau) = 2\tilde{Q}_\sigma((\frac{\partial}{\partial \tau})^2)(\frac{1}{(e^t - e^{-t})^2}),
$$

which gives the claim. For $\sigma$ odd we similarly get

$$\Theta_d(\tau) = Q_\sigma(-\frac{\partial}{\partial \tau})\frac{1}{e^t - e^{-t}},$$

which gives the claim in that case. ∎

**Theorem 8.1.1** *The theta series* $\Theta = \Theta_{\sigma, \varphi}$ *admits a continuation to a meromorphic function on* $\mathbb{C}$. *It satisfies a functional equation:*

$$\Theta(\tau) + \Theta(-\tau) = 2(-1)^{m+1}\dim(\varphi)\frac{\chi(X_\Gamma)}{\chi(X^d)}\Theta_d(i\tau).$$

*The function* $\Theta(\tau)$ *has simple poles at* $\tau = \pm i l_\gamma$, $\gamma \neq 1$ *and the residue at* $\tau = \pm i l_\gamma$ *is*

$$-\frac{1}{2\pi}\frac{l_\gamma}{\mu_\gamma}\frac{\operatorname{tr}\varphi(\gamma)\operatorname{tr}\sigma(m_\gamma)}{\det(1 - \gamma^{-1}|\mathfrak{n})}.$$

*The function* $\Theta(\tau)$ *further has poles of order* $\dim(X)$ *at the nonpositive poles of* $\Theta_d(i\tau)$. *This concludes the list of poles of* $\Theta$.

**Proof:** The argument principle shows that

$$\Theta(\tau) = n_0 + \frac{1}{2\pi i}\int_C e^{i\tau z}(\frac{F'}{F}(z) - \frac{2n_0}{z})dz,$$

where $n_0 = \operatorname{ord}_{z=0}F(z)$ and $C$ is the contour shown at the following picture:

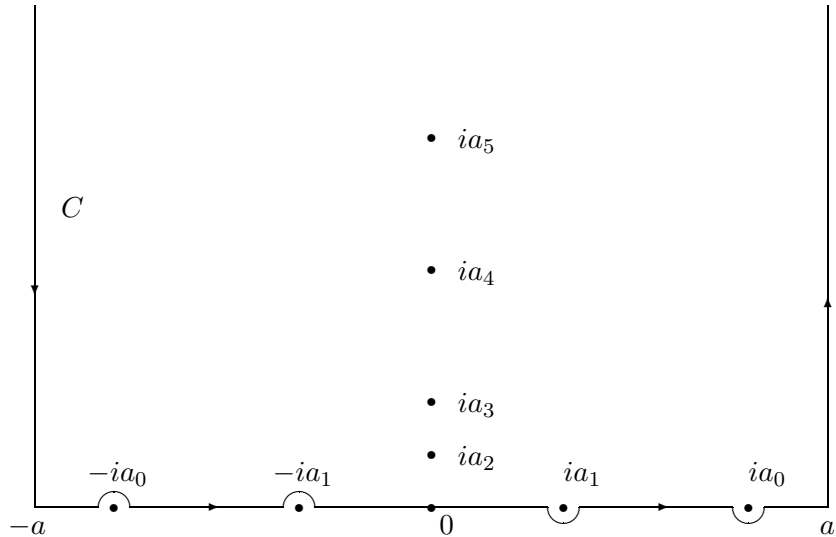

At first the contour goes from $-a + i\infty$ to $-a$, then from $-a$ to $a$ circumventing the $\pm i a_j$ which are nonzero as shown and then it goes from $a$ to $a + i\infty$. If one of the $a_j$ happens to be zero, which is equivalent to the fact that $n_0 \neq 0$, the contour goes right trough it.

Let $C^+ = C \cap \{\text{Re } (z) > 0\}$ and $C^- = C \cap \{\text{Re } (z) < 0\}$. Since $F$ is even it follows that the logarithmic derivative $\frac{F'}{F}$ is odd and hence

$$\Theta(\tau) = n_0 + \frac{1}{2\pi i} \int_{C^+} e^{i\tau z} (\frac{F'}{F}(z) - \frac{2n_0}{z}) dz.$$

Now by the determinant formula we get

$$\frac{F'}{F}(z) - \frac{2n_0}{z} = \frac{Z'}{Z}(\rho_0 + z) - \frac{2n_0}{z} + 2\alpha z P'(z^2) + \frac{d}{dz}\text{logdet}(D_\sigma + z),$$

where $\alpha = 2(-1)^{m+1} \dim(\varphi) \frac{\chi(X_\Gamma)}{\chi(X^d)}$.

Therefore $\Theta(\tau)$ splits as follows:

$$\Theta(\tau) = n_0 + \Theta_Z(\tau) + \Theta_{pol}(\tau) + \Theta_{dual}(\tau).$$

Let's start with the simplest part

*1. The polynomial contribution*
This is

$$
\begin{aligned}
\Theta_{pol} &= \frac{\alpha}{\pi i} \int_{C^+} e^{i\tau z} z P'(z^2) dz \\
&= \frac{\alpha}{\pi} \frac{\partial}{\partial \tau} P'(-(\frac{\partial}{\partial \tau})^2)(\int_{C^+} e^{i\tau z} dz - \int_{-C^-} e^{-i\tau z} dz) \\
&= \frac{\alpha}{\pi} \frac{\partial}{\partial \tau} P'(-(\frac{\partial}{\partial \tau})^2)(\frac{i}{\tau} - \frac{i}{\tau}) \\
&= 0.
\end{aligned}
$$

So the polynomial contribution vanishes.

*2. The zeta contribution*
Here we consider

$$
\begin{aligned}
\Theta_Z(\tau) &= \frac{1}{2\pi i} \int_{C^+} e^{i\tau z} (\frac{Z'}{Z}(\rho_0 + z) - \frac{2n_0}{z}) dz + \frac{1}{2\pi i} \int_{-C^-} e^{-i\tau z} (\frac{Z'}{Z}(\rho_0 + z) - \frac{2n_0}{z}) dz \\
&= \frac{1}{\pi} \int_0^a \sin(\tau z) (\frac{Z'}{Z}(\rho_0 + z) - \frac{2n_0}{z}) dz \\
&\quad + \frac{1}{2\pi i} \int_a^{a+i\infty} e^{-i\tau z} (\frac{Z'}{Z}(\rho_0 + z) - \frac{2n_0}{z}) dz \\
&\quad + \frac{1}{2\pi i} \int_{a-i\infty}^a e^{-i\tau z} (\frac{Z'}{Z}(\rho_0 + z) - \frac{2n_0}{z}) dz.
\end{aligned}
$$

The first integral should be read as a path integral over a path circumventing the singularities in the lower half plane. The first summand thus defines an odd entire function in $\tau$. Therefore it does not contribute to the functional equation.

We may assume $a$ to be large enough so that the Euler product for the zeta function $Z$ converges at Re $(z) = a$. Then we get

$$
\begin{aligned}
\log Z(\rho_0 + z) &= -\sum_{[\gamma]\,\mathrm{prime}} \sum_{N \geq 0} \sum_{k=1}^{\infty} \frac{1}{k} e^{-zkl_\gamma} \mathrm{tr}\varphi(\gamma^k)\mathrm{tr}\sigma(m_{\gamma^k})\mathrm{tr}\mathrm{Ad}^N((m_{\gamma^k}a_{\gamma^k})^{-1}|\mathfrak{n}) \\
&= -\sum_{[\gamma]\neq 1} \sum_{N \geq 0} \frac{e^{-zl_\gamma}}{\mu_\gamma} \mathrm{tr}\varphi(\gamma)\mathrm{tr}\sigma(m_\gamma)\mathrm{tr}\mathrm{Ad}^N((m_\gamma a_\gamma)^{-1}|\mathfrak{n}) \\
&= -\sum_{[\gamma]\neq 1} \frac{e^{-zl_\gamma}}{\mu_\gamma} \frac{\mathrm{tr}\varphi(\gamma)\mathrm{tr}\sigma(m_\gamma)}{\det(1-\gamma^{-1}|\mathfrak{n})},
\end{aligned}
$$

and so

$$
\frac{d}{ds}\log Z(\rho_0 + z) = \sum_{[\gamma]\neq 1} b_\gamma e^{-zl_\gamma},
$$

where

$$
b_g a = \frac{l_\gamma}{\mu_\gamma} \frac{\mathrm{tr}\varphi(\gamma)\mathrm{tr}\sigma(m_\gamma)}{\det(1-\gamma^{-1}|\mathfrak{n})}.
$$

From this we get

$$
\frac{1}{2\pi i}\int_a^{a+i\infty} e^{-i\tau z}(\frac{Z'}{Z}(\rho_0 + z) - \frac{2n_0}{z})dz + \frac{1}{2\pi i}\int_{a-i\infty}^a e^{-i\tau z}(\frac{Z'}{Z}(\rho_0 + z) - \frac{2n_0}{z})dz
$$

$$
= \frac{1}{2\pi i}\sum_{[\gamma]\neq 1} b_\gamma(\frac{e^{(-i\tau - l_\gamma)a}}{-i\tau - l_\gamma} - \frac{e^{(i\tau - l_\gamma)a}}{i\tau - l_\gamma}).
$$

which is an odd meromorphic function.

To finish the zeta contribution it remains to consider

$$
f(\tau) \stackrel{\mathrm{def}}{=} -\frac{1}{2\pi i}\int_{a-i\infty}^a e^{-i\tau z}\frac{2n_0}{z}dz - \frac{1}{2\pi i}\int_a^{a+i\infty} e^{i\tau z}\frac{2n_0}{z}dz.
$$

We see that $f$ is a primitive of the function $\frac{2n_0}{\pi}\frac{\sin(a\tau)}{\tau}$ and that $\lim_{\tau\to+\infty} f(\tau) = 0$. From this it follows that $f(0) = \frac{2n_0}{\pi}\int_0^\infty \frac{\sin(x)}{x}dx = -n_0$. Since $f$ is the primitive of an even entire function we conclude that $f$ is entire and

$$
f(\tau) + f(-\tau) = -2n_0,
$$

which makes the $2n_0$ drop out of the functional equation.

### 3. The dual space contribution

We consider

$$\Theta_{dual}(\tau) = \frac{\alpha}{2\pi i}\int_{C^+} e^{i\tau z} d\mathrm{logdet}(D_\sigma + z) + \frac{\alpha}{2\pi i}\int_{-C^-} e^{-i\tau z} d\mathrm{logdet}(D_\sigma + z).$$

We replace $-C^-$ by $C^-$ again and move the contour on to the imaginary axis to get

$$\begin{aligned}
\Theta_{dual}(\tau) &= \frac{\alpha}{2\pi i}\int_0^{i\infty} e^{i\tau z}(d\mathrm{logdet}(D_\sigma + z) - d\mathrm{logdet}(D_\sigma - z))\\
&= \frac{\alpha}{2\pi i}\int_0^{infty} e^{-\tau t}(\frac{d\mathrm{logdet}(D_\sigma + z)}{dz}|_{z=it} - \frac{d\mathrm{logdet}(D_\sigma - z))}{dz}|_{z=-it})dt
\end{aligned}$$

Now we have

$$\begin{aligned}
\frac{d}{dz}\mathrm{logdet}(D_\sigma + z) &= -\frac{d}{dz}\zeta'_{D_\sigma + z}(0)\\
&= -\frac{d}{dz}\frac{d}{ds}|_{s=0}\zeta_{D_\sigma + z}(s)\\
&= -\frac{d}{ds}|_{s=0}\frac{d}{dz}\zeta_{D_\sigma + z}(s)\\
&= \frac{d}{ds}|_{s=0}s\zeta_{D_\sigma + z}(s+1)\\
&= \zeta_{D_\sigma + z}(1),
\end{aligned}$$

since $\zeta_{D_\sigma + z}(s)$ is holomorphic at $s = 1$. Therefore we conclude

$$\Theta_{dual}(\tau) = \frac{\alpha}{2\pi i}\int_0^\infty e^{-\tau t}(\zeta_{D_\sigma + it}(1) - \zeta_{D_\sigma - it}(1))dt.$$

Now assume $\sigma$ odd and recall that

$$\begin{aligned}
\zeta_{D_\sigma + z}(1) &= \int_0^\infty t^{s-1}e^{-tz}\Theta_d(t)dt|_{s=1}\\
&= \int_0^\infty t^{s-1}e^{-tz}Q_\sigma(-\frac{\partial}{\partial t})\frac{1}{e^t - e^{-t}}dt|_{s=1}\\
&= \int_0^\infty t^{s-1}\frac{1}{e^t - e^{-t}}Q_\sigma(\frac{\partial}{\partial t})e^{-tz}dt|_{s=1}\\
&= Q_\sigma(-z)\int_0^\infty \frac{t^{s-1}e^{-tz}}{e^t - e^{-t}}dt|_{s=1}\\
&= Q_\sigma(-z)\sum_{n=0}^\infty \int_0^\infty t^{s-1}e^{-t(z+1)}e^{-2nt}dt|_{s=1}\\
&= Q_\sigma(-z)\sum_{n=0}^\infty \frac{1}{(z+2n+1)^s}|_{s=1}.
\end{aligned}$$

And thus

$$\zeta_{D_\sigma+it}(1) - \zeta_{D_\sigma-it}(1) = Q_\sigma(-it)(\sum_{n=0}^{\infty} \frac{1}{(2n+1+it)^s} - \sum_{n=0}^{\infty} \frac{1}{(2n+1-it)^s})|_{s=1}$$

$$= Q_\sigma(-it) \sum_{n=0}^{\infty} \frac{2it}{(2n+1)^2+t^2}$$

$$= \frac{\pi}{2} Q_\sigma(-it) i \tanh(\frac{\pi t}{2}).$$

So for $\sigma$ odd we get

$$\Theta_{dual} = \frac{\alpha}{4} \int_0^{\infty} e^{-\tau t} Q_\sigma(-it) \tanh(\frac{\pi t}{2}) dt$$

$$= \frac{\alpha}{4} Q_\sigma(i\frac{\partial}{\partial \tau}) \int_o^{\infty} e^{-\tau t} \tanh(\frac{\pi t}{2}) dt.$$

Since

$$\tanh(\frac{\pi t}{2}) = 1 + 2 \sum_{n=1}^{\infty} (-1)^n e^{-\pi n t},$$

we get

$$\Theta_{dual}(\tau) = \frac{\alpha}{4i} \tilde{Q}_\sigma(-(\frac{\partial}{\partial \tau})^2)(\frac{1}{\tau^2} + 2\sum_{n=1}^{\infty} (-1)^n \frac{1}{(\tau+\pi n)^2}),$$

which completes the meromorphic continuation. As for the functional equation one gets

$$\Theta_{dual}(\tau) + \Theta_{dual}(-\tau) = \frac{\alpha}{2i} Q_\sigma(i\frac{\partial}{\partial \tau}) \sum_{n \in \mathbb{Z}} (-1)^n \frac{1}{\tau+\pi n}$$

$$= \frac{\alpha}{2i} Q_\sigma(i\frac{\partial}{\partial \tau}) \frac{1}{e^t + e^{-t}},$$

which gives the claim.

For even $\sigma$ one gets analogously

$$\Theta_{dual}(\tau) = \frac{\alpha}{4} \int_0^{\infty} e^{-\tau t} Q_\sigma(-it)(\coth(\frac{\pi t}{2}) - \frac{2}{\pi t}) dt.$$

The summand with $\frac{2}{\pi t}$ gives an odd rational function thus does not contribute to the functional equation. The rest is analogous to the odd case. ∎

## 8.2   Higher Rank

The main ingredients of the proof of the functional equation of the theta series in the rank one case were the Euler product and the determinant formula for the zeta function. It is easy to see that the theorem generalizes to all the various situations in the preceding chapters where we had established these two ingredients. We are not going to bore the reader in listing all the results. We will just restrict to the case $SL_2(\mathbb{R})^m$.

We take up the situation of 5.2 again, so $G = SL_2(\mathbb{R})^m$ for some natural number $m$, $\Gamma \subset G$ is a nice cocompact discrete subgroup and $(\varphi, V_\varphi)$ is a finite dimensional unitary representation of $\Gamma$. The representation $\varphi$ defines a holomorphic flat hermitian vector bundle $E_\varphi$ over $X_\Gamma = \Gamma \backslash X = \Gamma \backslash G / K$, We denote by $\triangle_{p,q,\varphi}$ the Hodge-Laplace operator on $E_\varphi$-valued $(p, q)$-forms. We consider the **theta function** :

$$\Theta_\varphi(\tau) = \sum_{q=0}^{m} \operatorname{tr} e^{-\tau \sqrt{\triangle_{0,q,\varphi} - \frac{1}{4}}}, \quad \operatorname{Re}(\tau) > 0.$$

To establish a functional equation of this theta series we also need a **dual theta series**:

$$\begin{aligned}
\Theta_d(\tau) &= \operatorname{tr} e^{-\tau P} \\
&= \sum_{j=0}^{\infty} (2j+1) e^{-\tau(j + \frac{1}{2})}.
\end{aligned}$$

Then we get

**Theorem 8.2.1** *The dual theta function* $\Theta_d(\tau)$ *extends to an even meromorphic function. The poles of* $\Theta_d$ *are of order 2 and are located at* $\tau = \pi i k$, $k \in \mathbb{Z}$. *The theta function* $\Theta_\varphi(\tau)$ *extends to a meromorphic function on the entire plane which satisfies the functional equation*

$$\Theta_\varphi(\tau) + \theta_\varphi(\tau) = 2m \dim(\varphi) \operatorname{Ar}(X_\Gamma) \Theta_d(i\tau).$$

*The poles of* $\Theta_\varphi(\tau)$ *are the nonnegative poles of* $\Theta_d(i\tau)$ *and simple poles at* $\tau = \pm i l_\gamma$ $\gamma \in \Gamma - \{1\}$, $\mathbb{R} - \operatorname{rank}(\gamma) = 1$ *and the residue at* $\tau = \pm i l_\gamma$ *is*

$$-\frac{1}{2\pi} \frac{l_\gamma}{\mu_\gamma} \frac{\operatorname{tr}(\varphi(\gamma))}{\det(1 - \gamma^{-1}|\mathbf{n})} \chi_1(X_\gamma).$$

∎

The remaining cases of higher rank are left as an exercise to the interested reader.

# Index